\documentclass[10pt,twoside]{amsart}

\usepackage{amssymb}
\usepackage{amsxtra}
\usepackage{graphicx}

\theoremstyle{plain}

\begingroup\makeatletter\ifx\SetFigFont\undefined
\def\x#1#2#3#4#5#6#7\relax{\def\x{#1#2#3#4#5#6}}%
\expandafter\x\fmtname xxxxxx\relax \def\y{splain}%
\ifx\x\y   
\gdef\SetFigFont#1#2#3{%
  \ifnum #1<17\tiny\else \ifnum #1<20\small\else
  \ifnum #1<24\normalsize\else \ifnum #1<29\large\else
  \ifnum #1<34\Large\else \ifnum #1<41\LARGE\else
     \huge\fi\fi\fi\fi\fi\fi
  \csname #3\endcsname}%
\else
\gdef\SetFigFont#1#2#3{\begingroup
  \count@#1\relax \ifnum 25<\count@\count@25\fi
  \def\x{\endgroup\@setsize\SetFigFont{#2pt}}%
  \expandafter\x
    \csname \romannumeral\the\count@ pt\expandafter\endcsname
    \csname @\romannumeral\the\count@ pt\endcsname
  \csname #3\endcsname}%
\fi
\fi\endgroup

\newtheorem{thm}{Theorem}
\newtheorem{theorem}{Theorem}[section]

\newtheorem{defn}[theorem]{Definition}

\newtheorem{remark}[theorem]{Remark}

\newtheorem{lemma}{Lemma}
\newtheorem{cor}{Corollary}

\newcommand{\la}{\lambda}
\newcommand{\al}{\alpha}
\newcommand{\be}{\beta}
\newcommand{\ga}{\gamma}
\newcommand{\de}{\delta}

\newcommand{\thn}{\theta_{k+1}}

\newcommand{\Q}{\mathbb{Q}}
\newcommand{\R}{\mathbb{R}}
\newcommand{\T}{\mathbb{T}}

\newcommand{\Z}{\mathbb{Z}}

\newcommand{\Zmod}[1]{\Z/{#1}\Z}

\newcommand\ra{\rightarrow}
\newcommand\lra{\longrightarrow}

\newcommand\seta{\{\alpha_i\}_{i=1}^g}
\newcommand\setb{\{\beta_i\}_{i=1}^g}
\newcommand\setg{\{\gamma_i\}_{i=1}^g}

\newcommand\setd{\{\delta_i\}_{i=1}^g}

\newcommand{\Si}{\Sigma}

\newcommand{\ben}{\begin{enumerate}}
\newcommand{\een}{\end{enumerate}}

\newcommand{\Oz}{P.\ Ozsv{\'a}th\,}
\newcommand{\Sz}{Z.\ Szab{\'o}\,}

\begin{document}
\title{Heegaard-Floer Homology and String Links}
\author{Lawrence P. Roberts}
\maketitle

\section{Introduction}

In \cite{Knot} \Oz and \Sz use the technology of Heegaard-Floer homology 
to refine the Alexander-Conway polynomial of a marked knot in $S^{3}$. 
In particular, they define Knot-Floer homology groups for relative 
$Spin^{c}$ structures that correspond to 
the 
terms in the polynomial: the Euler characteristic of the homology
with 
rational coefficients corresponding to $i$ gives the coefficient of 
$t^{i}$. 
They have since shown that the 
non-vanishing of these groups characterizes the genus of the knot, 
\cite{Genu}. In 
\cite{Alte} 
they employ Kauffman's state summation approach to the Alexander-Conway 
polynomial 
to give a concrete realization of the generators of the chain complex for 
each $i$ and 
their gradings (though not, unfortunately, for the differential). 
Furthermore,
the techniques extend to {\em null-homologous} knots in an arbitrary three 
manifold, $Y$,
where the knot may also be interpreted as giving a filtration of the 
Heegaard-Floer
chain groups for $Y$ that is also an invariant of the isotopy class of the 
knot. 
They also refine the one variable Alexander-Conway polynomial of an 
$m$-component
link in $S^{3}$ by converting the link, 
in a way preserving isotopy classes, to a null-homologous 
knot in a $\stackrel{m-1}{\#} S^{1} \times S^{2}$.  

In this paper we simultaneously generalize the preceding picture in two 
ways: first, by removing
the restriction on the homology class of the embedding, and second, by 
defining chain complexes for a ``string link'' thereby obtaining the Heegaard-Floer 
analog of the {\em multi-variable} torsion of the string link for its 
universal Abelian covering space. There is a relationship between the torsion
of a string link and the multi-variable Alexander polynomial of a simple
link closure of the string link. Recently, \Oz and \Sz have announced a 
version of Heegaard-Floer homology for links which enhances the multi-variable
Alexander polynomial as the knot Floer homology enhanced the single-variable version.
While there should be a relationship between the string link homology and
the link homology, it should be noted that the string link homology appears
to be a different beast; it does not share the symmetry under change of the
components' orientation, for example. On the other hand, we will show that in most 
respects the string link homology is a natural generalization of the knot Floer homology.

In general, we start with a basing of a link in a three manifold, $Y$: 
we require an oriented disc $D$ embedded in $Y$ so that the link 
components 
intersect the disc once, positively. This configuration is known as a 
$d$-base, 
following the work of N. Habegger and X.S. Lin, \cite{Habb}. There are 
many $d$-basings
of the same link, and our invariant will be sensitive to these. In $S^{3}$ there is a more perspicuous description of the
configuration, called a string link:

\begin{defn} Choose $k$ points $p_{1}, \ldots, p_{k}$ in  $D^{2}$. A 
$k$-stranded ``string link'' in $D^{2} \times I$ is a proper embedding, 
$\coprod_{i=1}^{k} f_{i}$ of $\coprod_{i = 1}^{k} I_{i}$ into $D^{2} 
\times I$, where $f_{i}: I_{i} \ra D^{2} \times I$, such that  $f_{i}(0) =  
p_{j} \times 0 $ and $f_{i}(1) = p_{s} \times 1 $. The string link 
is called ``pure'' if $j = s$ for each interval. 
\end{defn}

\noindent A neighborhood of a $d$-base, $D$, is a copy of $D^{2} \times I$, and its 
complement 
in $S^{3}$ is also a copy of $D^{2} \times I$. In this $D^{2} \times I$, 
the 
$d$-based link appears as $k$ copies of $I$ extending from one end to the 
other.
For more general three manifolds, we will present our string link as a 
string link
in $D^{2} \times I$ with the data of a framed link diagram where all the 
surgeries
occur in $D^{2} \times I$. \\
\ \\
For a $d$-based link $\mathbb{L} \cup D$ embedded in a general 
three manifold $Y$, we measure the homology of the
components of the embedding by a lattice, $\Lambda$, in $\Z^{k}$ spanned 
by the all vectors of the form $([h] \cap L_{1}, \ldots, [h] \cap L_{k})$, 
where  $[h] \in H_{2}(Y;\Z)$ and the $L_{i}$ are the components of the 
associated link. We will prove:

\begin{thm}
Let $\mathbb{L} \cup D$ be a $d$-based link in $Y$. 
Then, for each $Spin^{c}$ structure, $\mathfrak{s}$, on $Y$, there 
 is a relatively $\Z^{k}/\Lambda$-indexed abelian group 
$\oplus_{\lambda} \widehat{HF}(Y, \Gamma; \mathfrak{s}, \lambda)$ 
where each of the factors is an isotopy invariant of the based link.
\end{thm}
When $\Lambda \equiv 0$, most of the results for knots transfer 
straightforwardly. In particular, the presence of the $d$-based link
imposes a filtration upon the Heegaard-Floer chain complexes of the 
ambient manifold that is  filtered chain homotopy invariant 
up to isotopy of based link. The case where $\Lambda \not\equiv 0$ occurs 
naturally 
when trying to define chain maps from cobordisms of knots in cobordisms of
three manifolds. 
\\
\\
{\bf Note:} In \cite{Knot} the Knot-Floer homology is denoted by 
$\widehat{HFK}(Y, K; \cdot)$. We will assume that the presence of a 
$\Gamma = \mathbb{L} \cup D$ or $K$ implies the use of the data determined 
by that object. Thus, we will use $\widehat{HF}(Y, K; \cdot)$ for the 
knot-Floer homology. When we wish to refer to the Heegaard-Floer 
homology of the ambient three manifold (ignoring the information provided 
by $K$) we will simply omit reference to the knot or bouquet. 
However, these are not relative homology groups according to the
classical axioms, nor do they solely depend upon the complement of 
the knot or link. 
\\
\\
When $\mathbb{L} \cup D \subset S^{3}$, we may re-state the theorem in 
terms of the associated string link, $S$:

\begin{cor}
For each $\overline{v} \in \Z^{k}$ there is an isotopy invariant 
$\widehat{HF}(S; \overline{v})$ of the string link $S$.
\end{cor}

In section \ref{sec:gencomb}, we realize the generators of the chain 
complex for a string link $S$ in $D^{2} \times I$ from a projection of 
$S$. 
They are identified with a sub-set of maximal forests -- 
satisfying specific constraints imposed by the meridians -- in 
a  planar graph constructed from the projection of $S$. This description 
generalizes the 
description of generators, their indices, and their gradings given by \Oz 
an \Sz in 
\cite{Alte}. 

\begin{lemma}
There are vector weights assigned to crossings so that for 
each tree, adding the weights calculates in which index, $\overline{v}$, 
the corresponding generator occurs. Furthermore, there are weights
assigned to crossings which likewise calculate the grading of the 
generator.
\end{lemma}

For the specific weights see Figure \ref{fig:weights}. 
This lemma requires a generalization of L. 
Kauffman's ``Clock Lemma'' to maximal forests that describes the 
connectivity of the set of 
maximal forests in a planar graph under two natural operations. 
\\
\\
The weights and gradings are enough to form the Euler characteristic of 
the homology 
groups with rational coefficients, 
which 
is related to the Alexander-Conway polynomials of 
the link components in section \ref{sec:AlexCon}. However, the Euler 
characteristic can 
also be interpreted
as a polynomial arising from the first homology of a covering space. 

We let $X = D^{2} \times I - \mathrm{int}(N(S))$ and $E = \partial X - 
D^{2} \times \{0\}$. Consider the
$\Z^{k}$-covering space, $\widetilde{X}$, determined by the Hurewicz map 
$\pi_{1}(X) \ra H_{1}(X;\Z)$ $\cong \Z^{k}$ taking meridians to 
generators. Let $\widetilde{E}$ be the pre-image of $E$ under the 
covering map. Then there is a presentation, \cite{Thet}, $\Z^{k} \ra 
\Z^{k} \ra H_{1}(\widetilde{X}, \widetilde{E}; \Z) \ra 0$ of $\Z[t^{\pm 1}_{1}, 
\ldots, t^{\pm 1}_{k}]$-modules whose $0^{th}$ elementary ideal is 
generated
by a $\mathrm{det}\, M$, a polynomial called the {\em torsion} of the 
string link.

\begin{thm}
Let $S$ be a string link. The Euler characteristic of $\widehat{HF}(S; 
\overline{v}; \Q)$ is the coefficient of $t_{1}^{v_{1}}\cdots 
t_{k}^{v_{k}}$ 
in a polynomial $p(t_{1}, \cdots, t_{k})$ describing the torsion of 
the string link, $\tau(S)$, \cite{Livi}.
\end{thm}
R. Litherland appears to have originated the study of the module 
$H_{1}(\widetilde{X}, \widetilde{E}; \Z)$ as a source of Alexander 
polynomials, \cite{Thet}. He used it to study generalized $\theta$-graphs, 
which, 
once we pick a preferred edge, correspond to the string links above.

Many results follow from trying to replicate known properties of the 
torsion. 
Braids are a special sub-class of string links, for which it is known that 
the torsion is always trivial. Likewise, it follows easily that

\begin{lemma}
If the string link $S$ is isotopic to a braid, then $\widehat{HF}(S) \cong 
\Z_{(0)}$ 
where $\widehat{HF}(S) \cong \oplus_{\overline{v}} \widehat{HF}(S; 
\overline{v})$
\end{lemma}
\noindent
The subscript in $\Z_{(0)}$ designates the grading.

This result should be likened to the analogous result for $1$-stranded 
string links, i.e. marked knots, that are also braids: the unknot has 
trivial knot Floer homology.
While string links are usually considered up to isotopy fixing their 
endpoints
on $D^{2} \times \{0, 1\}$, this result has the implication
that our invariant will be unchanged if we move the ends of the strings 
on $D^{2} \times \{0, 1\}$ (but not between ends).

Furthermore, as in \cite{Alte}, 
\begin{lemma}
Alternating string links have trivial differential
in each index, $\overline{v}$, for the Heegaard decomposition arising
from an alternating projection. 
\end{lemma}
\noindent The proof may be found in section \ref{sec:Altern} 

As \Oz and \Sz can extend the knot Floer homology to links, we may extend 
the constructions for
string links to a sub-class of colored tangles in $D^{2} \times I$. For a 
tangle, $T$,
we allow closed components, in addition to open components requiring that
the open components independently form 
a string link. To each open or closed component we assign a color $\{1, 
\ldots, k\}$ which corresponds to the variable $t_{i}$ 
used for that component. We require that each color be applied to one and 
only
one open component. We may then use the colors to construct a string 
link, $S(T)$, in a second manifold, $\stackrel{n}{\#} S^{1} \times S^{2}$, 
where $n$ is the number of closed components in $T$. The isotopy class of
$S(T)$ is determined by that of $T$ in $D^{2} \times I$, 
allowing us to consider $\widehat{HF}(Y, S(T); \mathfrak{s}_{0}, i)$
as an isotopy invariant of $T$. With this definition, we may extend the
skein exact sequence of \cite{Knot} to crossings where each strand has the
same color, see section \ref{sec:gencomb}.
  
Finally, we analyze how the homology changes for three types of string 
link compositions. Each has the
form of connect sum in its Heegaard diagram, and the proofs roughly follow 
the approach for 
connect sums taken by 
\Oz and \Sz. We picture our three manifolds as given by surgery on 
framed links in $D^{2} \times I$
with an additional string link, $S$, with $k$ components. Given such 
diagrams for
$S_{1}$ in $Y_{1}$ and $S_{2}$ in $Y_{2}$ we may 1) place
them side by side to create a string link, $S_{1} + S_{2}$ with $k_{1} 
+ k_{2}$ components, 
2) when $k_{1} = k_{2}$ we
may stack one diagram on top of the other (as with composition of braids) 
to obtain the string
link $S_{1} \cdot S_{2}$, 
and 3) we may replace
a tubular neighborhood, i.e. a copy of $D^{2} \times I$, of the $i^{th}$ 
strand in 
$S_{1}$ with the entire picture
for $S_{2}$ to obtain a string satellite, $S_{1}(i, S_{2})$. The analysis 
of the second
type proceeds differently than in \cite{Knot}: we consider it as a closure 
of $S_{1} + S_{2}$
found by joining the ends of $S_{1}$ on $D^{2} \times \{1\}$ with the ends 
of $S_{2}$ on 
$D^{2}\times \{0\}$ in a particular way. We prove the following
formulas for the homologies, where $\mathfrak{s} = \mathfrak{s}_{0} \# 
\mathfrak{s}_{1}$:

$$
\widehat{HF}(Y, S_{1} + S_{2}; \mathfrak{s}, [\overline{j}_{0}] 
\oplus [\overline{j}_{1}]) \cong H_{\ast}(\widehat{CF}(Y_{0}, S_{0}; 
\mathfrak{s}_{0}, [\overline{j}_{0}]) \otimes \widehat{CF}(Y_{1}, 
S_{1}; \mathfrak{s}_{1}, [\overline{j}_{2}]))
$$

$$
\widehat{HF}(Y, S_{1} \cdot S_{2};\mathfrak{s}, [\overline{k}]) = 
\bigoplus_{[\overline{k}_{1}] + [\overline{k}_{2}] = [\overline{k}]\, 
\mathrm{mod }\Lambda} H_{\ast}( \widehat{CF}(Y_{1}, S_{1}; 
\mathfrak{s}_{1}, [\overline{k}_{1}]) \otimes \widehat{CF}(Y_{2}, S_{2}; 
\mathfrak{s}_{2}, [\overline{k}_{2}]))
$$

$$
\begin{array}{c}
\widehat{HF}(Y, S_{1}(i, S_{2}); \mathfrak{s}, [(l_{1}, \ldots, 
l_{k_{1}+k_{2} - 1})] \cong  \\ \ \\ \bigoplus_{[\overline{v}'] + 
[\overline{w}'] = 
[\overline{l}] \mathrm{mod\ }\Lambda'} H_{\ast}(\widehat{CF}(Y_{0}, S_{0}; 
\mathfrak{s}_{0}, [\overline{v}]) \otimes \widehat{CF}(Y_{1}, S_{1}; 
\mathfrak{s}_{1}, [\overline{w}]))
\end{array}
$$
\noindent
where $\overline{v}' = (v_{1}, \ldots, v_{i-1}, v_{i}, \ldots, v_{i}, 
v_{i+1}, \ldots, v_{k_{1}})$, repeating $v_{i}$ a total of $k_{2}$-times,
and $\overline{w}' = (0, \ldots, 0,$ $w_{1}, \ldots, w_{k_{2}}, 0, \ldots, 
0)$ with zero entries except for places $i, \ldots, i+k_{2}-1$.
and $\Lambda' = \Lambda + \overline{0} 
\oplus \Lambda_{1} \oplus \overline{0}$.
\\
\\
{\bf Summary by section:}
\\
\\ \noindent
In section \ref{sec:dbase} we describe Heegaard decompositions for string 
links. In
particular, we analyze the Heegaard equivalences preserving the structure 
of a $d$-based
link and relate them to the types of diagrams -- multi-point Heegaard 
diagrams -- used
to define the chain complexes. 
Section \ref{sec:multi} describes
the indices that replace the exponents of the 
Alexander polynomials in the general setting and defines the chain
complexes we will use. We
define the homology theory and record its basic properties. In section 
\ref{sec:Inv} 
we provide a summary of invariance 
slanted towards
dealing with the effects of an ambient manifold with non-zero second Betti 
number.  Section 
\ref{sec:Chain} 
develops the necessary theory
for chain maps and is here largely for completeness. Section 
\ref{sec:Chainlamb} 
illustrates the first substantive
examination of the effect of $\Lambda$. Section \ref{sec:stringcomp} 
records the proofs of 
the results cited for
combining string links.  In section
\ref{sec:Euler} we review the study of Alexander invariants related to 
string links in $D^{2} \times I$. In \ref{sec:gencomb} we describe
the state summation approach to the Euler characteristic and its
relationship with the generators of the chain complex. We also
describe the results for braids and alternating string links; the 
relationship
with the Alexander-Conway polynomial of the link components; and the skein
exact sequence.
\\
\\
{\bf Note:} We do not address issues of orientation of moduli spaces in 
the paper. However, nothing we say will alter the existence of the 
coherent orientations. It is merely convenient to suppress this 
information. As usual we may work with $\Zmod{2}$-coefficients to 
avoid these issues.

\setcounter{tocdepth}{1}
\tableofcontents

\section{String Links and Heegaard Diagrams}
\label{sec:dbase}

\subsection{String Links in $D^{2} \times I$}

\begin{figure}
\begin{center}
\begin{picture}(5957,3534)
\includegraphics{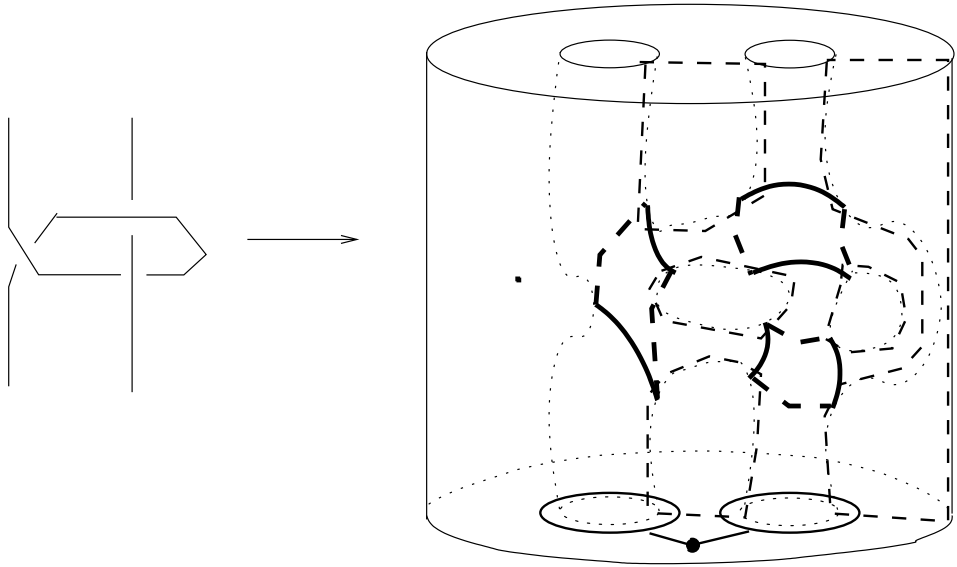}
\end{picture}%
\end{center}
\caption{The Heegaard Diagram for a String Link. We place the meridians at 
the bottom
of the diagram and $\al$ curves at each of the crossings, reflecting the 
type of crossing.}
\label{fig:String}
\end{figure}

In this section, the topological input for our invariant will be laid out. 
To start
we consider the situation in $S^{3}$:

\begin{defn} Choose $k$ points $p_{1}, \ldots, p_{k}$ in  $D^{2}$. A 
$k$-stranded ``string link'' in $D^{2} \times I$ is a proper embedding, 
$\coprod_{i=1}^{k} f_{i}$ of $\coprod_{i = 1}^{k} I_{i}$ into $D^{2} 
\times I$, where $f_{i}: I_{i} \ra D^{2} \times I$, such that  $f_{i}(0) =  
p_{j} \times 0 $ and $f_{i}(1) = p_{s} \times 1 $.  The string link is 
{\textrm pure} if $j = s$ for each interval. We orient 
the strands ``down'' from $1$ to $0$.
\end{defn}
We consider these up to isotopy, preserving the ends $D^{2} \times \{i\}$ 
for $i = 0, 1$. Then braids form a sub-set of string links. However, the 
invariant we define will allow us to move the ends of the string link on 
their end of the cylinder, but not between ends. With this additional 
freedom we can always undo any braid.\\
\ \\
We consider $D^{2} \times I$ embedded inside $S^{3}$. A projection of the 
string link onto a plane provides the data for a Heegaard decomposition of 
$S^{3}$. For a string 
link whose strands are oriented downwards, we may draw a Heegaard diagram 
as in Figure \ref{fig:String}.

We take a small tubular neighborhood of each strand in $D^{2} \times I$ 
and glue it
to the three ball that is $S^{3} - i(\,D^{2} \times I\,)$. This is a
handlebody to which we attach $1$ handles at each of the crossings. These 
handles occur
along the axis of the projection, between the two strands; when we 
compress the handle 
to obtain the tubular neighborhood of an ``X", the attaching circles 
appear as in
Figure \ref{fig:String}, one for each type of crossing, cf. \cite{Alte}. 
The attaching
circles for the handles from the strands are placed in $D^{2} \times 
\{0\}$ as meridians
for each strand. 

The complement of this in $S^{3}$ is also a handlebody. It can be 
described by taking two
$0$-handles in $D^{2} \times I$, above and below the plane we projected 
onto (thought
of as cutting $D^{2} \times I$ through the middle). We then attach handles 
through the
faces of image of the projection: a copy of the $D^{2}$ factor in the 
handle should be
the face. We use $1$-handle for the face on the leftmost side of the 
diagram, called $U$, to cancel one of the $0$-handles. As these two 
handlebodies have the same boundary, they
must have the same number of handles, and the decomposition is a Heegaard 
decomposition.

To each strand of a string link, we may associate a knot: ignore the other 
strands and join
the two ends of our strand with an unknotted arc in the complement of 
$D^{2} \times I$. Likewise we may associate a link to a pure string link 
by using $k$ unknotted, unlinked arcs
as for the closure of a braid. Furthermore, by retaining $D^{2} \times 
\{0\}$, oriented as the boundary of $D^{2} \times I$, we have an embedded 
disc which intersects the strands of the link in one point 
with $L_{i} \cap D = + 1$. This is a $d$-base in the language of Habbegger 
and Lin, \cite{Habb}.

\subsection{String Links in $Y$}

For general three manifolds, $Y$, we define a string link to be a 
$d$-based link:

\begin{defn}
A $d$-base for an oriented link $\mathbb{L}$ is an oriented disc, $D$, 
embedded in 
$Y$ whose interior
intersects each component of $\mathbb{L}$ exactly once, positively.
\end{defn}

By thickening the disc, $D$, we obtain an embedded ball with the structure
of a cylinder, $D^{2} \times I$. $\mathbb{L}$ in the complement of
this ball comports more with our intuition for a string link. We consider 
these complexes up to isotopy of the complex $D \cup \mathbb{L}$, 
preserving the orientations on the parts. However, the invariant will 
allow the strands of  $\mathbb{L}$ not to return to the same point, and 
these ends to move independently on, but not between, the two sides of 
$D$.
\\
\\
{\bf Example:} For knots in $S^{3}$ a $d$-basing amounts to picking
a point on the knot, the intersection point with $D$. This is a marked
knot as used in \cite{Knot}. The string link associated to a knot comes 
from dividing the knot at that point and pulling the ends apart. 
\\
\\
We will use Heegaard diagrams for $Y$ where the $\al$-handlebody extends 
the handlebody that is a neighborhood of $D \cup \mathbb{L}$. If we 
think of this handlebody up to isotopy we would have great freedom in 
moving around the handles corresponding to the strands. When we 
transfer to the Heegaard picture the additional freedom
allows us to braid the ends of the strands through isotopies of the 
Heegaard surface,
without altering the disc. We may even disjoin the link components: if 
the endpoints of a strand correspond, we have a $d$-based link, but we may 
braid 
the ends on one side of the disc independently of the other. However, 
the data from $D$ prevents the ends from moving between the 
different sides of $D$. 
\\
\\
{\bf Example:} To draw the Heegaard diagram for a string link in $Y$, we 
think of $D \cup \mathbb{L}$ embedded in a framed link diagram for $Y$ in 
$S^{3}$. We may
then draw a Heegaard diagram from a projection of the whole configuration, 
attaching the
framed components by paths to $D \cup \mathbb{L}$. 
The $\al$ curves on the new components are no longer the meridians, but 
instead are chosen to be framing curves for the surgery, oriented 
opposite the longitudes from $S^{3}$. All of this can be done
away from the disc $D^{2} \times \{0\}$ containing the meridians. By 
inverting a neighborhood
of the disc, we can present the string link in $Y$ by a string link
in a framed link diagram in $D^{2} \times I$ with $D^{2} \times \{0\}$ as 
preferred disc.

\subsection{Heegaard Diagrams and Based Links}
We describe the general approach to Heegaard diagrams subordinate
to a string link, $D \cup \mathbb{L}$. Above we gave 
specific
examples of such diagrams; we would now like to characterize all such
diagrams and describe the Heegaard equivalences between them. In 
particular,
we will be concerned with retaining the data provided by $D$. 
\\
\\
First we establish some conventions. Let $Y$ be a three manifold. Assume,
as given, 
a Heegaard decomposition, $Y = H_{\al} \cup_{\Sigma^{g}} H_{\be}$, where 
$\partial H_{\al}$ $= \Sigma$ $= -\partial H_{\be}$ and the gradient flow 
corresponds to the 
outward-pointing normal to $H_{\al}$. We denote this decomposition 
by  $(\Sigma, \seta, \setb)$ with $\seta$ being the co-cores of the one 
handles in $H_{\al}$ and likewise for $\setb$. 
\\
\\
To make our diagram reflect a string link $\mathbb{L} \cup D$, we
additionally require that our Heegaard decomposition satisfies:

\ben
\item The three manifold $(\Sigma, \{\alpha_{i}\}_{i=k+1}^{g}, \setb)$ is 
$Y^3 - N(\mathbb{L} \cup D)$. The complement of \newline $(\Sigma, 
\{\alpha_i\}_{i \neq j}, \setb)$ for $j \leq k$, should be
homeomorphic to a tubular neighborhood of  $L_{j}$ in $Y$. We require 
$\al_j$ to be an oriented meridian for this tubular neighborhood.

\item There is a disc $D'$ in $\Sigma-\seta-\setb$ whose boundary contains
one connected segment from each of $\al_{1}, \ldots, \al_{k}$. The string 
link 
formed from  $D'$ glued to the attaching discs
for $\al_{1}, \ldots, \al_{k}$ and the components of $\mathbb{L}$ from (1)
is isotopic to $\mathbb{L} \cup D$. 
\een

The first condition requires that the diagram be subordinate to 
the string link. The second requires that there be a disc in the 
Heegaard
data that produces a $d$-basing for $\mathbb{L}$ and that induces a string
link isotopic to the original one.  Such diagrams 
are said to 
be {\em subordinate} to $\mathbb{L} \cup D$. Furthermore, a diagram with
such a choice of $\al$'s and a disc $D'$ determines a string link as in 
item (2). Note that $D'$ must be oriented opposite to $\Si$.
\\
\\
We may relate diagrams subordinate to a configuration $\mathbb{L} \cup D$, 
embedded inside a three manifold, $Y$, by the following lemma:

\begin{lemma}{\label{lem:equiv}\textrm{[cf. Lemma 4.5 \cite{Four}]}} Let 
$Y$ be a 
closed, oriented three-manifold. Let
$\mathbb{L} \cup D \subset Y$ be an embedded $d$-base.
 Then there is a Heegaard diagram subordinate to 
$\mathbb{L} \cup D$  
and any two such subordinate diagrams may be connected by a sequence of 
the following moves:

\ben
\item Handleslides and isotopies among the elements of 
$\{\al_i\}_{i=k+1}^{g}$
\item Handleslides and isotopies of $\setb$.
\item Stabilization introducing $\al_{g+1}$, $\be_{g+1}$ intersecting in 
a single point.
\item Isotopies of $\{\al_i\}_{i=1}^{k}$ and handleslides of 
	them over elements of $\{\al_i\}_{i=k+1}^{g}$.

\een
where we disallow isotopies and handleslides of any attaching 
circles resulting in a curve intersecting the disc $D'$.
\end{lemma}
 
\noindent {\bf Proof:} Let $\Gamma = \mathbb{L} \cup D$. $N(\Gamma)$ is a handlebody
which we may extend to a Heegaard decomposition for $Y$. The disc $D$
corresponds to a disc $D'$ in $\partial N(\Gamma)$ at one end of the 
thickening of $D$, determined by the orientation requirements. The 
meridians
are chosen to occur in the disc, and thus abut $D'$ in accordance with
requirements. The additional $\al$'s and $\be$'s may be chosen to avoid
$D'$ as $D'$ is contractible in $\Si$. Thus, such a diagram does exist. 
It is shown in  \cite{Hom3}, Prop. 7.1,  that any isotopy across a 
contractible 
region in $\Sigma$ may be obtained by handleslides and isotopies not 
crossing that 
contractible region; any $\be$ isotopy which crosses the entire star may 
be accomplished by handleslides not intersecting the star. This assures us 
that the
choice made in pushing the $\be$'s away from $D'$ does not affect the 
outcome.  

Given the diagrams for two isotopic embeddings of $\Gamma$, we must see 
that they can be related by the moves described above. These moves 
preserve the region $D$ in the original diagram. This region and a small 
neighborhood of it in $Y$ act as the disk for the $d$-base. On the other
hand, the isotopy carries a neighborhood of 
the disc into a neighborhood of the new disc; these neighborhoods are
all homeomorphic and may be used as a $3$-handle for each of the
Heegaard diagrams. That the isotopy preserves $\mathbb{L}$ outside this 
ball allows us to fix the meridians and consider the additional handles,
describing $Y - N(\Gamma)$. 

If we consider $\partial N(\Gamma)$ as $\partial_{+}(Y - 
N(\Gamma))$ 
then the existence of the Heegaard diagram follows from the existence of a 
relative Morse function that is equal
to $1$ on the boundary and  that the isotopy class of 
$\{\al_i\}_{i=1}^{k}$ is
determined by their being meridians of the knots determined by the
strands of $\mathbb{L}$. As usual, we may cancel $1$-handles with 
$0$-handles until there is only one $0$-handle. Similarly we may cancel 
off $3$-handles until
there are none. The relative version of Cerf's theorem states that any two 
such diagrams can 
be linked through the first three moves and the introduction of new index 
$0/1$ cancelling pairs or new index $2/3$ cancelling pairs. 

However, we would like to ensure that the path can be chosen through 
diagrams with only one index $0$ handle and no index $3$ handles. As we 
introduce a new $3$-handle, we also introduce a cancelling $2$-handle. The 
new $2$-handle  will have one end of its co-core on $\partial N(\Gamma)$, 
since 
there are no other $3$-handles. In the diagram for $Y$, the one with the 
prescribed meridians, this $2$-handle has a core that is a homological 
linear combination of the $\al$'s. If we cut $\Sigma$ along $\{\al_{k+1}, 
\ldots, \al_{g}\}$, and cap the new boundaries with discs,  the image of 
the core will be null-homotopic: it will be homotopic to the boundary of 
the $B^2$ at the end of the co-core lying on $\partial N(\Gamma)$.  Since 
it is 
null-homotopic, the core cannot have a non-zero coefficient for a meridian 
for its homology class.  Thus, it is linear combination of $\{\al_{k+1}, 
\ldots, \al_{g}\}$. According to lemma 2.3 of \cite{Hom3}, the core curve 
can be obtained as the image of a $\{\al_{k+1}, \ldots, \al_{g}\}$ under 
handleslides. Thus, any diagram obtained from the diagram after a $2/3$ 
pair is added, could be obtained from the old diagram as handleslides over 
the new core can be given by handleslides over $\al$'s, not using any 
meridians.  We may bypass $2/3$-handle pairs. The same argument applies 
for $0$-handles, as there are fewer restrictions on the $\be$'s.

Allowing the meridians to move, the last equivalence follows from 
considering the 
bouquet after surgering out $\{\al_i\}_{i=k+1}^{g}$. The meridians are 
determined up to isotopy, and as those isotopies cross surgery discs 
there are corresponding handleslides. Likewise, two curves for the same 
meridian, abutting the same disc $D'$ at a specified point, will be 
isotopic. Thus,
any choice of meridians may be moved to the one chosen above. 
This
takes care of all our choices, so any two diagrams subordinate to the same
isotopy class of string link can be related with the moves in the lemma.
$\Diamond$
\\
\\
Without disallowing equivalences that intersect $D'$, the moves above 
preserve the handlebody
neighborhood of $\mathbb{L} \cup D$. This is weaker than preserving the 
isotopy class of $\mathbb{L} \cup D$. However, with the additional data 
provided by $D'$, any isotopy will pull the disc along, preventing 
the ends of strands from twisting at the intersection point with $D$.
\\
\\
{\bf Example:} We return to our method for drawing a Heegaard diagram for
a string link in $Y$. Using the reduced Heegaard equivalences, 
we argue that the construction does not depend upon how the framed 
components are joined to 
$N(\mathbb{L} \cup D)$. The argument comes from \cite{Four}, 
for the same result for a bouquet, with minor alteration. 

Consider two distinct arcs, $s_{1}$ and $s_{2}$, joining a framed 
component to  
$\Gamma$, and form a regular neighborhood of the graph provided by $\Gamma 
\cup s_{1} \cup s_{2}$. We extend this to a diagram for $S^{3}$ (by adding 
handles for crossings, etc).  We draw a subordinate diagram using $s_{1}$ 
to attach the framed component by placing $\al_{k+1}$ as a meridian to 
$s_{2}$. To obtain a diagram subordinate to the second choice of paths we 
erase $\al_{k+1}$ 
and replace it with $\al'_{k+1}$, a meridian to $s_{1}$. 
We surger out all the $\al_{i}$ for $i > k+1$ to obtain a 
genus $k+1$ surface. We wish to see that $\al_{k+1}$ and $\al'_{k+1}$ are 
isotopic; if they are we may move the correponding curves in the original 
diagram through isotopies and handleslides one into the other. The two 
curves, along with some non-meridional $\al$'s, bound a punctured torus 
coming
from the framed component. After surgering, the other boundaries are 
filled, and surgering the framing attaching circle transforms the torus 
into a cylinder. 
Therefore, the $\al_{k+1}$ and $\al'_{k+1}$ now bound a cylinder 
which does not involve the disc $D'$. 
Each time the isotopy of curves determined by the cylinder crosses a disc 
coming 
from the surgered handles, there is a corresponding handleslide in 
the original picture. This provides a sequence of Heegaard equivalences 
that are allowed under the reduced equivalences of $d$-based links.

Sliding a strand in the string link over a framed curve
produces a new string link which, along with the framed 
components, produces a second Heegaard diagram for $Y$ 
equivalent to that from before handlesliding (a simple process, but
lengthy, and not provided here). Handleslides of framed components over 
each other can be effected by a bouquet with a path joining the two 
components, which we have seen is available, and then using the same 
argument as above. Likewise, adding $\pm 1$ framed unlinked, unknots 
(blowing up/down)  can be effected using the reduced Heegaard 
equivalences. Hence, using two different framed link descriptions for $Y$ 
will not change the equivalence of the Heegard diagram.

\subsection{Marked Diagrams}

It is cumbersome to retain the disc, $D'$, in our Heegaard diagrams. 
Furthermore, since the restricted Heegaard equivalences 
eliminate handleslides over the meridians they often impede the 
simplification of the Heegaard diagram. We give 
another interpretation of the topology making our diagrams more tractable.

Using the embedded star, $D'$, we may introduce marked points into 
$\Si$. We choose $w$ to be in the interior of $D'$ and $z_i$ to be 
on the other side of $\al_{i}$ in the region of $\Si - \seta - \setb$
abutting the same segment as $D'$. From a subordinate diagram for the 
$d$-based link we have realized a {\em multi-pointed diagram}. 
The equivalences for multi-pointed diagrams are the standard Heegaard 
equivalences -- with 
no restriction on handleslides -- but with the caveat that no isotopy of 
an attaching circle or marked point may allow one to cross the other. 

Conversely,  for any Heegaard diagram with a choice of marked points $w, 
z_{1}, 
\ldots, z_{k}$ in $\Sigma - \seta - \setb$ we can construct a Heegaard 
diagram subordinate to a string link. First, choose paths from each 
$z_{i}$ to $w$ crossing only $\al$'s. Then take neighborhoods of the 
gradient flow lines, in $H_{\be}$, joining the index $0$ critical point 
and the marked points $\Si$, remove these neighborhoods from $H_{\be}$ and 
add them to $H_{\al}$. The complement in $H_{\be}$ is still a handlebody, 
since we 
have removed the neighborhood of $k+1$ segments. Adding the
neighborhoods to $H_{\al}$ creates a new handlebody $H_{\al'}$. 
The new $\al$'s are meridians of the gradient flow lines, and the new 
$\be$'s 
are loops following the flow line from $w$ to the critical point, then 
to $z_{i}$ and back along the path we chose in $\Si$, crossing only 
$\al$'s. It is 
straightforward to find the region $D'$: it is a disc in the portion of 
$\partial H_{\al'}$
coming from the exchange.  Thus, a multi-point diagram gives us a diagram 
subordinate
to a string link through the preferred disc $D$. Furthermore, if we use 
the preferred disc, $D$, to produce a multi-point diagram, after some 
handleslides of the $\al$'s over the new meridians, we can de-stabilize 
the new $\al$'s and $\be$'s to obtain the original 
multi-point diagram. 

The relationship between the equivalences for multi-point diagrams and 
those subordinate to string links comes from noticing that when going 
from a 
multi-point diagram to a string link diagram, performing an illicit 
isotopy 
over
a marked point $z_{i}$ corresponds to an illicit handleslide 
over a meridian, according to the construction above. In fact, 
were we to surger out the meridians, the point $z_{i}$ would correspond to 
one of the two discs used to replace that meridian (the other would be 
close 
to $w$ inside the region $D$). Isotoping across it would be the same as a 
handleslide across a meridian.

This construction may depend upon the choice of the new $\be$-paths. 
If we surger all the $\be$ attaching circles in the 
multi-point diagram, the chosen paths become a star in $S^{2}$ joining 
$z_{i}$ to $w$ for all $i$. If two such stars  are isotopic in the 
complement of the marked points, the resulting diagrams for the $\thn$ 
graph are equivalent. Each time the isotopy of a segment crosses a disc 
introduced by the surgery, we should think of our new  $\be$ being slid 
over an old $\be$. Braiding of the marked points in $S^{2}$, carrying 
along the star,
will not, in general, produce a diagram isotopic to the one with the star 
before braiding. However, these produce equivalent diagrams as they are 
both 
stabilizations of the same diagram. This discrepancy, once again, reflects 
the
inability to detect braiding once we switch to considering Heegaard 
diagrams.
\\
\\
{\bf Note:} We may perform handleslides to ensure that in our Heegaard 
diagrams the meridians each intersect only one $\be$-curve. 
When the meridian intersects only one $\be$, the 
boundary of such a region must include multiples of the full meridian
When there is more than one $\be$ intersecting the meridian, we will 
instead count the intersections with the the embedded $\thn$-graph found
from the gradient flow lines through the marked points.

\subsection{Admissibility and $d$-Based Links}
\ \\
In Heegaard-Floer homology, when $H_{2}(Y ; \Z) \not \cong 0$, we must use
diagrams submitting to certain admissibility requirements, \cite{Hom3}. 
We argue here that presenting a $d$-based link as a string link in $D^{2} 
\times I$, with an additional framed link defining $Y$, can be made admissible
without disrupting the disc $D$ or the structure of $D^{2} \times I$.

We make use of the lemmas in section 5 of \cite{Hom3}.

\begin{lemma}\label{lem:tnorm} Let $\mathbb{L} \cup D \subset Y$ be a 
$d$-based 
link. Let $\mathfrak{s}$ be a $Spin^{c}$ structure on $Y$. Then there 
is a strongly/weakly admissible diagram for $(Y, \mathfrak{s}, \mathbb{L} 
\cup D)$ 
presented as a framed link in the complement of a string link in $D^{2} 
\times I$. 
\end{lemma}

\noindent {\bf Proof:} Suppose $Y$ is presented as surgery on a link in $S^{3}$, and 
$\Gamma$ is $\mathbb{L} \cup D$ in this diagram. We adjust this, as above, 
to
be a framed link diagram in $D^{2} \times I$ with a string link. Recall
that we must join the framed components to $\mathbb{L} \cup D$ by paths 
which
we assume do not touch $D^{2} \times \{0\}$. With the 
framing curves as $\al$-curves, this provides a diagram for $Y$. However, 
it need not be admissible; we may need to wind the attaching circles to 
make it so. We must ensure that the winding paths do not affect the
discs $D^{2} \times \{0\}$ or $D^{2} \times \{1\}$. We make two 
observations. 
\ben
\item First, any doubly periodic domain must have at least one boundary 
containing multiples of a framing curve. Otherwise, by replacing framing 
curves with the meridians of the link components, we would obtain a 
periodic 
domain in a diagram for $S^{3}$. Furthermore, two periodic domains may not 
produce the same linear combination of framing curves in their boundaries. 

\item Second, a meridian of a framed component may be chosen to 
intersect the framing curve which replaces it, once and only once, and 
intersect no other $\al$'s. Each will, however, intersect at least 
one $\be$ curve in the projection. By Proposition 5.3.11 of  \cite{Kirb}, 
these meridians generate all of $H_{1}(Y; \Z)$. By winding along them we 
may obtain intersection points representing any $Spin^{c}$ structure; we 
have an intersection point which employs the framing curves intersecting 
the same $\be$'s as the meridian.  
\een
These are the conditions necessary to draw the conclusion of lemmas 
5.2, 5.4, and 5.6 of \cite{Hom3}. These lemmas guarantee the results in 
the proposition.
$\Diamond$
\\
{\bf Remark:}
There are other embedded
objects lurking in the background of our Heegaard diagrams. The first:

\begin{defn}\cite{Four} An oriented bouquet, $\Gamma$, is a one-complex 
embedded in 
$Y$ which is the union of a oriented link $\mathbb{L} =$ $\cup_{i=1}^{k} 
L_i$ with a collection of $k$ embedded segments, $\ga_{i}$, each 
connecting a point on 
$L_i$ to a fixed reference point in $Y$, and otherwise disjoint from  
$\mathbb{L}$ and each other.
\end{defn}
\noindent  
We will consider such objects up to isotopy in $Y$, preserving the graph 
structure. They are also known as clover links, \cite{Levi}. Such
embedded graphs underly the construction of maps in Heegaard-Floer 
homology from four-dimensional cobordisms formed by $2$-handle additions.

Given a $d$-based link, we may form such an object by choosing
our reference point in $D$ and joining it to the intersections with
$\mathbb{L}$ by a tree. There are many such trees, found, for instance, by
braiding in $D$, and these are not necessarily isotopic. However,
our invariant does not distinguish them. Furthermore, two different
string links can give the same isotopy class of bouquet corresponding
to twisting the bouquet along the reference paths (or moving the
ends of the string link back and forth between the ends
of the cylinder). We can record this twisting by adorning a 
bouquet with a $d$-base formed from a disc neighborhood of its reference 
paths. 
\\
\\
Additionally, in an associated multi-pointed diagram, the 
critical points corresponding to  the $3$-cell and the $0$-cell for a 
Morse function compatible with the 
Heegaard decomposition along with the gradient flow lines though the  
marked points gives an embedded copy of a $\thn$-graph with a preferred
edge arising from $w$: a graph with two vertices and $k+1$ edges oriented 
from one 
vertex to the other. By taking a tubular neighborhood of the preferred 
edge we have a copy of $D^{2} \times I$; in its 
complement is a string link determined by the other edges. The region 
around $v_{-}$ together 
with discs bounded by the meridians of the non-preferred edges forms the 
$d$-base. The $\thn$ graph may be obtained from a diagram subordinate to a 
string link by stretching $\Sigma$ along the boundary of $D' \cup_{i \leq 
k} D(\al_{i})$ where $D(\al_{i})$
is the attaching disc for $\al_{i}$. The result is a diagram for a tubular 
neighborhood of the $\thn$-graph with one edge linked by the embedded 
circle. The region 
bounded by the meridians and the new curve, containing $D'$, 
corresponds to $v_{-}$, the vertex into which the edges point. Isotopies 
of $d$-based links correspond to isotopies of the associated 
$\thn$-graphs.

\section{Multi-point Heegaard-Floer Homology}
\label{sec:multi}

\subsection{Background and Notation}
The reader should consult \cite{Hom3} for the notation used below. \\
\ \\
Let $H_{\al} \cup_{\Si^{g}} H_{\be}$ and a choice of $w, z_{1}, \cdots, z_{k} \in \Sigma$
be a multi-pointed strongly/weakly admissible Heegaard decomposition of $Y$. 
We will denote the additional marked point data by $\Gamma$.
Let $\seta$ be a set
of $g$ disjoint, simple, closed curves in $\Si$ whose images are linearly 
independent in $H_{1}(\Sigma; \Z)$ and which bound compression discs in 
$H_{\al}$. Define $\setb$ similarly for $H_{\be}$. We assume
the curves in  $\seta$ and $\setb$ are in general position. 
Then choose a path of generic nearly 
symmetric almost-complex
structures, $J_{s}$, on $Sym^{g}(\Si)$, in accordance with the 
restrictions in 
\cite{Hom3}. Furthermore, choose an equivalence class of intersection 
points, 
$\mathfrak{s}$, for 
$Y$ and 
a coherent system of orientations for the equivalence class, 
\cite{Hom3}\footnote{Those
willing to work with $\Zmod{2}$-coefficients may ignore this requirement}. \\
\ \\
Let $\mathcal{I}(\mathfrak{s})$ the set of all intersection points, ${\bf x} \in \T_{\al} \cap \T_{\be}$
which represent the $Spin^{c}$ structure, $\mathfrak{s}$. For $\mathfrak{s}$, define 
$CF^{\infty}_{\Gamma}(Y; \mathfrak{s})$  as the $\Z$-module  

$$ 
{\textrm Span}_{\Z}\{[{\bf x}, i, \overline{v}] | {\bf 
x}\in\mathcal{I}(\mathfrak{s}),\ i \in\Z,\ \overline{v} \in \Z^{k}\}
$$
\\
There is a natural map on $CF^{\infty}_{\Gamma}$:
$$
U ([{\bf x}, i, \overline{v}]) = [{\bf x}, i - 1, v_{1} - 1, \ldots, v_{k} 
- 1]
$$
which makes $CF^{\infty}_{\Gamma}$ into a module over $\Z [U]$. 
\\
\\
As in \cite{Hom3}, we may define other groups 
by taking $CF^{-}_{\Gamma}(Y ; 
\mathfrak{s})$ to be the sub-group of $CF^{\infty}_{\Gamma}$ where $i < 
0$; by taking $CF^{+}_{\Gamma}(Y; \mathfrak{s})$ to be the resulting 
quotient group (with $i \geq 0$);  and by taking  
$\widehat{CF}_{\Gamma}(Y, \mathfrak{s})$ to be that sub-group spanned by 
those generators with $i=0$. 

The formulas for relative gradings follow as in \cite{Hom3}, \cite{3Man}. In particular, when
$c_{1}(\mathfrak{s}_{w}({\bf x}))$ is not torsion there is a
relative $\Z/\!\delta(\mathfrak{s})\Z$-grading on the chain complexes, where
$$
\delta(\mathfrak{s}) = \gcd_{\xi \in H_{2}(Y; \Z)}< c_{1}(\mathfrak{s}), 
\xi> 
$$
The map $U$ reduces this grading by $2$, \cite{Hom3}. If $\mathfrak{s}$ is a torsion $Spin^{c}$ structure,
the chain complexes inherit the  absolute 
$\Q$-grading, $gr_{\Q}$ defined in 
\cite{Four}. We may present $Y$ as surgery on a link 
which includes the components of $\mathbb{L}$, which receive an
$\infty$-framing. Then we can assign the absolute grading to the 
generators using the formula in \cite{Four}.

\subsection{Filtration Indices}
The additional marked points 
$z_{1}, \ldots, z_{k}$ provide indices for the generators related to the 
grading. 

\subsubsection{Complete Sets of Paths}
\ \\

For the equivalence class $\mathfrak{s}$, we choose a complete set of 
paths for 
$\mathfrak{s}$ as in Definition 3.12 of \cite{Hom3}:

\ben
\item An enumeration $\mathcal{S} = \{ {\bf x_0}, {\bf x_1}, \ldots, {\bf 
x_m}\}$ of the intersection points in   $\mathcal{I}(\mathfrak{s})$.
\item A collection of homotopy classes $\phi_i \in \pi_2({\bf x_0},{\bf 
x_i})$ with $n_w(\phi_i)=0$
\item Periodic domains $\Xi_1, \ldots, \Xi_{b_2}$ $\in \pi_2({\bf x_0}, 
{\bf x_0})$ representing a basis for  $H_2(Y; \Z)$.
\een
\noindent
Any path in $\pi_2({\bf x_i},\, {\bf x_j})$ can then be written uniquely 
as splicings of the $\Xi_i$ and the paths $\phi_i$, and any periodic 
domain in $\pi_{2}({\bf 
x}, {\bf x})$ can be identified with one in $\pi_{2}({\bf x_{0}}, {\bf 
x_{0}})$.
\\
\\
As in \cite{3Man}, given a complete set of paths we may find a map
$$
\mathcal{A} : \pi_2({\bf x}, {\bf y}) \ra H_2(Y;\mathbb{Z})
$$
by taking $\phi_y^{-1} \ast \phi \ast \phi_x \in \pi_2({\bf x_0}, {\bf 
x_0})$ $\cong \mathbb{Z} \oplus H_2(Y; \mathbb{Z})$, and projecting to the 
second factor. This map has the property that $\mathcal{A}(\phi_{1} \ast 
\phi_{2}) =$ 
$\mathcal{A}(\phi_{1}) + \mathcal{A}(\phi_{2})$. The 
procedure provides an identification
of $\pi_{2}({\bf x}, {\bf x})$ with $\pi_{2}({\bf x_{0}}, {\bf x_{0}})$, 
and the action of  $\mathcal{P} \in \pi_{2}({\bf x}, {\bf x})$ on $\phi 
\in \pi_{2}({\bf x}, {\bf y})$ produces $\phi'$ with $\mathcal{A}(\phi') = 
$ $\mathcal{A}(\phi) + \mathcal{A}(\mathcal{P})$.  
\\
\\
 Furthermore, a choice of  basepoint ${\bf x_0}$, a basis for $H_{2}(Y, 
\Z)$, and a choice of an additive $\mathcal{A}$ which maps 
$\pi_2({\bf x_0},{\bf x_0})$ surjectively onto $H_{2}(Y; \Z)$ and 
is invariant under the action of $[S]$, gives a complete set of paths. 
For a doubly periodic 
domain $\mathcal{P} \in \pi_2({\bf x_0},{\bf x_0})$, $\mathcal{A}$ assigns 
it a value in $H_2(Y; \Z)$. For $\phi \in \pi_2({\bf x_0},{\bf y})$ we 
have 
two quantities $\mathcal{A}(\phi)$ and $n_w(\phi)$. By subtracting the 
periodic domain $\mathcal{P}$ with the same value under $\mathcal{A}$ as 
$\phi$ and subtracting the right number of $[S]$'s, we find a $\phi'$ 
where both quantities are zero. This $\phi'$ is unique since $\mathcal{A}$ 
is surjective, so we choose it as our element, $\phi_y$, in a
complete set of paths.

\subsubsection{Definition of Filtration Indices}

We may define a filtration index for the complete set of paths as a map
$\overline{\mathcal{F}} : \mathcal{I}(\mathfrak{s}) \ra \Z^{k}$ satisfying 
the basic relation

$$
\overline{\mathcal{F}}({\bf y}) - \overline{\mathcal{F}}({\bf x}) = (n_w - 
n_{\overline{z}})(\phi) + n_{\overline{z}}(\mathcal{A}(\phi))
$$
where $n_{\overline{z}}(\phi) = (n_{z_1}(\phi), \ldots, n_{z_k}(\phi))$ 
and $n_{w}(\phi) = (n_{w}(\phi), \ldots, n_{w}(\phi))$. 
These compare 
the information found at the preferred point, $w$, to that at any other 
marked 
point $z_{k}$.
 
When we add a periodic domain $\mathcal{P}$ to $\phi$ the right hand side 
of the filtration relation changes
by $-n_{\overline{z}}(\mathcal{P}) + n_{\overline{z}}(\mathcal{P}) = 0$. 
Thus, with the complete set of paths, the  relation determines the value 
of 
$\overline{\mathcal{F}}$ on every intersection point, independent of 
$\phi$, up to a vector: any other filtration index 
for ${\bf x_{0}}$ and $\mathcal{A}$ has the form $\overline{\mathcal{F}}  
+  v$
for some vector $v \in \Z^{k}$. 

The last term on the right can be re-written in the $i^{th}$ coordinate as 
$L_{i} 
\cap \mathcal{A}(\phi)$. The 
orientation on the meridians is the one induced from the attaching disk 
oriented to intersect $L_{i}$ positively: $L_{i} \cap D_{i} = + 1$. If we 
choose $\mathcal{P}$ a periodic region representing the homology class $h 
\in H_{2}$ then  $\partial \mathcal{P}$ may contain multiples of the 
meridians. By drawing the diagrams we find that 
with this orientation convention $(n_{w} - n_{z_{i}})(\mathcal{P}) = $ $- 
L_{i} \cap h$. The quantity $n_{w} - n_{z_{i}}$ measures\footnote{When $\al_{i}$ 
intersects more than one $\be$, and with the 
gradient flowing with the outward 
normal of $H_{\al}$,  $n_{w} - n_{z_{i}}$ is  minus the
intersection  number with the edge corresponding to $L_{i}$ in the graph, 
$\Theta(\Gamma)$.}
the number of 
times the $i^{th}$ meridian occurs in  $\partial 
\mathcal{P}$.

Adding or 
subtracting elements of $\pi_{2}({\bf x}, {\bf x})$ thus alters the last 
term 
by vectors in  the lattice, $\Lambda$ $\subset \Z^{k}$, spanned by: 

$$
\begin{array}{c}
(n_{z_1}(\Xi_1),\, \ldots,\, n_{z_k}(\Xi_1)) \\
\vdots \\
(n_{z_1}(\Xi_k),\, \ldots,\, n_{z_k}(\Xi_k))
\end{array}
$$
since we assume $n_w(\Xi_j) = 0$. 
\\
\\
Differences in the filtration index can be calculated directly if we know 
explicitly the homotopy classes in the complete set of  paths. The 
difference 
between ${\bf x}$ and ${\bf y}$ is given by $n_{\overline{z}}(\phi_y \ast 
\phi_x^{-1})$. 
For this composite disc $\mathcal{A} = 0$ and $n_{w} = 0$. For rational 
homology 
spheres the choice of $\phi \in \pi_{2}({\bf x}, {\bf y})$  varies only 
by multiples of $[S]$. Thus,   $\phi_y \ast -\phi_x =$ $\phi + r[S]$. In 
this case, or in any case where $\Lambda \equiv 0$, the filtration index
takes values in $\Z^{k}$, and we recover the 
formulas used for knots in \cite{Knot}. 

Once we have fixed the value of ${\bf x_{0}}$, the choice of $\mathcal{A}$ 
prescribes a 
value in $\Z^{k}$ for each ${\bf y}$. Different choices of complete sets 
of paths prescribe different values for ${\bf y}$; however, all these 
values map to the same element of $\mathbb{Z}^k/\Lambda$. For example, 
consider ${\bf x}$ and ${\bf y}$ joined by a path 
$\phi$ with $n_w(\phi)=0$. If we
change to $\mathcal{A}'$ then $-n_z(\phi -\mathcal{A}'(\phi)) = 
-n_z(\phi'_{\bf y}) + n_z(\phi'_{\bf x})$ $= -n_z(\phi_{\bf y} \ast 
-\phi_{\bf x}\ast \mathcal{P}) =$ $-n_z(\phi - \mathcal{A}(\phi)) + 
\lambda $. Changing $\mathcal{A}$ to $\mathcal{A}'$ changes the relation 
by an element of $\Lambda$ for each pair of intersection points. 
\\
\\
Thus, we may remove the dependence upon $\mathcal{A}$ by considering
filtraion indices with values in the quotient  $\Z^{k}/\Lambda$. 
The intersection points are then relatively 
$\Z^{k}/\Lambda$-indexed. The $\Z^{k}$ indices are  
``lifts'' of these indices, which we use when needing to facilitate 
comparisons as $\Lambda$ changes. 

\begin{remark} From now on, we  assume, as chosen, a point ${\bf 
x_{0}} \in \mathcal{I}(\mathfrak{s})$ and a complete set of paths for 
$\mathfrak{s}$ and ${\bf x_{0}}$. Furthermore, we require that
if $z_{i}$ and $z_{j}$ are in the same component of  $\Sigma - \seta - 
\setb$ then $\mathcal{F}_{i}$ and $\mathcal{F}_{j}$ must be 
equal. If $z_{i}$ is in the same component as $w$ then  $\mathcal{F}_{i} 
\equiv  C_{i}$, 
a constant, which we require to be $0$, unless otherwise noted. 
\end{remark}

\subsubsection{A Special Case}
\ \\
For null-homologous knots and torsion $\mathfrak{s}$ there is a canonical 
choice of filtration index found from the first Chern class, \cite{Knot}. 
Suppose all 
the knots in $Y$ found by closing strands in $S$ are null-homologous, 
and that we have a Heegaard diagram where the intersection point on the 
each meridian
is fixed. Let $\la_i$ be a longitude for the closure of the $i^{th}$ 
strand, $L_{i}$. 
This curve can be realized in the Heegaard 
diagram for the string link as a curve crossing only one $\al$-curve, the 
$i^{th}$ meridian. Interchanging the meridian with this longitude gives a 
Heegaard
diagram for the manifold found by performing $0$-framed surgery on 
$L_{i}$. To each intersection point in ${\bf x} \in \T_{\al} \cap 
\T_{\be}$ we
can associate an intersection point, ${\bf x}'$, for the new Heegaard 
diagram, cf. \cite{Knot}.
The Seifert surface for the closure of the $i^{th}$ strand becomes a 
doubly periodic
domain, $\mathcal{P}_{i}$, in the new diagram. Following the argument in 
\cite{Knot}
shows that we may choose

$$
\mathcal{F}_{i}({\bf x}) = \frac{1}{2} < c_{1}(\mathfrak{s}_{w}({\bf x'}), 
[\mathcal{P}_{i}]>
$$
for our filtration index. This provides a canonical choice over different
$Spin^{c}$ structures on $Y$, which the axiomatic description lacks.

\subsection{The Relationship Between Gradings and Filtration Indices}

There is another interpretation of the filtration indices, which we 
describe for a knot. If we choose $\phi \in \pi_{2}({\bf x}, {\bf y})$ we 
can lift the relative $\Zmod{\delta(\mathfrak{s})}$-grading using the 
procedure 
for filtrations:
$$
gr_{w}( [{\bf x}, i], [{\bf y}, j])  = \mu(\phi) - 2\,n_{w}(\phi) + 2( i - 
j) - <c_{1}(\mathfrak{s}_{w}), \mathcal{A}(\phi)>
$$
where $\mathfrak{s}_{w}$ is the $Spin^{c}$ structure represented by $(w, 
{\bf x})$.  When we use $z$ as the basepoint for Heegaard-Floer homology,  
we have
$$
gr_{z}( [{\bf x}, i], [{\bf y}, j])  = \mu(\phi) - 2\,n_{z}(\phi) + 2( i - 
j) - <c_{1}(\mathfrak{s}_{z}), \mathcal{A}(\phi)>
$$
calculating $\mathcal{A}(\phi)$ using the complete set of paths for $w$.
If we add a periodic domain for $z$, the 
term with $n_{z}$ does not change. However, $\mathcal{P}_{z} = 
\mathcal{P}_{w} + r\,[S]$, so $\mathcal{A}(\mathcal{P}_{z}) =$ 
$\mathcal{A}(\mathcal{P}_{w})$ and this expression is independent of 
$\phi$. The difference 
$$
gr_{z}( [{\bf x}, i], [{\bf y}, j])   - gr_{w}( [{\bf x}, i], [{\bf y}, 
j])  = 2(n_{w} - n_{z})(\phi) - <c_{1}(\mathfrak{s}_{w}), 
\mathcal{A}(\phi)> + <c_{1}(\mathfrak{s}_{z}), \mathcal{A}(\phi)> 
$$
is equal to $2(\mathcal{F}({\bf y}) - \mathcal{F}({\bf x}))$  since  
$\mathfrak{s}_{z} - \mathfrak{s}_{w} =  [K]$ and 
$n_{z}(\mathcal{A}(\phi)) =$ $[K] \cap \mathcal{A}(\phi)$.
We can adapt this discussion to string links by considering
each coordinate separately. In effect, the 
filtration indices are measuring the difference in relative gradings for 
an intersection point induced by different basepoints.

\subsection{The Differential and the $\Gamma$-Sub-Complex}

As in \cite{Hom3} there is a differential, $\partial^{\circ}$ on 
$CF^{\circ}_{\Gamma}(Y; \mathfrak{s})$
for $\circ$ equal to $\infty$, $+$, $-$, or $\widehat{\cdot}$, defined by 
the linear extension
of

$$
\partial [{\bf x}, i, \overline{v}] = 
\sum_{{\bf y} \in \mathcal{I}(\mathfrak{s})} \sum_{\begin{array}{c} \phi 
\in \pi_2({\bf x},{\bf y}) \\ 
\mu(\phi)=1 \end{array}}
\# \widehat{\mathcal{M}}(\phi) [{\bf y}, i\,-\, n_w(\phi), \overline{v} - 
n_{\overline{z}}(\phi) + n_{\overline{z}}(\mathcal{A}(\phi))]
$$
where the signed count is made 
with respect to a choice of a coherent system of orientations. 
We verify that this is a differential below.  The differential
is a $\Z[U]$-module map when $\circ$ is $\pm$ or $\infty$. Since
$i - n_{w}(\phi) \leq i$ when 
$\phi$ admits holomorphic representatives, we see that the sub-group 
$CF^{-}_{\Gamma}$ is a sub-complex. The differential on $CF^{+}_{\Gamma}$ 
makes it into a quotient complex.  For 
$\widehat{CF}_{\Gamma}$  we restrict to those generators with $i=0$ 
and those $\phi$ with $n_{w}(\phi) = 0$. 
\\
\\
When the lattice $\Lambda_{(Y, \Gamma)} \equiv 0$ this complex is 
$\Z^{k+1}$ filtered by the relation

$$
(i, j_1, \ldots, j_k) < (i', j'_1, \ldots, j'_k)
$$
when $i < i'$ and $j_l < j'_l$ for all $l$. We have $\partial [{\bf x}, i, 
\overline{v}] \leq [{\bf x}, i, \overline{v}]$ since $n_{w}(\phi), 
n_{z_{i}}(\phi) \geq 0$ on classes represented by a holomorphic disc. (The 
partial ordering 
on non-zero linear combinations of generators, $\sum {\bf y}_{i} \leq$ 
$\sum {\bf x}_{j}$, occurs when every generator ${\bf y}_{i}$ $ \leq {\bf 
x}_{j}$ for each $j$). When $\Lambda \not\equiv 0$ the additional terms in 
the differential can disrupt the monotonicity of the indices, precluding 
similar 
filtrations. For $\widehat{CF}_{\Gamma}$ and $\Lambda \equiv 0$ there
is a $\Z^{k}$ filtration defined analogously for the indices from the
vector, $\overline{v}$.  
\\
\\
We choose a complete set of paths and a multi-point filtration 
index $\overline{\mathcal{F}}$.  Consider the sub-group of 
$CF^{\infty}_{\Gamma}(Y; \mathfrak{s})$
generated by those $[{\bf x}, i, \overline{v}]$ with
$$
 (v_1,\, v_2,\, \ldots\, ,\, v_k)  - \, (i,\, i, \, \ldots\, 
,\, i) 
 = \mathcal{F}({\bf x})
$$
For a choice of $\mathcal{A}$,  there is a unique $k$-tuple $(v_1, \ldots, 
v_k)$ associated
with $[{\bf x}, i]$ giving an element of this sub-group.  When $\Lambda 
\equiv 0$ this $k$-tuple induces a $\Z^{k}$ filtration on  the chain 
complexes for the Heegaard-Floer homology, $CF^{\circ}(Y, \mathfrak{s})$ 
in \cite{Hom3}. 
We define the filtration by $(v_1, 
\ldots, v_k)$ $\le (c_1, \ldots, 
c_k)$ for each fixed $k$-tuple $(c_1, \ldots, c_{k})$.
\\
\\
The filtration index relation ensures that this sub-group is a 
sub-complex. In particular,
if ${\bf y}$ is in the boundary of ${\bf x}$, i.e. $<\partial {{\bf x}}, 
{\bf y}> \neq 0$, 
then
$$
\begin{array}{c}
\ \\
\overline{\mathcal{F}}({\bf y}) = \overline{\mathcal{F}}({\bf x}) + 
(n_{\overline{w}} - n_{\overline{z}})(\phi)  + 
n_{\overline{z}}(\mathcal{A}(\phi)) \\
\ \\
=  (v_1,\, \ldots\, ,\, v_k) - (i,\, \ldots\, ,\, i) + ((n_{w} -n 
_{z_1})(\phi), \ldots, 
(n_{w}-n_{z_1})(\phi)) + n_{\overline{z}}(\mathcal{A}(\phi))\\
\ \\
= \big( \overline{v} - n_{\overline{z}}(\phi)+ 
n_{\overline{z}}(\mathcal{A}(\phi))\big) - \big( i - 
n_{\overline{w}}(\phi)\big) \\
\ \\
\end{array}  
$$
Thus $[{\bf y}, i\,-\, n_w(\phi), \overline{v} - n_{\overline{z}}(\phi) + 
n_{\overline{z}}(\mathcal{A}(\phi))]$ still satisfies the sub-complex 
condition. 
Furthermore, the action of $U$ preserves the sub-complex, affording it 
the structure of a $\Z[U]$-module. We call this sub-complex 
$CF^{\infty}(Y, \Gamma; \mathfrak{s})$. 

\begin{lemma}
The map $\partial : CF^{\infty}(Y, \Gamma; \mathfrak{s}) \ra 
CF^{\infty}(Y, \Gamma; \mathfrak{s}) $ is a differential.
\end{lemma}
\noindent
{\bf Proof:} This follows from \cite{Hom3} with almost no alteration. We 
consider $\psi \in 
\pi_2({\bf x_0},{\bf x_2})$ with 
a moduli space of holomorphic representatives of dimension $\mu(\psi) = 
2$. We already know that the components of the boundary of  
$\mathcal{M}(\psi)$ which contribute to $\partial^{2}$ are of the form 
$\mathcal{M}(\phi_1) \times \mathcal{M}(\phi_2)$ for classes $\phi_{1} \in 
\pi_{2}({\bf x_0}, {\bf x_1})$ and $\phi_{2} \in \pi_{2}({\bf x_1}, {\bf 
x_2})$. The other possible ends are eliminated as in \cite{Hom3}. We 
choose our
almost complex structures to exclude the bubbling of spheres at the 
intersection points. The introduction of marked points in the diagram does 
not alter this. On the other hand, we know that boundary bubbles will 
cancel in the summation of coefficients. This occurs for counts in  
specified degenerate homotopy classes, and all representatives of the 
homotopy class will induce the same change in the indices: $i, v_{1}, 
\ldots, v_{k}$. All other degenerations are excluded for dimension 
reasons.

Since every homotopy class of discs satisfies the filtration index 
relation so will each of the classes $\psi$, $\phi_{1}$ and $\phi_{2}$. In 
particular, every boundary component of $\mathcal{M}(\psi)$ of the form 
$\phi_{1} \ast \phi_{2}$ contributes to $\partial^{2}$. We know that
$$
n_w(\psi) = n_w(\phi_1) + n_w(\phi_2) 
$$
$$
n_{z_i}(\psi) = n_{z_i}(\phi_1) + n_{z_i}(\phi_2)
$$
$$
\mathcal{A}(\psi) = \mathcal{A}(\phi_1) + \mathcal{A}(\phi_2)
$$
As a result, different ends of the compactification of 
$\widehat{\mathcal{M}}(\psi)$ are taken to the same element $[{\bf y}, i, 
\overline{v}]$ in $\partial^{2}$ and still cancel after a choice of a 
coherent orientations for the moduli spaces. So $\partial^{2}= 0$ and 
defines a differential on both our complex and sub-complex.\\
$\Diamond$\\
\\
The action of $H_{1}(Y, \Z)/{\mathrm Tors}$ on the Heegaard-Floer 
homology, \cite{Hom3}, 
extends to an action on the homology of the sub-complex. 
Let $\ga \subset \Sigma$ be a simple, closed curve 
representing  the 
non-torsion class $h \in H_1$ and missing every intersection point between 
an $\al$ and a $\be$. Let $a(\ga,\phi)$  be the intersection number in 
$\T_{\al}$ of 
$\ga \times Sym^{g-1}(\Sigma) \cap \T_{\al}$ and $u(1 + it)$ 
where $u$ represents $\phi$. 
This induces a map $\zeta \in Z^{1}(\Omega(\T_{\al}, 
\T_{\be}), \Z)$. The action of such a co-cycle is defined by the 
formula:
$$
A_{\zeta}([{\bf x}, i, \overline{v}]) = 
\sum_{\bf y} \sum_{\{\phi:\mu(\phi) =1\}} \zeta(\phi) \cdot \big( \# 
\widehat{\mathcal{M}}(\phi)\big)[{\bf y}, i - n_w(\phi), 
\overline{v} - n_{\overline{z}}(\phi  - \mathcal{A}(\phi))] 
$$
If $\Lambda_{Y}$ is trivial then the map $A_{\zeta}$ is a filtered chain 
morphism.  In addition, $A_{\zeta}$ preserves the sub-complex 
$CF^{\infty}(Y,\Gamma;\mathfrak{s})$, as $\partial$ does. 

Following the argument presented above that $\partial^{2} = 0$, we can 
verify that this is a chain map. We have the analog of Proposition 4.17 of 
\cite{Hom3}:

\begin{thm} 
There is a natural action of the exterior algebra, 
$\Lambda^{\ast}(H_1(Y,\Z)/{\mathrm Tors})$ on the homology $HF^{\infty}(Y, 
\Gamma; 
\mathfrak{s})$, where $\zeta \in H_1(Y,\Z)/{\mathrm Tors}$ lowers degree 
by $1$ and 
induces a filtered morphism of the chain complex when $\Lambda_{Y} \equiv 
0$.
\end{thm} 
\noindent
{\bf Proof:} To see that this is a chain map note that
the formula in lemma 4.17 of  \cite{Hom3} 
for the coefficients of $\partial A_{\zeta} \pm A_{\zeta}\partial$ still 
applies as it only depends upon the $\phi$'s and not upon the 
additional indices. As for the differential, any 
$\phi$ used in $\pi_2({\bf x}, {\bf z})$ with $\mu(\phi)=2$ will give the 
same set of  indices for ${\bf z}$. The same observation applies to lemmas 
4.18 and 4.19 of 
\cite{Hom3}. 
\\
\\
Various other chain complexes may be defined from the above construction. 
When $\Lambda \equiv 0$ we may require any sub-set of the indices to be 
less than zero to get a sub-complex. Or we may take the quotient by this
sub-complex. We can require that $i=0$ and look at holomorphic discs with 
$n_w = 0$, or also
disallow  discs that cross  some of the $z_i$. When we require 
$n_{z_{i}}(\phi) =0$ for
some or all of the marked points, we denote the resulting complex by 
$\widehat{CF}$.  We then examine  the sub-complexes generated by the 
intersection points inducing a given filtration index for the $z_i$ we 
have disallowed. If we allow the disc to intersect the marked points 
associated with the components in a sub-link $\mathbb{L}_{2}$, but not 
all marked 
points,
then the filtration indices from  $\mathbb{L}_{2}$ define a 
$\Z^{|\mathbb{L}_{2}|}$- 
filtration  on each complex $\widehat{CF}(Y, \Gamma - \mathbb{L}_{2}; 
\overline{v}_{1})$, where $\overline{v}_{1}$ records the fixed filtration 
indices for the complementary sub-complex. When $\Lambda \not\equiv 0$ the 
only 
obvious sub-complexes are those that depend upon
the $i$ index, as in Heegaard-Floer homology. However, it is still 
meaningful to look at the sub-complex where $i=0$ and the differential 
includes only those $\phi$ with $n_{w}(\phi)=0$ and $n_{z_i}(\phi)=0$ for 
 all of the $z_i$. The differential in this complex is

$$
\widehat{\partial}[{\bf x}, \overline{j}] = \sum_{y}\sum_{\phi} 
\#\widehat{\mathcal{M}}(\phi)[{\bf y}, \overline{v} + 
n_{\overline{z}}(\mathcal{A}(\phi))]
$$
The $\overline{v}$ index may change, but only by an element of $\Lambda$. 
Thus
the sub-complex differential preserves the summands:

$$
\bigoplus_{\overline{v}\equiv\overline{v}'\, {\mathrm mod}\, \Lambda} 
CF(Y,\Gamma; \mathfrak{s}, \overline{v})
$$
\noindent
For null-homologous links, this definition reduces to the complex  
$\widehat{CF}(Y, \Gamma; \mathfrak{s}, \overline{v})$. The
equivalent when $\Lambda \not\equiv 0$ should be a complex 
$\widehat{CF}(Y, 
\Gamma; [\overline{v}])$ where the
last entry  is a coset of  $\Z^{k}/\Lambda$. The specific filtration index 
prescribes a specific coset; however, only the relative difference between 
cosets will be invariant. 
\\
\\
{\bf Convention:} For the rest of the paper, $\widehat{CF}(Y, \Gamma; 
\mathfrak{s}, [\overline{v}])$
is the complex with differential counting $\phi$'s with $n_{w}(\phi)=0$ 
and $n_{z_i}(\phi)=0$ for 
all of the $z_i$. When $\Lambda \equiv 0$, it is the $E^{1}$ term in the 
spectral 
sequence defined by the $\Z^{k}$ filtration on $\widehat{CF}_{\Gamma}(Y; 
\mathfrak{s})$. 
\\
\\
Altogether we will have the following theorem, an extension of the theorem 
in \cite{Knot} which provided the statement for null-homologous knots. The 
union $\Gamma_{1} \cup \Gamma_{2} = \Gamma$ requires that $\Gamma_{j}$ be 
found from $\Gamma$ by ignoring a set of  link components, sets that are 
disjoint for $j = 1, 2$. 

\begin{thm}
For $\Gamma \subset Y$ coming from a $d$-based link, the 
homology $\widehat{HF}(Y, 
\Gamma; \mathfrak{s})$ is a relatively $\Z^{k}/\Lambda$ indexed invariant 
of $\mathbb{L} \cup D$ up to isotopy which is a direct sum \\ 
$\oplus_{[\lambda]} 
\widehat{HF}(Y, 
\Gamma; \mathfrak{s}, [\lambda])$ of invariant  
sub-groups. There is a natural action of $H_{1}(Y; \Z)/ \mathrm{Tors}$ on 
each of the factors. When $\Gamma_{1} \cup \Gamma_{2} = \Gamma$, and 
$\Lambda_{\Gamma_{2}} \equiv 0$, the 
filtration indices for $\Gamma_{2}$ may be defined using the first Chern 
class. In this case, the $\Z^{|\Gamma_{2}|}$ filtered chain homotopy type 
of $CF^{\Gamma_{2}}(Y,  \Gamma_{1}; \mathfrak{s}, [\overline{j}])$ is an 
invariant of $\Gamma$ for each coset, $[\overline{j}]$, of 
$\Z^{|\Gamma_{1}|}/\Lambda_{\Gamma_{1}}$.
\end{thm} 
\noindent
{\bf Comments:}
\ben
\item Although it seems plausible, when $\Lambda \not\equiv 0$, to 
consider 
those $\phi$ with $n_{\overline{z}}(\phi) \in \Lambda$, this does not 
prescribe a differential. That 
$\phi = \phi_{1} \ast \phi_{2}$ does not
imply that $n_{\overline{z}}(\phi_{j}) \in \Lambda$. Hence the terms in 
$\partial^{2}$ may 
give rise to complementary boundaries one of which does not arise 
from the definition of the differential.  

\item The complexes for $\Gamma$ are similar to the twisted 
coefficients of \cite{3Man}. Indeed, the proof of invariance parallels 
that for the twisted coefficients. However, twisted coefficients are used 
to distiguish homotopy classes that are otherwise indistiguishable; we
cancel distinctions which would otherwise appear.  

\item The action of $H_{1}(Y; \Z)$ may be extended to the various 
sub-complexes discussed. The definition mimics that for the differential. 
In particular, we may limit ourselves to $\phi$ with $n_{w}(\phi) = 0$ and 
$n_{z_{i}}(\phi) = 0$. When $\Lambda \equiv 0$, this action will be 
natural as filtered morphisms of chain complexes.
\een

\section{Maps of Multi-Pointed Homologies}
\label{sec:Chain}

We may extend the theory of maps on $HF^{\circ}$ induced by cobordisms of 
three manifolds.   

\subsection{More about Complete Sets of Paths}
\label{sec:complete}
Let $X$ be a Heegaard triple $\Sigma_{\al\be\ga}$. We wish to choose an 
assignment $\mathcal{A}$ from triangles representing certain $Spin^c$ 
structures to $H_{2}(X, \Z)$ that is 

\ben
\item[] I. Additive under splicing and addition of periodic domains.
\item[] II. Restricts to each boundary of the triple as a specified 
additive assignment from a complete set of paths for that component.
\een
Let $\mathfrak{u}$ be a $Spin^{c}$ structure for $X$. 
We show how to define an additive assigment for the homotopy classes of triangles
representing an orbit of $\mathfrak{u}$ under the action of $H^{2}(X, \partial X; 
\Z)$.

Given a single triangle 
$\psi_{0} \in \pi_{2}({\bf x_{0}}, {\bf  y_{0}}, {\bf z_{0}})$ with 
$n_{w}(\psi_{0})=0$, representing $\mathfrak{u}$, and  complete sets of 
paths $\mathcal{A}_{\al\be}$, 
$\mathcal{A}_{\ga\al}$, and $\mathcal{A}_{\ga\be}$ for the $Spin^{c}$ 
structures it restricts to at the boundary, we may define a triply 
periodic domain for any triangle $\psi \in \pi_{2}({\bf x}, {\bf y}, {\bf 
z})$ which restricts to the boundaries as the same $Spin^{c}$ structures 
as $\psi_{0}$. Take $\phi_{\bf x}$, $\phi_{\bf y}$, and $\phi_{\bf z}$ in 
their respective complete set of paths and consider the triply periodic 
domain such that 

$$\phi_{\bf z}  \ast \psi_{0} + \mathcal{T} = \psi \ast (\phi_{\bf x} 
\otimes  
\phi_{\bf y}) + r [S]$$
Then $\mathcal{A}_{\al\be\ga}(\psi) = [\mathcal{T}]$ in $H_{2}(X; \Z)$.  

Alteration of $\psi$ by a class $\phi$ in one of the boundaries induces 
the
relation:
$$
\mathcal{A}_{\al\be\ga}(\psi \ast 
\phi) = \mathcal{A}_{\al\be\ga}(\psi) + \mathcal{A}_{\circ}(\phi)
$$
where $\circ$ should be replaced with the pair of subscripts 
designating that boundary. For the particular case that $\psi 
= \phi_{\ga\be} \ast \psi_{0} \ast (\phi_{\al\be} \otimes 
\phi_{\ga\al})$ 
where $\phi_{\al\be} \in \pi_{2}({\bf x}, {\bf x_{0}})$, $\phi_{\ga\al} 
\in \pi_{2}({\bf y}, {\bf y_{0}})$, $\phi_{\ga\be} \in \pi_{2}({\bf 
z_{0}}, {\bf z})$, we have:
$$
\mathcal{A}_{\al\be\ga}(\psi) = \mathcal{A}_{\ga\be}(\phi_{\ga\be}) + 
\mathcal{A}_{\al\be}(\phi_{\al\be})   +
\mathcal{A}_{\be\ga}(\phi_{\ga\al})
$$
Finally, alteration of $\psi_{0}$ without changing the basepoints on the 
boundary components changes the identification with homology classes by 
adding the class of a triply periodic domain. 

\subsection{Push-Forward Filtration Indices}

In analogy to the three dimensional case, let $\Lambda_{\al\be\ga}$ be 
the lattice in $\mathbb{Z}^k$ of vectors $(n_{z_{1}}(\mathcal{T}), \ldots, 
n_{z_{k}}(\mathcal{T}))$ where $\mathcal{T}$ is any triply periodic 
domain. With additive identifications on the boundaries compatible with 
that on 
the homotopy classes of triangles, we can push forward the filtration 
indices from $Y_{\al\be}$ and $Y_{\ga\al}$. We define

$$
\overline{\mathcal{G}}({\bf z}) = \overline{\mathcal{F}}_{\al\be}({\bf x}) 
+ \overline{\mathcal{F}}_{\ga\al}({\bf y}) + (n_w - 
n_{\overline{z}})(\psi) + n_{\overline{z}}(\mathcal{A}(\psi))
$$
\noindent
where $\psi \in \pi_{2}({\bf x}, {\bf y}, {\bf z})$. The expression on the 
right does not change under the addition of any triply periodic domains or 
the class $[S]$. Thus it does not depend upon the specific homotopy class 
of triangles joining three given intersection points. Nor indeed does it 
change if we use a triangle abutting ${\bf z}$ but with different intial 
intersection points. Since

$$
\overline{\mathcal{F}}_{\al\be}({\bf x}) - 
\overline{\mathcal{F}}_{\al\be}({\bf x'}) = (n_w - 
n_{\overline{z}})(\phi_{1}) + 
n_{\overline{z}}(\mathcal{A}_{\al\be}(\phi_{1}))
$$
$$
\overline{\mathcal{F}}_{\ga\al}({\bf y}) - 
\overline{\mathcal{F}}_{\ga\al}({\bf y'}) = (n_w - 
n_{\overline{z}})(\phi_{2}) + 
n_{\overline{z}}(\mathcal{A}_{\ga\al}(\phi_{2}))
$$
when $\psi' = \psi \ast \phi_{1} \ast \phi_{2}$, we find

$$
\overline{\mathcal{G}}({\bf z}) = \overline{\mathcal{F}}_{\al\be}({\bf 
x'}) + \overline{\mathcal{F}}_{\ga\al}({\bf y'}) + (n_w - 
n_{\overline{z}})(\psi') + n_{\overline{z}}(\mathcal{A}(\psi'))
$$
since 
$$
\mathcal{A}_{\al\be\ga}(\psi)  + \mathcal{A}_{\al\be}(\phi_{1}) + 
\mathcal{A}_{\ga\al}(\phi_{2}) = \mathcal{A}_{\al\be\ga}(\psi')
$$
\noindent

We may then check that $\overline{\mathcal{G}}$ is a filtration index for 
the $Spin^{c}$ structure on the $\ga\be$ boundary. Let $\phi \in 
\pi_{2}({\bf z}, {\bf z'})$ and $\psi \in \pi_{2}({\bf x_{0}}, {\bf 
y_{0}}, {\bf z})$, then 

$$
\overline{\mathcal{G}}({\bf z'}) - \overline{\mathcal{G}}({\bf z})  = (n_w 
- 
n_{\overline{z}})(\phi) + n_{\overline{z}}(\mathcal{A}_{\al\be\ga}(\psi')) 
- n_{\overline{z}}(\mathcal{A}_{\al\be\ga}(\psi))
$$
where $\psi' = \phi \ast \psi$. But 
$\mathcal{A}_{\al\be\ga}(\psi') = \mathcal{A}_{\al\be\ga}(\psi) + 
\mathcal{A}_{\ga\al}(\phi)$.
\\
\\
{\bf Note:} If we require the filtration index to be $0$ on the 
basepoints, ${\bf x_0}$ and ${\bf y_{0}}$, we may use $\psi_{0}$ to 
calculate the value of $\overline{\mathcal{G}}({\bf z_{0}}) =$ 
$-n_{\overline{z}}(\psi_{0})$. We will assume, unless otherwise stated, 
that
$n_{\overline{z}}(\psi_{0}) = \overline{0}$.   
\\

\subsection{Constructing Additive Assignments}

Suppose we have a four manifold $X_{\al\be\ga}$ defined by a Heegaard 
triple. We can consider the long exact sequence 
$$
\cdots \ra H_{3}(X, Y_{\ga\be}) \ra H_{2}(Y_{\ga\be}) \ra H_{2}(X) \ra 
H_{2}(X, Y_{\ga\be}) \ra \cdots
$$
$H_{2}(Y_{\ga\be})$ injects into $H_{2}(X)$ for a Heegaard 
triple because each doubly periodic domain is also a triply periodic 
domain. Furthermore,  the image of $H_{2}(X)$ inside $H_{2}(X, 
Y_{\ga\be})$ is a free group as it consists of those triply periodic 
domains with non-trivial $\al$-boundary; it is finitely generated and 
torsion free, since no multiple except $0$ can eliminate that boundary.  
We may choose a splitting $H_{2}(X) \cong  H_{2}(Y_{\ga\be}) \oplus C$, 
which we hope will reflect our Heegaard diagram. 

Suppose we have additive assignments $\mathcal{A}_{\al\be}$, 
$\mathcal{A}_{\ga\al}$, and a specified basepoint ${\bf z_{0}}$ in 
$\T_{\ga}\cap \T_{\be}$. Suppose further that we have chosen a single 
triangle $\psi_{0} \in \pi_{2}({\bf x_{0}}, {\bf  y_{0}}, {\bf z_{0}})$ 
and  for each ${\bf z} \in \T_{\ga} \cap \T_{\be}$, there is a preferred 
triangle, $\psi_{\bf z} \in \pi_{2}({\bf x}, {\bf  y}, {\bf z})$, abutting 
${\bf z}$ and  representing a $Spin^{c}$ structure whose restrictions to 
the three boundary components are the same as those of $\psi_{0}$. We will 
assume that all these triangles have $n_{w}(\psi_{\bf z}) = 0$. This induces a 
complete set of paths on $\Sigma_{\ga\be}$ for the 
$Spin^{c}$ structure represented by ${\bf z_{0}}$.  

For each ${\bf z}$ and 
$\phi \in \pi_{2}({\bf z_{0}},{\bf z})$, $n_{w}(\phi) = 0$, there is a 
triply periodic domain satisfying $\psi_{\bf z} \ast (\phi_{\bf x} \otimes 
\phi_{\bf y}) =$ $\phi \ast \psi_{0} + \mathcal{T}$. Using the splitting 
we may divide $[\mathcal{T}] = [\mathcal{P}] \oplus [\mathcal{T}']$. 
If we subtract $\mathcal{T}'$ from   $\psi_{\bf z} 
\ast (\phi_{\bf x} \otimes \phi_{\bf y}) $, we have a triangle which 
differs from $\phi \ast \psi_{0}$ by a doubly periodic domain in 
$Y_{\ga\be}$. We add this periodic domain to $\phi$ to get $\phi_{\bf z}$. 
Different choices of $\phi$ produce the same choice for $\phi_{\bf z}$ 
relative to the splitting. Suppose $\phi' = \phi + \mathcal{P}_{\ga\be} + 
r[S]$. We can ignore the term involvong $[S]$.  Then $\mathcal{T}_{\phi'} 
= \mathcal{T}_{\phi} - \mathcal{P}$. Projecting to $H_{2}(Y_{\ga\be})$, we 
have $\mathcal{P}_{\phi'} = \mathcal{P}_{\phi} - \mathcal{P}$. Adding this 
to $\phi'$ shows that $\phi'+  \mathcal{P}_{\phi'} =$ $\phi+  
\mathcal{P}_{\phi}$.

\subsection{Chain Maps}

As in \cite{Hom3} and \cite{Four}, we start with a Heegaard triple 
defining a four manifold 
$X_{\al\be\ga}$. We assume that we have complete sets of paths on the 
boundaries for the restriction of a $Spin^{c}$ structure, $\mathfrak{u}$, 
and that 
there is a compatible 
map $\mathcal{A}$ for $\Si_{\al\be\ga}$. We may define a multi-point chain 
map for 
the $Spin^{c}$ structure:
$$
F_{\mathfrak{s}}([{\bf x}, i_{1}, \overline{j}_{1}], [{\bf y}, i_{2}, 
\overline{j}_{2}]) =
\sum_{\bf z} \sum_{\psi} \#\mathcal{M}(\psi) [{\bf z}, i_{1}+ i_{2} - 
n_{w}(\psi), \overline{j}_{1} + \overline{j}_{2} - n_{\overline{z}}(\psi) 
+ n_{\overline{z}}(\mathcal{A}_{\al\be\ga}(\psi))]
$$
where $\psi$ is class representing $\mathfrak{s}$ with $\mu(\psi) = 0$.

That this is a chain map follows from the usual arguments by examining 
ends of moduli spaces with $\mu(\psi') = 1$. The identities 
for compatibility of additive assignments imply that for an end modelled 
upon
$$
\mathcal{M}(\psi_{\al\be\ga}) \times \mathcal{M}(\phi_{\ga\al})
$$
we have 
$$
\mathcal{A}_{\al\be\ga}(\psi') = \mathcal{A}_{\al\be\ga}(\psi) + 
\mathcal{A}_{\ga\al}(\phi_{\ga\al})
$$
As for the differential these identities allow us to conclude
that the various cancellations necessary for this to specify a chain map 
will still occur.  

Furthermore, the map is $U$ invariant as the moduli spaces do not depend 
upon the particular $i$ or $j_{l}$. Finally, the map preserves the 
sub-complexes defined from the filtration relations, where the filtration 
on the $\ga\be$-boundary is the push-forward of those on the $\al\be$ and 
$\ga\al$-boundaries for a choice of a complete set of paths.
\\
\\
{\bf Note:} The next section outlines the invariance of the homology 
groups
under the reduced Heegaard equivalences. 
We should additionally check the invariance of the cobordism maps on 
homology 
under the various alterations: invariance of the almost complex 
structure, isotopies of attaching circles, handleslides, and 
stabilizations. Using the observations of the next section, these proofs
follow directly. The compatibility of the maps, 
$\mathcal{A}$, ensure that the new indices do not disrupt the chain 
homotopy identities found in \cite{Hom3}. The following sub-section 
relates 
the associativity properties of these maps. Once that is done, we can 
recover 
the arguments in the first sections of \cite{Four} by slight 
modifications. The naturality of the strong equivalence maps follows from 
these considerations. Specifically, under our restricted set of Heegaard 
equivalences, the chain maps for  two presentations of $(Y, \Gamma)$ 
commute with the maps induced by strong equivalences on homology.

\subsection{Associativity with Multi-point Filtrations}

\ \\
We consider two triples $(\Sigma, \seta, \setb, \setg; w)$ and 
$(\Sigma, 
\setg, \setb, \setd; w)$ where we have equipped $\Sigma$ with marked 
points $z_{1}, \ldots, z_{k}$. We assume a choice of $Spin^{c}$ 
structures, $\mathfrak{s}_{1}$ and $\mathfrak{s}_{2}$, for the respective 
triples and a complete sets of paths for $\Sigma_{\al\be}$, 
$\Sigma_{\ga\be}$, $\Sigma_{\de\be}$, $\Sigma_{\ga\al}$, and 
$\Sigma_{\ga\de}$, as well as $\psi_{\al\be\ga} \in \pi_{2}({\bf x_{0}}, 
{\bf y_{0}}, {\bf u_{0}})$ and $\psi_{\ga\be\de} \in  \pi_{2}({\bf u_{0}}, 
{\bf z_{0}}, {\bf w_{0}})$ representing the $Spin^{c}$ structures. Let 
$\mathfrak{s}$ be the $Spin^{c}$ structure, $\mathfrak{s}_{1} \# 
\mathfrak{s}_{2}$ on the quadruple, $X_{\al\be\de\ga}$.  

In this paper we consider chain maps resulting from diagrams which 
represent
surgeries, handleslides, or the other Heegaard equivalences. For these 
diagrams, each of the three manifolds described by $\Sigma_{\ga\al}$, 
$\Sigma_{\de\ga}$, and 
$\Sigma_{\de\al}$ will be homeomorphic to a connected sum of $S^{3}$'s and 
$S^{1}\times S^{2}$'s, although the additional basepoints may prevent the 
reduced Heegaard equivalences from converting the diagram into the 
standard picture. 
The only triply periodic domains for this triple will be sums of doubly 
periodic domains from $\Sigma_{\ga\al}$ and $\Sigma_{\ga\de}$. 
Furthermore, every doubly periodic domain, $\mathcal{P}_{\al\de} \subset$ 
$\Sigma_{\al\de}$, will be a linear combination of doubly periodic domains 
from $\Sigma_{\ga\al}$ and $\Sigma_{\ga\de}$. 

Each such diagram will possess only one torsion $Spin^{c}$ structure which will 
admit a special, closed generator for its chain group: $\Theta^{+}$. All 
our homotopy classes of triangles and quadrilaterals will be required to 
use those special generators when available. These generators will always have filtration index 
equal to $\overline{0}$ and be basepoints for their respective complete 
set of paths. We denote by $\mathfrak{s}_{0}$ the unique $Spin^{c}$ 
structure represented
by classes in $\pi_{2}(\Theta^{+}_{\al\de}, \Theta^{+}_{\de\ga}, 
\Theta^{+}_{\ga\al})$. The geometry of the Heegaard diagram will often 
provide a splitting $H_{2}(Y_{\ga\al}) \oplus H_{2}(Y_{\ga\de}) $ as 
$H_{2}(Y_{\al\de}) \oplus \Z^{m}$ and a homotopy class of triangles, 
$\psi_{\Theta}$ abutting the special intersection points in 
$X_{\al\de\ga}$. We will use this data in the construction of a complete 
set of paths for the $Spin^{c}$ structure represented by 
$\Theta^{+}_{\de\al}$. 

Given $\mathfrak{s}$  let 
$\mathfrak{G}$ be its orbit under the action of  $i_{\ast}\big( 
H_{2}(Y_{\be\ga}; \Z) \big)$. 
We may then extend the coherent systems of orientations for the triples to 
one for quadruple 
$X_{\al\be\de\ga}$ and this orbit, see \cite{Hom3}. We also assume that 
the quadruple is strongly/weakly admissible, as necessary. 

The choices made above define maps, $\mathcal{A}$, on the boundaries of 
the 
quadruple and on the two original triples. As before, we may define such a 
map on the quadruple: $\mathcal{A}_{\al\be\ga\de}(\psi)$ as the homology 
class of
the sum of doubly periodic domains from pairs of $\{\al, \be, \ga, \de \}$ 
which must be added to $\phi_{\bf w} \ast \psi_{0}$ to get $\psi \ast 
\big( \phi_{\bf x} 
\otimes \phi_{\bf y} \otimes \phi_{\bf z}\big) + r[S]$. Define the map

$$
\begin{array}{c}
H_{\mathfrak{G}}([{\bf x}, i_{1}, \overline{j}_{1}] \otimes [{\bf y}, 
i_{2}, \overline{j}_{2}] \otimes [{\bf z}, i_{3}, \overline{j}_{3}] ) =
$$
\\
\\
$$
\sum_{\bf w} \sum_{\{ \psi | \mu(\psi) = 0\} } \#\mathcal{M}(\psi) [{\bf 
w}, i_{1} + i_{2} + i_{3} - n_{w}(\psi), \overline{j}_{1} + 
\overline{j}_{2} + \overline{j}_{3} - n_{\overline{z}}(\psi) + 
n_{\overline{z}}(\mathcal{A}_{\al\be\ga\de}(\psi))]
\end{array}
$$
where $\psi$ represents a $Spin^{c}$ structure from $\mathfrak{G}$. 
This map induces a chain homotopy
$$
F_{\mathfrak{s}_{1}} \circ F_{\mathfrak{s}_{2}} - \sum_{\{ \mathfrak{u} 
\in \mathfrak{G}\} } F_{\mathfrak{u}|X_{\al\be\de}} \circ 
F_{\mathfrak{s}_{0}} = \partial_{\de\be} \circ H \pm H \circ 
\partial_{\al\be \otimes \ga\al \otimes \de\ga}
$$
in the standard way, \cite{Hom3}, thereby establishing associativity for 
our cobordism maps
if we can find compatible additive assignments. (In the above formula, 
$F_{\mathfrak{s}_{0}} = F(\Theta^{+}_{\al\de} \otimes 
\Theta^{+}_{\ga\al})$ which, under our assumptions, if
almost always $\pm \Theta^{+}_{\ga\de}$; however, with the additional 
marked points, this 
will need to be verified.) This is done as with in the proof of
associativity for twisted coefficients.

In our setting we need also to check that the push forward filtrations are compatible. 
If we push forward filtrations from $Y_{\al\be}$ to 
$Y_{\ga\be}$
and then to $Y_{\de\be}$, using the basepoints for the complete sets of 
paths in the other boundaries, we find that
for the homotopy classes $\psi$ representing a $Spin^{c}$ structure on  
$\al\be\de$, the relationship is:

$$
 \overline{\mathcal{F}}_{\de\be}({\bf w}) = 
\overline{\mathcal{F}}_{\al\be}({\bf x})  
-n_{\overline{z}}(\psi_{\Theta}) + (n_w - n_{\overline{z}})(\psi) + 
n_{\overline{z}}(\mathcal{A}_{\al\be\de}(\psi))
$$
Implicit here is the 
calculation of $-n_{\overline{z}}(\psi' - \mathcal{A}(\psi'))$ 
$=-n_{\overline{z}}(\psi_{\Theta})$ since $\psi'$ abuts the three 
basepoints. When $n_{\overline{z}}(\psi_{\Theta}) = 0$ we recover
the push forward from the $\al\be\de$-diagram, so the filtraion indices
are also correct.

\subsection{Filtration Changes under Chain Maps in Various Settings}
\label{sec:Chainlamb}

\ \\Suppose $(\Sigma, \seta, \setb, \setg)$ defines
a cobordism from $Y_0$ to $Y_1$ presented as surgery on a framed link in 
$Y_0$. We will think of this as a diagram in $Y_{0}$  
with the components of $\mathbb{L}$ receiving $+\infty$ framing. We now 
elucidate
the effect of chain maps on the $\Lambda$-lattices and the filtration 
indices
in the common situations in which they arise. We will assume that any 
initial triangle, 
$\psi_{0}$, for building $\mathcal{A}$ has $n_{w}(\psi_{0})= 0$ and 
$n_{z_{i}}(\psi_{0}) = 0$ for every $i$. This condition ensures that the 
push-forward filtration is $0$ 
on the basepoint abutting $\psi_{0}$. 
\\
\\
I. When considering the complexes $\widehat{CF}(Y_0, \Gamma; 
\mathfrak{s}_{0})$ and $\widehat{CF}(Y_1, \Gamma; \mathfrak{s}_{1})$ we 
restrict to homotopy classes of triangles, $\psi$, with $n_{w}(\psi) =$  
$n_{\overline{z}}(\psi) =0$. The chain complexes are relatively 
$\Z^{k}/\Lambda_{\al\be}$ and $\Z^{k}/\Lambda_{\ga\be}$ -indexed, 
respectively, and the homology group for each filtration index is an 
invariant. Let $F$ be the chain map defined as before for a $Spin^{c}$ 
structure restricting to $\mathfrak{s}_{0}$ and $\mathfrak{s}_{1}$. Both 
$\Lambda_{\al\be}$ and $\Lambda_{\ga\be}$ inject into 
$\Lambda_{\al\be\ga}$. We may consider, by taking direct sums of the 
homologies at each end, that $\widehat{CF}(Y_{i}, \Gamma; 
\mathfrak{s}_{i})$ is $\Z^{k}/\Lambda_{\al\be\ga}$- indexed. The direct 
sums occur over the pre-images of the maps $\Z^{k}/\Lambda_{\al\be}$, 
$\Z^{k}/\Lambda_{\ga\be} \ra \Z^{k}/\Lambda_{\al\be\ga}$. The chain map
preserves the relative $\Z^{k}/\Lambda_{\al\be\ga}$ structure as the 
filtrations now satisfy

$$
\overline{\mathcal{G}}({\bf x}) = \overline{\mathcal{F}}({\bf x}) + 
n_{\overline{z}}(\mathcal{A}_{\al\be\ga}(\psi))
$$   
\noindent
II. $\Lambda_{\al\be\ga} \equiv 0$: For example, when we have a 
 string link in $S^{3}$ and the cobordism is generated by 
surgeries on curves which  are algebraically split from $S$. In 
this case, the push-forward filtration satisfies

$$
\overline{\mathcal{G}}({\bf y}) = \overline{\mathcal{F}}({\bf x}) + (n_{w} 
- n_{\overline{z}})(\psi)
$$
for every $\psi$ representing a $Spin^{c}$ structure on the cobordism 
restricting in a specified way to the ends, and for any choice of a 
complete set of paths. The filtrations on the ends are $\Z^{k}$ 
filtrations and the chain map is a filtered map for the 
$\Z^{k+1}$-filtration. This situation occurs in the long exact skein 
sequence of \cite{Knot}. In $S^{3}$, a filtration index on each component
can be calculated using the first Chern class of a $Spin^{c}$ 
structure on the manifold obtained by performing $0$ surgery on the knot.
In \cite{Knot}, \Oz\ and \Sz\ show that, in the case under consideration,
the push forward of this filtration is the one determine by the first 
Chern
class calculation on $Y_{1}$ and the formula above corresponds to their 
identity for $c_{1}$.
\\
\\
III. We will assume that no periodic domain from $\Sigma_{\ga\al}$ has 
$n_{z_i}(\mathcal{P}) \neq 0$. The diagrams encoding legitimate 
handleslides of attaching circles satisfies this assumption. Additionally, 
surgery on components of a bouquet, but not on any of the components of 
$\Gamma$, produces such a diagram. These periodic domains arise from the 
$S^1 \times S^2$ connected
sum components in $\Sigma_{\ga\al}$. Then the filtration on 
$\Sigma_{\ga\al}$ for the torsion $Spin^{c}$ structure is equivalent to 
$0$. We require our $\psi$'s to restrict to this $Spin^{c}$ structure. The 
filtration index will now satisfy

$$
\overline{\mathcal{F}}_{\ga\be}({\bf z}) = 
\overline{\mathcal{F}}_{\al\be}({\bf x}) +  (n_w - n_{\overline{z}})(\psi) 
+ n_{\overline{z}}(\mathcal{A}_{\al\be\ga}(\psi))
$$
The last term may be improved when we are considering only $\psi$'s 
representing a given $Spin^{c}$ structure. For then, 
$\mathcal{A}_{\al\be\ga}$ takes values in those triply periodic domains 
formed by summing doubly periodic domains from the boundary components. 
The value of  $n_{\overline{z}}(\mathcal{A}_{\al\be\ga}(\psi))$ is  
$n_{\overline{z}}([\mathcal{P}_{\al\be}]) + n_{\overline{z}}( 
[\mathcal{P}_{\ga\be}])$. Upon taking the quotient by $\Lambda_{\al\be} + 
\Lambda_{\ga\be} \subset \Z^{k}$ we have the same relationship as in case 
I. In the proof of invariance in the next section we will give example of 
the 
general push-forward filtration index and the construction of complete 
sets of paths for a triple.

\section{Invariance}
\label{sec:Inv}

To see that homology of these complexes are truly an invariant of the 
triple $(Y, \Gamma; \mathfrak{s})$, we need to show that performing any of 
the following moves will produce a chain homotopy equivalent complex:

\ben
\item[] $\bullet$ Handleslides and isotopies of $\setb$.
\item[] $\bullet$ Handleslides and isotopies of $\{\al_i\}_{i=k+1}^{g}$.
\item[] $\bullet$ Stabilization.
\item[] $\bullet$ Istopies of $\{\al_i\}_{i=1}^{k}$ and handleslides of 
them over element of $\{\al_i\}_{i=k+1}^{g}$.
\een

Furthermore, we are not allowed to isotope or slide over any portion of
the disc $D'$. We can, however, arrange for a $\be$ curve to isotope past 
the entire disc. The resulting diagram can be achieved by 
allowable handleslides in $\setb$ because the disc is contractible, 
\cite{Hom3}. 

We develop the proof of invariance through the multi-pointed diagrams; it 
precisely mimics the proof for Heegaard-Floer 
homology, \cite{Hom3}, and uses the same technical results.  
 
\subsection{Invariance Under Alteration of the Complete Set of Paths}

It does not matter which additive identification one uses. Let ${\bf 
x}_{0}$
be the basepoint for both $\mathcal{A}$ and $\mathcal{A}'$. 
This data determines the filtration index using the defining relation 
after
choosing a value for $\overline{\mathcal{F}}({\bf x}_{0})$, which we also 
take
to be $\overline{\mathcal{F}}'({\bf x}_{0})$. 
We define a map between the two complexes by  $[{\bf x}, i , \overline{v}] 
\ra [{\bf x}, i , \overline{v} 
+ n_{\overline{z}}(\mathcal{A}'(\phi_x))]$ where $\phi_{x}$ is the
special path found from $\mathcal{A}$. This is an isomorphism on the 
generators and is a chain map. Furthermore, the map takes the sub-complex 
defined by $\mathcal{A}$ to that defined by 
$\mathcal{A}'$. Thus it is an isomorphism on the homology groups 
determined by
the different sets of data. Note that altering $\mathcal{A}$ changes the 
filtration index for an intersection point by an element of $\Lambda$. The 
identification with cosets
remains unalterred by this isomorphism.
\\
\\
For changing ${\bf x}_{0}$ to ${\bf y}_{0}$, we may use $\phi_{y_0}^{-1}$ 
for 
$\phi'_{x_0}$ and $\phi_z \ast \phi_{y_0}^{-1}$ as $\phi'_z$ to define a 
complete set of paths and associated identification $\mathcal{A}'$. This 
does not alter $n_{\overline{z}}(\mathcal{A}(\phi))$ and changes 
$\overline{\mathcal{F}}$ 
by $\overline{\mathcal{F}}({\bf y}_0)$. The map $[{\bf x}, i, 
\overline{v}] \ra$ $[{\bf x}, i, \overline{v} +
n_{\overline{z}}(\phi_{y_0})]$ will induce a chain isomorphism from the 
homologies
defined by one complete set of paths to that defined by the other, 
although
shifting in the assignment to cosets of $\Z^{k}/\Lambda$ will occur. 

Likewise, changing the value of $\overline{\mathcal{F}}({\bf x}_{0})$ 
shifts the values
of the filtration index, and the assigment to a coset, but does not change 
any of the
homology groups or their relation to each other in the relatively indexed 
groups. 

\subsection{Results on Admissibility}

Strong/Weak admissibility can be achieved for all our diagrams without
disrupting the assumptions coming from $\Gamma$. We have seen the 
existence of
such diagrams already. 
In section 5 of \cite{Hom3} \Oz\  and \Sz\  show how to ensure that 
isotopies,
handlelsides and stabilizations can be realized through such diagrams.  
In each case this is achieved by finding a set $\setg \subset \Sigma$
of disjoint, simple closed curves with the property that $\#(\be_i \cap 
\ga_j) = \delta_{i,j}$ and that $\T_{\al} \cap 
\T_{\ga} \neq \emptyset$ (or the same but with the roles of $\al$ and 
$\be$ switched). 

We convert, through stabilization, our multi-point diagram into a diagram 
with an embedded disc. We only need require that $w$ and the entire disc 
do not intersect the winding region. 
However, we may always 
choose our $\ga$'s to lie in the disc's complement. If we wish to wind 
$\al$'s we can do the 
same, requiring only that each $\ga_i$ that intersects a meridian does so 
away from the segment in $\partial D'$. With this arrangement of 
$\setg$ the proofs of lemmas 5.4, 5.6, and 5.7 of \cite{Hom3} 
carry through. 

\subsection{Invariance of Complex Structure and Isotopy Invariance}

In these cases, we re-write the chain maps defined in \cite{Hom3}
to incorporate the new indices. For example, 
\Oz\  and \Sz\ define a chain map for a homotopy of paths of almost 
complex
structures which we adjust to be, cf. \cite{Hom3}:

$$
\Phi^{\infty}_{J_{s,t}} [{\bf x}, i \, \overline{j}] =
\sum_{\bf y} \sum_{\phi} \#\mathcal{M}_{J_{s,t}}(\phi)[{\bf y}, i - 
n_w(\phi), \overline{j} - n_{\overline{z}}(\phi - \mathcal{A}(\phi))]
$$

\noindent where the sum is over all $\phi$ with $\mu(\phi)=0$ and the 
moduli space consists of sutiable holomorphic representatives of $\phi$. 
\\
\\
Note that in this case, $\mathcal{A}$ requires no adjustment, as 
alteration of the almost complex structure does not change the 
homotopy classes of discs. That the filtration relation
holds for all homotopy classes of discs ensures this is still a chain map 
and that
the map preserves the
$\Gamma$-sub-complex. When $\Lambda_{(Y, \Gamma)} \equiv 0$ the above map 
is a filtered chain morphism by the positivity of
$n_{z_{i}}$ and the absence of the $\mathcal{A}$ term. 
\\
\\
Similar alterations ensure that the map does induce an isomorphism on 
homology (we
need to adjust the chain homotopy in \cite{Hom3} which shows that the map 
has
an inverse on homology). We need only that the trivial homotopy
class of discs from ${\bf x}$ to ${\bf 
x}$ has $\mathcal{A}(\phi)=0$ since $-\phi_x \ast \phi_x \sim 
0$ in $\pi_2({\bf x}_{0}, {\bf x}_{0})$. The invariance of the action of
$H_{1}(Y, \Z)/\mathrm{Tors}$ follows as in \cite{Hom3} adjusting the
maps as above and using our previous observations. 

For the isotopy invariance, the same argument applies (as the proofs are
roughly parallel). We write the chain map coming from the 
introduction/removal
of a pair of intersection points as:

$$
\Gamma^{\Psi}[{\bf x}, i, \overline{j}] = \sum_y \sum_{\phi \in 
\pi^{\Psi}_2({\bf x}, {\bf y}}
\# \mathcal{M}^{\Psi}(\phi) [{\bf y}, i - n_w(\phi), \overline{j} - 
n_{\overline{z}}(\phi -\overline{\mathcal{A}}(\phi))]
$$
where we count holomorphic representatives with moving boundary, 
\cite{Hom3}.
Making these adjustments as necessary, we mimic the proof in \cite{Hom3}.

Two new features occur: first,  the isotopy may
remove ${\bf x}_{0}$, and second, new intersection points need to be 
included in the complete set of paths. If we must change the basepoint for 
the complete set of paths, the homologies will be unalterred. However,
the identification with the cosets of  $\Z^{k}/\Lambda$ will alter as 
the filtration index changes by the constant vector 
$\overline{\mathcal{F}}({\bf x}'_{0})$. This explains why we have 
only a relatively indexed homology group.

When we have a pair creation, with a fixed basepoint, we get new
intersection points in pairs ${\bf q}_{+}$ and ${\bf q}_{-}$ with an
obvious holomorphic disc in $\pi_{2}({\bf q}_{+}, {\bf q}_{-})$. The 
homotopy classes of discs joining intersection points from the original 
diagram do not change in this process. Thus we may 
extend $\mathcal{A}$ by choosing $\phi_{{\bf q}_{+}}$ for each ${\bf 
q}_{+}$ and amalgamating with the newly created disc to get $\phi_{{\bf 
q}_{-}}$. We use the extended $\mathcal{A}$
in our definition of $\Gamma^{\Psi}$ as $\pi^{\Psi}_2 \cong \pi'_2$. This  
restricts to the original identification on the homotopy classes for the 
original intersection points. For a pair annihilation we do not 
necessarily have the
arrangement of preferred paths as we have just described, but we may 
change the
complete set of paths without changing the homologies to obtain it. This 
merely 
changes $\mathcal{F}({\bf q}_{-})$ by an element of $\Lambda$.  

\subsection{Invariance under Handleslides}

The standard proof for handleslide invariance applies in this context, 
using
the chain maps induced by a Heegaard triple as in the previous section. We
show below that the element corresponding to $\Theta^{+}$ is still closed
and explain which additive assignment for homotopy classes of triangles
will work.

We will describe this solely for the $\al$'s. Away from the 
curves involved in the handleslide, the resulting boundary 
$\Sigma_{\ga\al}$ is the connected sum of genus $1$ diagrams for $S^{1} 
\times S^{2}$. Because we have not moved a curve across a marked point, 
the 
corresponding multi-point diagram has all the basepoints, $w$ and $z_i$, 
in the same domain $D$ of $\Sigma - \seta - \setg$.  
\\
\\
If we calculate the multi-point filtration indices for $\Sigma_{\ga\al}$  
we see that 1) $\Lambda \equiv 0$ and 2)all $2^g$ representatives of 
the torsion $Spin^c$ structure
have filtration index $0$. None of the holomorphic discs cross any marked 
points, so the homology
is the standard homology for the connected sum of $S^1 \times S^2$'s but 
with 
generators of the form $[ \Theta^{+}, i, i, \ldots, i]$. As usual, we will 
use the canonical generator $\Theta^{+}$, the maximally graded 
generator with $i=0$. We will use strongly admissible diagrams for the 
$Spin^c$ 
structures on the ends of the cobordism. The cobordism, $X_{\al\be\ga}$, 
is $Y \times I$ so each triangle will represent $\mathfrak{s} \times I$. 

As in the proof of handleslide invariance in Heegaard-Floer homology, 
\cite{Hom3}, 
our diagrams may be drawn so that 
each intersection point in $\T_{\al} \cap \T_{\be}$ can be 
joined to an intersection point ${\bf x'} \in T_{\ga} \cap T_{\be}$ by a 
unique holomorphic triangle $\psi_{\bf x} \in \pi_{2}({\bf x}, \Theta^{+}, 
{\bf x'})$ with domain contained in the sum of the periodic domains from 
$\Sigma_{\ga\al}$.

We choose additive assignments $\mathcal{A}_1$ and $\mathcal{A}_2$ and 
their corresponding complete set of paths on both $\Sigma_{\al\be}$ and 
$\Sigma_{\ga\be}$, respectively. We require that  ${\bf x_0'}$ on 
$\Sigma_{\ga\be}$ be the intersection point abutting $\psi_{{\bf x}_{0}}$. 
Because the handleslide does not cross a meridian we have that 
$n_{w}(\psi_{\bf x_{0}}) = n_{z_{i}}(\psi_{\bf x_{0}}) = 0$. 
If we identify $\pi_2({\bf x}_{0}, {\bf x}_{0})$ 
with $\pi_2({\bf x'_0}, {\bf x'_0})$ using $\psi_{{\bf x}_{0}}$ then we 
implicitly
have an identification $H_2(Y_0) \cong H_2(Y\times I) \cong H_2(Y_1)$. Any 
other 
triangle representing the same $Spin^c$ structure can be found from this 
triangle by splicing discs in the boundaries $\Sigma_{\al\be}$ and  
$\Sigma_{\ga\be}$, by splicing doubly periodic domains from $\Si_{\ga\al}$ 
(we always require our triangles include $\Theta^{+}$), and by adding 
copies of 
$\Sigma$. If
we think of $\psi = \phi' \ast \psi_{{\bf x}_{0}} \ast \phi^{-1} + 
\mathcal{P}_{\ga\al}+ k[S]$ 
then, using the formulas from section \ref{sec:complete},

$$ \label{eqn:As}
\mathcal{A}_{Y\times I}(\psi) = \mathcal{A}_2(\phi') - \mathcal{A}_1(\phi) 
+ H(\mathcal{P}_{\ga\al})
$$  
\noindent
To the handleslide cobordism we associate the multi-point chain map 
determined
by this additive assignment on $\Si_{\al\be\ga}$. The last term 
$H(\mathcal{P})$
plays no role in the chain map as the doubly periodic domains from that 
boundary
contain no marked points in their support. With these observations the 
usual
proof applies directly.

\subsection{Stabilization Invariance}

Stabilization changes the surface $\Si$ to $\Sigma'=\Sigma \# \T^2$ and
adds an additional $\al$ and an additional $\be$ curve, intersecting in a 
single
point $c$. Stabilization does not alter $H_2$ nor does it affect the 
structures of 
$\pi_2({\bf x},{\bf y})$. Given an additive identification $\mathcal{A}$, 
we may 
extend $\mathcal{A}$ to the stabilization. Since any intersection point
must involve the additional point $c$, 
we need only that $\mathcal{A}(\phi') = \mathcal{A}(\phi)$
and  the choice ${\bf x_0'} = {\bf x}_{0}\times c$. We make the necessary 
alterations for the gluing result, Theorem 10.4 in \cite{Hom3}, to hold. 
This theorem provides the invariance under stabilization in the standard 
case. Our extension of the additive identification ensures that the 
$\Gamma$-sub-complex is preserved. Thus, we have an isomorphism from the
homology for $(\Si, \seta, \setb, {\bf x}_{0}, \mathcal{A})$ to that for 
$(\Si', \seta \cup \al_{g+1}, \setb \cup \be_{g+1}, {\bf x}_{0} \times c, 
\mathcal{A}')$. In addition, in \cite{Hom3}, \Oz\, and \Sz\, show that the 
action of $H_1(Y,\Z)/\mathrm{Tors}$ is invariant under stabilization. That 
this also applies to multi-point diagrams follows analgously to the case of the 
differential.

\section{Basic Properties of $\widehat{HF}(Y, S; \mathfrak{s})$}

\subsection{Examples}
\begin{figure}
\begin{center}
\begin{picture}(6691,3727)
\put(3607,3032){\makebox(0,0)[lb]{\smash{{\SetFigFont{10}{12.0}{rm}$K$}}}}
\put(6020,970){\makebox(0,0)[lb]{\smash{{\SetFigFont{10}{12.0}{rm}$K$}}}}
\put(6395,2945){\makebox(0,0)[lb]{\smash{{\SetFigFont{10}{12.0}{rm}$K$}}}}
\put(2973,658){\makebox(0,0)[lb]{\smash{{\SetFigFont{10}{12.0}{rm}$z$}}}}
\put(2235,68){\makebox(0,0)[lb]{\smash{{\SetFigFont{10}{12.0}{rm}$\Theta^{+}$}}}}
\put(3294,23){\makebox(0,0)[lb]{\smash{{\SetFigFont{10}{12.0}{rm}$\Theta^{-}$}}}}
\put(2691,1235){\makebox(0,0)[lb]{\smash{{\SetFigFont{10}{12.0}{rm}$x_{1}$}}}}
\put(3415,1178){\makebox(0,0)[lb]{\smash{{\SetFigFont{10}{12.0}{rm}$x_{2}$}}}}
\put(4762,703){\makebox(0,0)[lb]{\smash{{\SetFigFont{10}{12.0}{rm}$w$}}}}
\put(2344,2781){\makebox(0,0)[lb]{\smash{{\SetFigFont{10}{12.0}{rm}$w$}}}}
\put(2171,3179){\makebox(0,0)[lb]{\smash{{\SetFigFont{10}{12.0}{rm}$z$}}}}
\put(5147,2949){\makebox(0,0)[lb]{\smash{{\SetFigFont{10}{12.0}{rm}$0$}}}}
\includegraphics{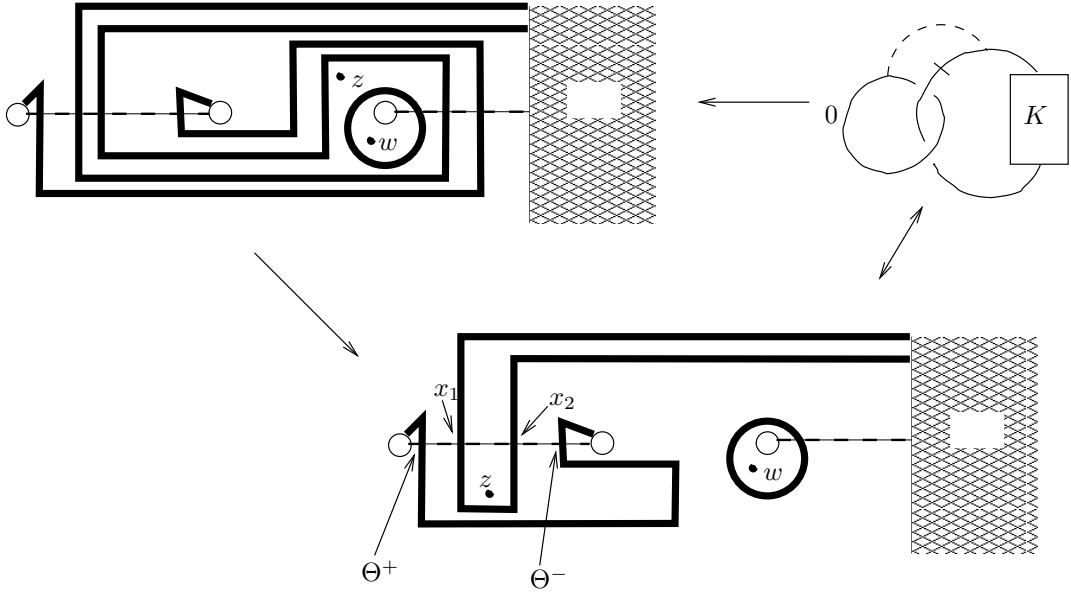}
\end{picture}%
\end{center}
\caption{$0$-surgery on an unknot linked with $K$. The $d$-base follows
the dashed ard without twisting. The hatched part of the diagram 
includes the effect of $K$ on the $\al$ (heavy) and $\be$ (dashed) 
curves. The portion shown in detail describes the linking and the 
surgery. Note that placing a meridian instead of the framing curve gives
a diagram for $K$. The intersection points never use $x_{1}$ or $x_{2}$, 
hence come from intersection points for the knot diagram. Destabilizing
the meridian for $K$ produces a diagram where the only holomorphic
discs are those from the complex for $K$ and for $S^{1} \times S^{2}$} 
\label{fig:zero}
\end{figure}

\ \\
\\
{\bf Example 1:}
In Figure \ref{fig:zero} we examine the homology of a knot in $S^{1} 
\times S^{2}$ which intersects $[S^2]$ precisely once. Regardless of the 
knot, $K$, we may find a diagram as in the figure. The intersection 
points, after handlesliding, giving generators in the chain complex are 
precisely $\Theta^{\pm} \times {\bf x}$, where ${\bf x}$ is a generator 
for the Knot Floer homology of $K$ in $S^{3}$. However, in the complex, 
only one of the two homotopy classes from $\Theta^{+}$ to $\Theta^{-}$ 
does not cross the marked point, $z$. It is straightforward to verify that 
$\widehat{\partial} \Theta^{+} \times {\bf x} =$ $ \Theta^{-} \times {\bf 
x} + \Theta^{+} \times \widehat{\partial}_{K}{\bf x}$. In each $Spin^{c}$ 
structure, the filtration index collapses; however, the differential as 
above also produces trivial homology in each $Spin^{c}$ structure.
 
\ \\
\\
{\bf Example 2:} Suppose we have a connected sum of standard genus $1$ 
diagrams for $S^{3}$ and $S^{1} \times S^{2}$. Then the canonical 
generator for the action of $H_{1}$ on the Heegaard-Floer group 
$\widehat{HF}(Y, \mathfrak{s}_{0})$, for the torsion $Spin^{c}$ structure, 
is represented by a 
single, closed intersection point: $\Theta^{+}$. We choose it as the 
basepoint for the complete set of paths, and choose a path in each 
component of the
connected sum to $\theta^{-}$. We can use these to connect $\Theta^{+}$ to 
any other intersection point. As there are only $2$ discs  in each 
component, which we assume do not cross $w$ and upon which the action of 
$\mathcal{A}$ takes one into the other, these each will evaluate to the 
same quantity on the $\overline{j}$ 
terms. However, since they
evaluate the same way on $\overline{j}$ we have seen that 
$[\Theta^{+}, 0, \overline{0}]$ and all other $\Theta$'s are closed in the 
complex $CF^{+}$. For the complex $\widehat{CF}$ with the additional 
marked points we will need to make further assumptions.

\subsection{Subtracting a Strand}

Suppose we remove a component strand of $S \subset (Y - B^{3})$ to obtain 
a new string link, 
$S'$. We may use  the complex and 
differential defined from a Heegaard diagram for $S$ by ignoring 
$z_k$. Without altering the complete set of paths or the point ${\bf 
x_{0}}$, the diagram without $z_{k}$ is a Heegaard diagram subordinate to 
$S'$, and the differential incorporates the 
same holomorphic discs. When $L_{k}$ does not algebraically 
intersect any homology class in $Y$, we can view the last coordinate, 
$j_{k}$, as
a filtration on the complex $\widehat{CF}(Y, S; \mathfrak{s})$ 
and use the associated spectral sequence with $E^{1}$ term $\oplus_{r} 
\widehat{HF}(Y, S; 
\mathfrak{s}, [j_{1}, \ldots, j_{k-1}], r)$ to calculate
$\widehat{HF}(Y, S'; \mathfrak{s}, [j_{1}, \ldots, j_{k-1}])$. It 
collapses in 
finitely many steps. 
\\
\\
Adding or subtracting an unknotted, unlinked, null-homologous component 
corresponds to adding or subtracting an
index which behaves like $i$.  The Heegaard diagram corresponds to 
stabilizing in the region containing $w$ and
placing a new point $z_{k+1}$ across the new $\al$, but in 
the same domain as $w$  

\subsection{Mirror String Links}

Let $S \subset Y- B^{3}$ be a string link in standard form, lying in the 
plane 
which defined the projection of our framed link diagram except in 
neighborhoods of the crossings. Let $\mathfrak{s}$ denote a $Spin^{c}$ 
structure on this manifold. 
Let $S'$ be the string link found through reflection in this plane, 
reflecting the framed components as well and switching the sign of their 
framing. Then $S'$ is the string link
induced by $S$ in $-Y$ under orientation reversal. Drawing the standard 
Heegaard diagram for $(Y, S)$, we may fix the $\be$-curves, and change the 
$\al$-curves for each crossing and framing to obtain a diagram for $(-Y, 
S')$. The meridians stay in their respective places;however, the marked 
points for each are reflected to the ``other side" of the strand to give 
marked points: $z'_{i}$. 
The intersection points in $\T_{\al} \cap \T_{\be}$ from the 
original diagram are in bijection with those of the new diagram. Each 
homotopy class $\phi$ is carried to a new homotopy class $\phi'$, but to 
join the same intersection points it must map in with reversed 
multiplicities. All this implies that we may calculate 
$\widehat{HF}_{\ast}( -Y, S'; \mathfrak{s}')$ by looking
at the intersection points for $(Y, S)$ and the differential for the 
complex using $-\Sigma$. As in Heegaard-Floer homology, this new complex 
calculates the co-homology $\widehat{CF}^{\ast}(Y, S; \mathfrak{s})$; 
there 
is thus an isomorphism $\widehat{HF}_{\ast}(Y, S; \mathfrak{s}) \ra 
\widehat{HF}^{\ast}( -Y, S'; \mathfrak{s}')$. 
(This isomorphism maps absolute degrees as $d \ra -d$ if $\mathfrak{s}$ is 
torsion). 
Using the same marked points, but the image of the basepoint and paths in 
the complete set of paths we find that $- \overline{\mathcal{F}}$ will be 
a filtration index for $S'$ when $\overline{\mathcal{F}}$ is for $S$ 
($\Lambda_{(-Y, S')} =$ $- \Lambda_{(Y, S)} = \Lambda_{(Y, S)})$. In 
particular, since each intersection point
has fixed image on each meridian, the boundary of a class $\phi$ must 
contain
whole multiples of the meridian and so $n_{z_{i}}(\phi) = - 
n_{z'_{i}}(\phi')$ since
the multiplicities reversed, but $n_{z'_{i}}(\phi') = n_{z_{i}}(\phi')$. 
In summary, there is an isomorphism (including the absolute grading when 
present):

$$
\widehat{HF}_{\ast}^{(-d)}(Y, S; \mathfrak{s}, [\overline{j}]) \ra 
\widehat{HF}^{\ast}_{(d)}(-Y , S'; \mathfrak{s}', [- \overline{j}])
$$

When $Y = S^{3}$ and $S$ is a normal string link, then $S'$ is the mirror 
image of $S$
found by switching all the crossings. The change in indices corresponds to 
the alteration
$t_{i} \ra t^{-1}_{i}$ in the Alexander polynomial (see section 
\ref{sec:Euler}).

\subsection{Three Operations on String Links}
\label{sec:stringcomp}

Given two string links, $S_{0}$ and $S_{1}$, in $Y_{0}$ and $Y_{1}$, there 
are 
three simple operations we can perform with two seperate string 
links, see Figure \ref{fig:Oper}. We always assume that the strands are 
oriented 
downwards. We will analyze the effects these operations have on the Floer 
homology. 

\setlength{\unitlength}{3200sp}%
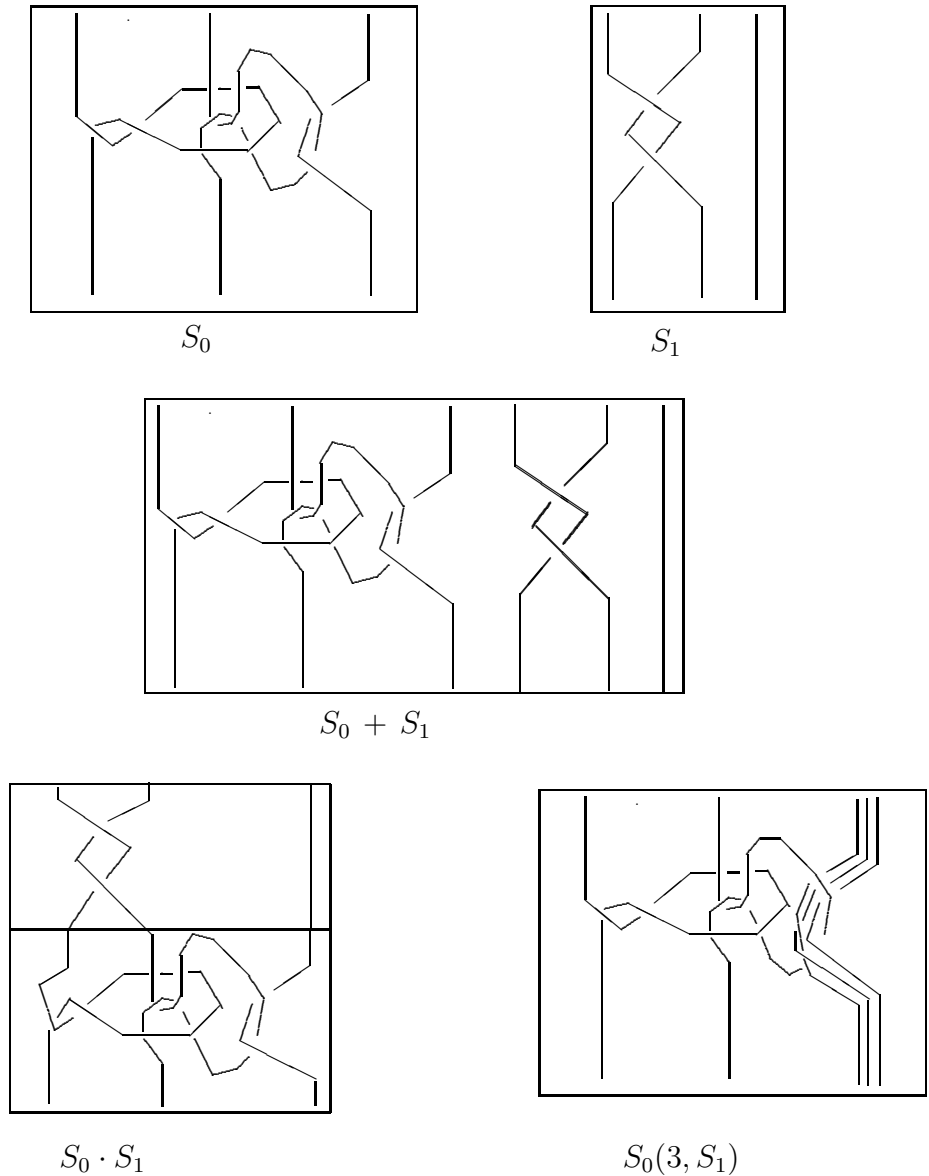
\begin{figure}
\begin{center}
\begin{picture}(7116,9049)(1108,-8231)
\thinlines
\put(2168,-4513){\framebox(4163,2275){}}
\put(1633,758){\line( 0,-1){799}}
\put(1633,-41){\line( 6,-5){282}}
\multiput(1915,-276)(7.05000,4.70000){21}{\makebox(1.6667,11.6667){\SetFigFont{5}{6}{rm}.}}
\put(2170,-70){\line( 6, 5){282}}
\put(2452,165){\line( 1, 0){188}}
\put(2745,168){\makebox(1.6667,11.6667){\SetFigFont{5}{6}{rm}.}}
\put(2951,180){\line( 1, 0){ 94}}
\multiput(3045,180)(4.34389,-7.23982){40}{\makebox(1.6667,11.6667){\SetFigFont{5}{6}{rm}.}}
\put(3195,-82){\line(-1,-1){235}}
\put(2958,-307){\line(-1, 0){517}}
\put(2441,-307){\line(-2, 1){470}}
\multiput(1971,-72)(-8.17391,-2.04348){24}{\makebox(1.6667,11.6667){\SetFigFont{5}{6}{rm}.}}
\put(1758,-207){\line( 0,-1){1222}}
\put(2670,755){\line( 0,-1){799}}
\multiput(2733,-120)(8.28378,1.38063){13}{\makebox(1.6667,11.6667){\SetFigFont{5}{6}{rm}.}}
\multiput(2833,-107)(4.32353,7.20588){15}{\makebox(1.6667,11.6667){\SetFigFont{5}{6}{rm}.}}
\put(2895, -7){\line( 0, 1){312}}
\multiput(2883,305)(4.47193,7.45321){23}{\makebox(1.6667,11.6667){\SetFigFont{5}{6}{rm}.}}
\multiput(2983,468)(8.49536,-2.12384){20}{\makebox(1.6667,11.6667){\SetFigFont{5}{6}{rm}.}}
\put(3145,430){\line( 1,-1){281}}
\multiput(3420,143)(4.81410,-7.22116){25}{\makebox(1.6667,11.6667){\SetFigFont{5}{6}{rm}.}}
\multiput(3533,-32)(-1.38235,-8.29412){35}{\makebox(1.6667,11.6667){\SetFigFont{5}{6}{rm}.}}
\multiput(2958,-232)(-3.91667,7.83333){13}{\makebox(1.6667,11.6667){\SetFigFont{5}{6}{rm}.}}
\multiput(2833,-45)(-8.77297,1.46216){11}{\makebox(1.6667,11.6667){\SetFigFont{5}{6}{rm}.}}
\multiput(2745,-32)(-6.92571,-5.19429){22}{\makebox(1.6667,11.6667){\SetFigFont{5}{6}{rm}.}}
\put(2598,-139){\line( 0,-1){141}}
\multiput(2608,-345)(5.03571,-6.71429){29}{\makebox(1.6667,11.6667){\SetFigFont{5}{6}{rm}.}}
\put(2749,-533){\line( 0,-1){893}}
\put(2733,168){\line( 1, 0){ 94}}
\put(3895,749){\line( 0,-1){517}}
\put(3895,232){\line(-3,-2){282}}
\multiput(3445,-73)(-2.68571,-8.05714){36}{\makebox(1.6667,11.6667){\SetFigFont{5}{6}{rm}.}}
\put(3351,-355){\line( 4,-3){564}}
\put(3915,-778){\line( 0,-1){658}}
\multiput(3420,-482)(-6.26667,-6.26667){16}{\makebox(1.6667,11.6667){\SetFigFont{5}{6}{rm}.}}
\multiput(3326,-576)(-8.17391,-2.04348){24}{\makebox(1.6667,11.6667){\SetFigFont{5}{6}{rm}.}}
\multiput(3138,-623)(-3.81081,7.62162){38}{\makebox(1.6667,11.6667){\SetFigFont{5}{6}{rm}.}}
\put(1283,-1557){\framebox(2987,2362){}}
\put(5751,756){\line( 0,-1){470}}
\put(5751,286){\line( 3,-2){564}}
\multiput(6315,-90)(-5.22222,-6.52778){37}{\makebox(1.6667,11.6667){\SetFigFont{5}{6}{rm}.}}
\put(6026,-432){\line(-5,-6){235}}
\put(5791,-714){\line( 0,-1){752}}
\put(6464,743){\line( 0,-1){282}}
\put(6464,461){\line(-1,-1){329}}
\multiput(6026, 
-7)(-5.03571,-6.71429){29}{\makebox(1.6667,11.6667){\SetFigFont{5}{6}{rm}.}}
\put(5914,-182){\line( 1,-1){564}}
\put(6478,-746){\line( 0,-1){705}}
\put(6901,743){\line( 0,-1){2209}}
\put(5626,-1557){\framebox(1488,2363){}}
\put(2668,-2351){\makebox(1.6667,11.6667){\SetFigFont{5}{6}{rm}.}}
\put(2268,-2286){\line( 0,-1){799}}
\put(2268,-3085){\line( 6,-5){282}}
\multiput(2550,-3320)(7.05000,4.70000){21}{\makebox(1.6667,11.6667){\SetFigFont{5}{6}{rm}.}}
\put(2805,-3114){\line( 6, 5){282}}
\put(3087,-2879){\line( 1, 0){188}}
\put(3380,-2876){\makebox(1.6667,11.6667){\SetFigFont{5}{6}{rm}.}}
\put(3586,-2864){\line( 1, 0){ 94}}
\multiput(3680,-2864)(4.34389,-7.23982){40}{\makebox(1.6667,11.6667){\SetFigFont{5}{6}{rm}.}}
\put(3830,-3126){\line(-1,-1){235}}
\put(3593,-3351){\line(-1, 0){517}}
\put(3076,-3351){\line(-2, 1){470}}
\multiput(2606,-3116)(-8.17391,-2.04348){24}{\makebox(1.6667,11.6667){\SetFigFont{5}{6}{rm}.}}
\put(2393,-3251){\line( 0,-1){1222}}
\put(3305,-2289){\line( 0,-1){799}}
\multiput(3368,-3164)(8.28378,1.38063){13}{\makebox(1.6667,11.6667){\SetFigFont{5}{6}{rm}.}}
\multiput(3468,-3151)(4.32353,7.20588){15}{\makebox(1.6667,11.6667){\SetFigFont{5}{6}{rm}.}}
\put(3530,-3051){\line( 0, 1){312}}
\multiput(3518,-2739)(4.47193,7.45321){23}{\makebox(1.6667,11.6667){\SetFigFont{5}{6}{rm}.}}
\multiput(3618,-2576)(8.49536,-2.12384){20}{\makebox(1.6667,11.6667){\SetFigFont{5}{6}{rm}.}}
\put(3780,-2614){\line( 1,-1){281}}
\multiput(4055,-2901)(4.81410,-7.22116){25}{\makebox(1.6667,11.6667){\SetFigFont{5}{6}{rm}.}}
\multiput(4168,-3076)(-1.38235,-8.29412){35}{\makebox(1.6667,11.6667){\SetFigFont{5}{6}{rm}.}}
\multiput(3593,-3276)(-3.91667,7.83333){13}{\makebox(1.6667,11.6667){\SetFigFont{5}{6}{rm}.}}
\multiput(3468,-3089)(-8.77297,1.46216){11}{\makebox(1.6667,11.6667){\SetFigFont{5}{6}{rm}.}}
\multiput(3380,-3076)(-6.92571,-5.19429){22}{\makebox(1.6667,11.6667){\SetFigFont{5}{6}{rm}.}}
\put(3233,-3183){\line( 0,-1){141}}
\multiput(3243,-3389)(5.03571,-6.71429){29}{\makebox(1.6667,11.6667){\SetFigFont{5}{6}{rm}.}}
\put(3384,-3577){\line( 0,-1){893}}
\put(3368,-2876){\line( 1, 0){ 94}}
\put(4530,-2295){\line( 0,-1){517}}
\put(4530,-2812){\line(-3,-2){282}}
\multiput(4080,-3117)(-2.68571,-8.05714){36}{\makebox(1.6667,11.6667){\SetFigFont{5}{6}{rm}.}}
\put(3986,-3399){\line( 4,-3){564}}
\put(4550,-3822){\line( 0,-1){658}}
\multiput(4055,-3526)(-6.26667,-6.26667){16}{\makebox(1.6667,11.6667){\SetFigFont{5}{6}{rm}.}}
\multiput(3961,-3620)(-8.17391,-2.04348){24}{\makebox(1.6667,11.6667){\SetFigFont{5}{6}{rm}.}}
\multiput(3773,-3667)(-3.81081,7.62162){38}{\makebox(1.6667,11.6667){\SetFigFont{5}{6}{rm}.}}
\put(5030,-2275){\line( 0,-1){470}}
\put(5030,-2745){\line( 3,-2){564}}
\multiput(5594,-3121)(-5.22222,-6.52778){37}{\makebox(1.6667,11.6667){\SetFigFont{5}{6}{rm}.}}
\put(5305,-3463){\line(-5,-6){235}}
\put(5070,-3745){\line( 0,-1){752}}
\put(5743,-2288){\line( 0,-1){282}}
\put(5743,-2570){\line(-1,-1){329}}
\multiput(5305,-3038)(-5.03571,-6.71429){29}{\makebox(1.6667,11.6667){\SetFigFont{5}{6}{rm}.}}
\put(5193,-3213){\line( 1,-1){564}}
\put(5757,-3777){\line( 0,-1){705}}
\put(6180,-2288){\line( 0,-1){2209}}
\put(5030,-2288){\line( 0,-1){470}}
\put(5030,-2758){\line( 3,-2){564}}
\multiput(5594,-3134)(-5.22222,-6.52778){37}{\makebox(1.6667,11.6667){\SetFigFont{5}{6}{rm}.}}
\put(5305,-3476){\line(-5,-6){235}}
\put(5070,-3758){\line( 0,-1){752}}
\put(5743,-2301){\line( 0,-1){282}}
\put(5743,-2583){\line(-1,-1){329}}
\multiput(5305,-3051)(-5.03571,-6.71429){29}{\makebox(1.6667,11.6667){\SetFigFont{5}{6}{rm}.}}
\put(5193,-3226){\line( 1,-1){564}}
\put(5757,-3790){\line( 0,-1){705}}
\put(6180,-2301){\line( 0,-1){2209}}
\put(2033,693){\makebox(1.6667,11.6667){\SetFigFont{5}{6}{rm}.}}
\put(1567,-6361){\line( 0,-1){282}}
\put(1567,-6643){\line(-5,-3){221.912}}
\multiput(1345,-6776)(2.65333,-7.96000){46}{\makebox(1.6667,11.6667){\SetFigFont{5}{6}{rm}.}}
\multiput(1465,-7134)(7.05000,4.70000){21}{\makebox(1.6667,11.6667){\SetFigFont{5}{6}{rm}.}}
\put(1720,-6928){\line( 6, 5){282}}
\put(2002,-6693){\line( 1, 0){188}}
\put(2295,-6690){\makebox(1.6667,11.6667){\SetFigFont{5}{6}{rm}.}}
\put(2501,-6678){\line( 1, 0){ 94}}
\multiput(2595,-6678)(4.34389,-7.23982){40}{\makebox(1.6667,11.6667){\SetFigFont{5}{6}{rm}.}}
\put(2745,-6940){\line(-1,-1){235}}
\put(2508,-7165){\line(-1, 0){517}}
\put(1991,-7165){\line(-3, 2){410.308}}
\multiput(1583,-6888)(-4.70000,-7.05000){21}{\makebox(1.6667,11.6667){\SetFigFont{5}{6}{rm}.}}
\put(2220,-6385){\line( 0,-1){517}}
\multiput(2283,-6978)(8.28378,1.38063){13}{\makebox(1.6667,11.6667){\SetFigFont{5}{6}{rm}.}}
\multiput(2383,-6965)(4.32353,7.20588){15}{\makebox(1.6667,11.6667){\SetFigFont{5}{6}{rm}.}}
\put(2445,-6865){\line( 0, 1){312}}
\multiput(2433,-6553)(4.47193,7.45321){23}{\makebox(1.6667,11.6667){\SetFigFont{5}{6}{rm}.}}
\multiput(2533,-6390)(8.49536,-2.12384){20}{\makebox(1.6667,11.6667){\SetFigFont{5}{6}{rm}.}}
\put(2695,-6428){\line( 1,-1){281}}
\multiput(2970,-6715)(4.81410,-7.22116){25}{\makebox(1.6667,11.6667){\SetFigFont{5}{6}{rm}.}}
\multiput(3083,-6890)(-1.38235,-8.29412){35}{\makebox(1.6667,11.6667){\SetFigFont{5}{6}{rm}.}}
\multiput(2508,-7090)(-3.91667,7.83333){13}{\makebox(1.6667,11.6667){\SetFigFont{5}{6}{rm}.}}
\multiput(2383,-6903)(-8.77297,1.46216){11}{\makebox(1.6667,11.6667){\SetFigFont{5}{6}{rm}.}}
\multiput(2295,-6890)(-6.92571,-5.19429){22}{\makebox(1.6667,11.6667){\SetFigFont{5}{6}{rm}.}}
\put(2148,-6997){\line( 0,-1){141}}
\multiput(2158,-7203)(5.03571,-6.71429){29}{\makebox(1.6667,11.6667){\SetFigFont{5}{6}{rm}.}}
\put(2299,-7391){\line( 0,-1){329}}
\put(2283,-6690){\line( 1, 0){ 94}}
\put(3445,-6344){\line( 0,-1){282}}
\put(3445,-6626){\line(-3,-2){282}}
\multiput(2995,-6931)(-2.68571,-8.05714){36}{\makebox(1.6667,11.6667){\SetFigFont{5}{6}{rm}.}}
\put(2901,-7213){\line( 2,-1){590.400}}
\put(3483,-7525){\line( 0,-1){188}}
\multiput(2970,-7340)(-6.26667,-6.26667){16}{\makebox(1.6667,11.6667){\SetFigFont{5}{6}{rm}.}}
\multiput(2876,-7434)(-8.17391,-2.04348){24}{\makebox(1.6667,11.6667){\SetFigFont{5}{6}{rm}.}}
\multiput(2688,-7481)(-3.81081,7.62162){38}{\makebox(1.6667,11.6667){\SetFigFont{5}{6}{rm}.}}
\put(3449,-5227){\line( 0,-1){1128}}
\put(1420,-7131){\line( 0,-1){564}}
\put(1490,-5248){\line( 0,-1){ 94}}
\put(1490,-5342){\line( 3,-2){564}}
\multiput(2054,-5718)(-5.22222,-6.52778){37}{\makebox(1.6667,11.6667){\SetFigFont{5}{6}{rm}.}}
\multiput(1765,-6060)(-4.70000,-7.05000){41}{\makebox(1.6667,11.6667){\SetFigFont{5}{6}{rm}.}}
\put(2196,-5347){\line(-2,-1){317.200}}
\multiput(1765,-5635)(-5.03571,-6.71429){29}{\makebox(1.6667,11.6667){\SetFigFont{5}{6}{rm}.}}
\put(1640,-5809){\line( 1,-1){579}}
\put(1120,-7763){\framebox(2475,1413){}}
\put(1120,-6350){\framebox(2475,1125){}}
\put(2195,-5350){\line( 0, 1){141}}
\put(5975,-5382){\makebox(1.6667,11.6667){\SetFigFont{5}{6}{rm}.}}
\put(5575,-5317){\line( 0,-1){799}}
\put(5575,-6116){\line( 6,-5){282}}
\multiput(5857,-6351)(7.05000,4.70000){21}{\makebox(1.6667,11.6667){\SetFigFont{5}{6}{rm}.}}
\put(6112,-6145){\line( 6, 5){282}}
\put(6394,-5910){\line( 1, 0){188}}
\put(6687,-5907){\makebox(1.6667,11.6667){\SetFigFont{5}{6}{rm}.}}
\put(6893,-5895){\line( 1, 0){ 94}}
\multiput(6987,-5895)(4.34389,-7.23982){40}{\makebox(1.6667,11.6667){\SetFigFont{5}{6}{rm}.}}
\put(7137,-6157){\line(-1,-1){235}}
\put(6900,-6382){\line(-1, 0){517}}
\put(6383,-6382){\line(-2, 1){470}}
\multiput(5913,-6147)(-8.17391,-2.04348){24}{\makebox(1.6667,11.6667){\SetFigFont{5}{6}{rm}.}}
\put(5700,-6282){\line( 0,-1){1222}}
\put(6612,-5320){\line( 0,-1){799}}
\multiput(6675,-6195)(8.28378,1.38063){13}{\makebox(1.6667,11.6667){\SetFigFont{5}{6}{rm}.}}
\multiput(6775,-6182)(4.32353,7.20588){15}{\makebox(1.6667,11.6667){\SetFigFont{5}{6}{rm}.}}
\put(6837,-6082){\line( 0, 1){312}}
\multiput(6825,-5770)(5.59965,6.71958){20}{\makebox(1.6667,11.6667){\SetFigFont{5}{6}{rm}.}}
\put(6931,-5642){\line( 1, 0){156}}
\put(7087,-5645){\line( 1,-1){281}}
\multiput(7362,-5932)(4.81410,-7.22116){25}{\makebox(1.6667,11.6667){\SetFigFont{5}{6}{rm}.}}
\multiput(7475,-6107)(-1.38235,-8.29412){35}{\makebox(1.6667,11.6667){\SetFigFont{5}{6}{rm}.}}
\multiput(6900,-6307)(-3.91667,7.83333){13}{\makebox(1.6667,11.6667){\SetFigFont{5}{6}{rm}.}}
\multiput(6775,-6120)(-8.77297,1.46216){11}{\makebox(1.6667,11.6667){\SetFigFont{5}{6}{rm}.}}
\multiput(6687,-6107)(-6.92571,-5.19429){22}{\makebox(1.6667,11.6667){\SetFigFont{5}{6}{rm}.}}
\put(6540,-6214){\line( 0,-1){141}}
\multiput(6550,-6420)(5.03571,-6.71429){29}{\makebox(1.6667,11.6667){\SetFigFont{5}{6}{rm}.}}
\put(6691,-6608){\line( 0,-1){893}}
\put(6675,-5907){\line( 1, 0){ 94}}
\put(7837,-5326){\line( 0,-1){517}}
\put(7837,-5843){\line(-3,-2){282}}
\multiput(7387,-6148)(-2.68571,-8.05714){36}{\makebox(1.6667,11.6667){\SetFigFont{5}{6}{rm}.}}
\put(7293,-6430){\line( 4,-3){564}}
\put(7857,-6853){\line( 0,-1){705}}
\multiput(7251,-6656)(-7.83333,-3.91667){13}{\makebox(1.6667,11.6667){\SetFigFont{5}{6}{rm}.}}
\multiput(7157,-6703)(-6.67570,5.56308){24}{\makebox(1.6667,11.6667){\SetFigFont{5}{6}{rm}.}}
\multiput(7006,-6572)(-3.15172,7.87931){21}{\makebox(1.6667,11.6667){\SetFigFont{5}{6}{rm}.}}
\put(5225,-7632){\framebox(2987,2362){}}
\put(7758,-5332){\line( 0,-1){470}}
\put(7758,-5802){\line(-3,-2){282}}
\put(7682,-5342){\line( 0,-1){423}}
\put(7682,-5765){\line(-5,-3){235}}
\multiput(7354,-6046)(-3.13333,-7.83333){31}{\makebox(1.6667,11.6667){\SetFigFont{5}{6}{rm}.}}
\multiput(7307,-5999)(-3.13333,-7.83333){31}{\makebox(1.6667,11.6667){\SetFigFont{5}{6}{rm}.}}
\multiput(7213,-6234)(1.38235,-8.29412){35}{\makebox(1.6667,11.6667){\SetFigFont{5}{6}{rm}.}}
\put(7204,-6365){\line( 0,-1){141}}
\put(7204,-6506){\line( 3,-2){572.077}}
\put(7767,-6901){\line( 0,-1){658}}
\multiput(7269,-6581)(3.23678,-8.09195){16}{\makebox(1.6667,11.6667){\SetFigFont{5}{6}{rm}.}}
\put(7316,-6703){\line( 5,-3){380.147}}
\put(7692,-6938){\line( 0,-1){611}}
\put(2439,-1836){\makebox(0,0)[lb]{\smash{{\SetFigFont{12}{14.4}{rm}$S_{0}$}}}}
\put(6076,-1861){\makebox(0,0)[lb]{\smash{{\SetFigFont{12}{14.4}{rm}$S_{1}$}}}}
\put(3514,-4824){\makebox(0,0)[lb]{\smash{{\SetFigFont{12}{14.4}{rm}$S_{0}\, 
+\, S_{1}$}}}}
\put(1501,-8186){\makebox(0,0)[lb]{\smash{{\SetFigFont{12}{14.4}{rm}$S_{0} 
\cdot S_{1}$}}}}
\put(5864,-8186){\makebox(0,0)[lb]{\smash{{\SetFigFont{12}{14.4}{rm}$S_{0} 
(3, S_{1})$}}}}
\end{picture}%
\end{center}
\caption{Three Simple Operations on String Links}
\label{fig:Oper}
\end{figure}
\setlength{\unitlength}{3947sp}%

\subsubsection{$S_{1} + S_{2}$}

We assume, for the first, that we have $(Y_{0}, \Gamma_{0})$ and $(Y_{1}, 
\Gamma_{1})$ presented as string links in framed surgery diagrams in 
$D^{2} \times I$. 
We assume that these have been put in standard form. This means we 
arrange that all the meridians, at the bottom of each 
diagram, intersect at most two $\be$'s. However, there is only one 
choice possible for every $g$-tuple of intersection points due to the 
presence of $U$. Amalgamating the second string link does not affect this 
property for the meridians. Alternately, we can wind in the complement of 
the star in $D^{2} \times \{0\}$ and the amalgamation region as their 
union is contractible in $\Si$.

Topologically, the amalgamation is a connect sum of $Y_{0}$ and $Y_{1}$ 
where the sums occur for balls removed outside the region depicted as 
$D^{2} \times I$. Thus for two $Spin^{c}$ structures, $\mathfrak{s}_{0}$ 
and $\mathfrak{s}_{1}$, there is a unique $\mathfrak{s} = \mathfrak{s}_{0} 
\# \mathfrak{s}_{1}$ on the amalgamated picture. Furthermore, $H_{2} \cong 
H_{2}(Y_{0}; \Z) \oplus H_{2}(Y_{1}; \Z)$. If the first string link has 
$k_{0}$ strands and the second $k_{1}$ strands, the amalgamation has 
filtration index taking values in $\Z^{k_{0}}/\Lambda_{Y_{0}} \oplus 
\Z^{k_{1}}/\Lambda_{Y_{1}}$. 

Counting $\al$'s and $\be$'s from the portion if $Y$ coming from $Y_{0}$ 
demonstrates that for an intersection point we must have an $\al$ from 
$Y_{0}$ pairing with a $\be$ from $Y_{0}$ and likewise for $Y_{1}$. Even 
if some $\be$'s extend from the $Y_{0}$ region to the $Y_{1}$ region 
(which we can avoid if we like), this remains true. In particular, the 
diagrams drawn from projections in $D^{2} \times I$ will have one such 
$\be$. Therefore, the 
generators for the new chain group are the product of generators for the 
two previous groups: as groups $\widehat{CF}(Y, S_{0} + S_{1}; 
\mathfrak{s}) \cong$ $\widehat{CF}(Y_{0}, S_{0}; \mathfrak{s}_{0}) 
\otimes \widehat{CF}(Y_{1}, S_{1}; \mathfrak{s}_{1})$. We may choose 
filtration indices for both links, choosing basepoints and complete sets 
of paths. The amalgamation will have $(\overline{\mathcal{F}}_{0}, 
\overline{\mathcal{F}}_{1})$ as a filtration index for the complete set of 
paths found by using the product of the two basepoints and the paths from 
complete sets for $Y_{0}$ and $Y_{1}$. That the domain containing $w$ 
corresponds to the outer boundary of $D^{2} \times I - (S_{0} + 
S_{1})$, and that the same is true for each string link individually, 
ensures that the paths in each complete set need not be alterred. 
Furthermore, this region separates the domains for classes, $\phi$, used 
in the  differentials in the two complexes; thus 
$\widehat{\partial}_{S_{0} + S_{1}} = 
\big(\widehat{\partial}_{S_{0}}\otimes I \big) \oplus 
\big( I\otimes\widehat{\partial}_{S_{1}}\big)$. We 
have verified that 

$$
\widehat{HF}(Y, S_{1} + S_{2}; \mathfrak{s}, [\overline{j}_{0}] 
\oplus [\overline{j}_{1}]) \cong H_{\ast}(\widehat{CF}(Y_{0}, S_{0}; 
\mathfrak{s}_{0}, [\overline{j}_{0}]) \otimes \widehat{CF}(Y_{1}, 
S_{1}; \mathfrak{s}_{1}, [\overline{j}_{2}]))
$$
up to gradings. The grading calculation follows as for connected sums 
\cite{3Man}. In particular, for two torsion $Spin^{c}$ structures on 
$Y_{0}$ and $Y_{1}$, the absolute grading satisfies $gr({\bf x} \otimes 
{\bf y}) =$ $gr_{S_0}({\bf x}) + gr_{S_{1}}({\bf y})$, which is all we 
require for string links in $S^{3}$. 
We may also establish this relation by using the Maslov index calculation 
for $Y_{0}$ or $Y_{1}$, presented as surgery on a link in $S^{3}$, found 
in the absolute gradings section of \cite{Four}. Our assumptions include 
such presentations for $Y_{0}$ and $Y_{1}$ and the connect sum provides 
one for $Y$. We may use triangles with $n_{w} = 0$ for $Y_{0}$ and 
$Y_{1}$. This allows use to use a product triangle in the calculation for 
$Y$. The first Chern class for the associated $Spin^{c}$ structure will be 
the sum of those for $Y_{0}$ and $Y_{1}$. Since the intersection form 
splits and the Euler characteristics and signatures of the cobordisms add, 
we see that the gradings for the complexes for $Y_{0}$ and $Y_{1}$ add to 
give that for the complex on $Y$. 

\subsubsection{$S_{1}\cdot S_{2}$}

The second operation is the composition of pure string links, the analog 
of composition for braids. The torsion of the composite is the product of 
the torsions of the two factors, \cite{Livi}. We may prove the analogous 
result  the homologies of the string links. For $n=1$ stranded string 
links composition corresponds to the connect sum of knots. 

Again we will work with $(Y_{0}, S_{0})$ and $(Y_{1}, S_{1})$ with the 
assumption that $S_{0}$ and $S_{1}$ have the same number, $n$, of strands 
going from top to bottom. In addition, we require each component of 
the string links to have one boundary on the top and one boundary on the 
bottom of the $D^{2} \times I$ region in their respective manifolds. No 
closed 
component is formed by the stacking operation.

We prove that 

$$
\widehat{HF}(Y, S_{0} \cdot S_{1};\mathfrak{s}, [\overline{k}]) = 
\bigoplus_{[\overline{k}_{0}] + [\overline{k}_{1}] = [\overline{k}]\, 
\mathrm{mod\,}\Lambda} H_{\ast}( \widehat{CF}(Y_{0}, S_{0}; 
\mathfrak{s}_{0}, [\overline{k}_{0}]) \otimes \widehat{CF}(Y_{1}, S_{1}; 
\mathfrak{s}_{1}, [\overline{k}_{2}]))
$$

\begin{figure}
\begin{center}
\begin{picture}(5003,4737)
\includegraphics{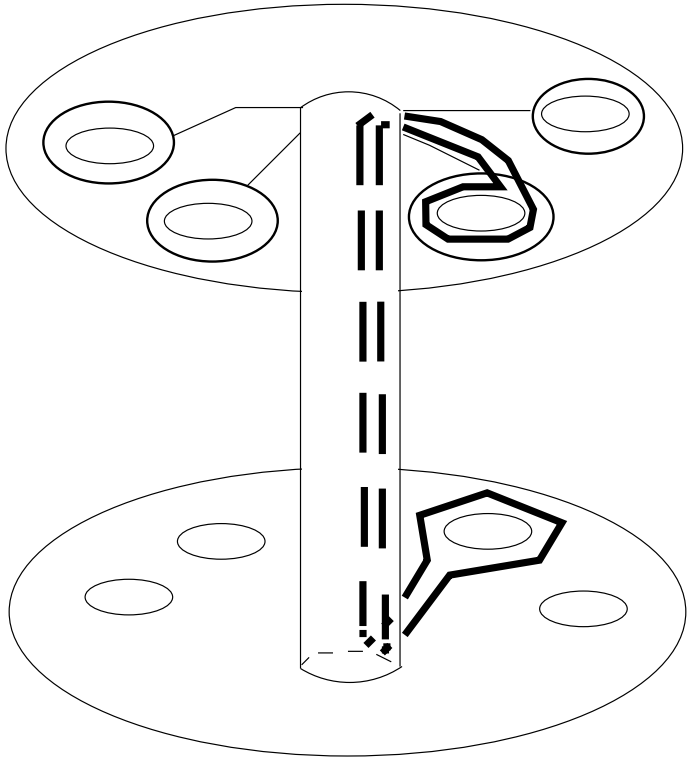}
\end{picture}%
\end{center}
\caption{Heegaard Diagram for the Composite of String Links}
\label{fig:compI}
\end{figure}

Let $\Si_{\al_{0}\be_{0}}$ be a weakly admissible Heegaard diagram for 
$(Y_{0}, S_{0})$ with marked points $w, z_{1}, \ldots, z_{n}$ and 
$\Si_{\al_{1}\be_{1}}$ be a weakly admissible Heegaard diagram for 
$(Y_{0}, 
S_{0})$ with marked points $w', z'_{1}, \ldots, z'_{n}$. In each case, we 
use a diagram in standard form. Thus, the region containing $w$ (or $w'$)
includes all of $\partial (D^{2} \times I)$ minus the strands. 
The diagram $\al_{0}\al_{1}, \be_{0}\be_{1}$ is formed as in 
Figure \ref{fig:compI}, by joining $\Si_{0}$ and $\Si_{1}$ with a tube. 
The end of the tube in $S_{1}$ should occur close to $w$. We have depicted 
the ``star'' for $S_{1}$ by the thin lines emanating from the tube. As the 
two $w$'s now occur in the same domain, we will consider only one $w$ 
point. We let $\al'_{i}$ be the small Hamiltonian isotopes of the 
$\al_{i}$ curves, except at the meridians. At the meridians we choose 
curves which traverse the tubeand loop around the $i^{th}$ strand in each 
diagram, as depicted
for the thick curve in Figure 
\ref{fig:compI}. 

\begin{figure}
\begin{center}
\begin{picture}(3660,2664)
\put(1941,239){\makebox(0,0)[lb]{\smash{{\SetFigFont{12}{14.4}{rm}$x_{1}$}}}}
\put(2298,2155){\makebox(0,0)[lb]{\smash{{\SetFigFont{12}{14.4}{rm}$x_{j-1}$}}}}
\put(2758,1103){\makebox(0,0)[lb]{\smash{{\SetFigFont{12}{14.4}{rm}$x_{2}$}}}}
\put(1528,2380){\makebox(0,0)[lb]{\smash{{\SetFigFont{12}{14.4}{rm}$x_{j}$}}}}
\put(541,1619){\makebox(0,0)[lb]{\smash{{\SetFigFont{12}{14.4}{rm}$\Theta^{+}$}}}}
\put(541,1018){\makebox(0,0)[lb]{\smash{{\SetFigFont{12}{14.4}{rm}$\Theta^{-}$}}}}
\put(1251,658){\makebox(0,0)[lb]{\smash{{\SetFigFont{12}{14.4}{rm}$z'_{i}$}}}}
\includegraphics{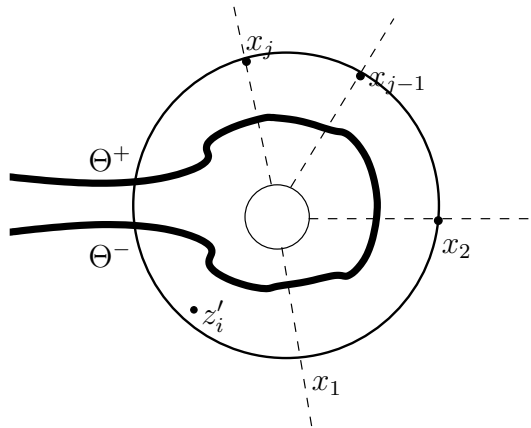}
\end{picture}%
\end{center}
\caption{The region near the $i^{th}$ meridian in the top diagram. The 
dashed curves
are $\be$'s and the thickest curve is the $\ga$ curve which replaces the 
circular
meridian}
\label{fig:hhhh}
\end{figure} 

We will analyze the cobordism generated by the triple $\al_{0}\al_{1}$, 
$\be_{0}\be_{1}$, and $\al_{0}\al_{1}'$, where the occurrence of repeated 
sets of curves in a diagram indicates using Hamiltonian isotopes in the 
standard way.  We show a closeup of the meridians in Figure 
\ref{fig:hhhh}. The thick curves should follow the ``star'' except that 
the domain containing $z_{i}'$ must also abut $\Theta^{-}$. 

Each of the ends of this cobordism has $2\,n$ marked points, so we will 
have filtration indices taking values in $\Z^{2n}$, modulo some lattice. 
First, we describe the various boundary components. 
\\
\\
\noindent
{\bf Boundary I:} $\al_{0}\al_{1}, \be_{0}\be_{1}$ \\ 

Topologically, this is a connect sum for $Y_{0}$ and $Y_{1}$ whose 
Heegaard diagram is drawn in the standard way. Each $Spin^{c}$ structure 
is therfore of the form $\mathfrak{s}_{0}$, $\mathfrak{s}_{1}$. Inclusion 
of the marked points puts us into the previous construction: amalgamation. 
Thus, the generators of the complex for $\mathfrak{s}$ are products of 
generators from $\mathfrak{s}_{0}$ and $\mathfrak{s}_{1}$, which we denote 
${\bf x} \otimes {\bf y}$. We may take the products of our basepoints to 
be the basepoint for a complete set of paths formed by the product of 
paths on the two separate diagrams (since $w$ excludes the tube joining 
the two). As in the amalgamation case, we find that $\Lambda_{I} \equiv$ 
$\Lambda_{(Y_{0}, S_{0})} \oplus \Lambda_{(Y_{1}, S_{1})}$ with filtration 
index $(\overline{\mathcal{F}}_{0}, \overline{\mathcal{F}}_{1})$. We use 
the previous result to identify 

$$
\widehat{CF}(Y_{0} \# Y_{1}, S_{0} + S_{1}; \mathfrak{s}, 
[\overline{j}_{0}] \oplus [\overline{j}_{1}]) \cong \widehat{CF}(Y_{0}, 
S_{0}; \mathfrak{s}_{0}, [\overline{j}_{0}]) \otimes \widehat{CF}(Y_{1}, 
S_{1}; \mathfrak{s}_{1}, [\overline{j}_{1}]) 
$$
\noindent
{\bf Boundary II:} $\al_{0}\al_{1}, \al_{0}\al'_{1}$ \\

Topologically, this boundary is a $g_{0} + g_{1}$ fold connect sums of 
$S^{1} \times S^{2}$'s, which may be seen by isotoping the new $\al'_{1}$ 
components down the strands of $S_{0}$, across the meridians found there, 
and back up the strands. However, with the additional marked points, it is 
unclear if the candidate for $\Theta^{+}_{std}$ is closed in this diagram. 
We can choose $\Theta^{+}_{std}$ as our basepoint, and use products of 
topological discs in $\Si$ from each connect sum component as our complete 
set of paths. If we denote by $e_{i}$, the $i^{th}$ basis vector in 
$\Z^{2n}$, the lattice for this component will be $\Lambda_{II} \equiv 
{\mathrm Span}\,\{\,e_{i} - e_{n+i}\}$. We now argue that 
$\Theta^{+}_{std}$ is indeed closed for the differential missing {\em all} 
marked points. 

There are precisely $2^{g_{0} + g_{1}}$ intersection possible for this 
diagram. We may use our complete set of topological discs to see that 
$\Theta^{+}_{std}$ has maximal grading. By the Heegaard-Floer homology 
theory, we must have that the {\em generator} $\Theta^{+}_{std}$ is 
closed for the differential 
only missing $w$. In particular, a holomorphic disc contributing to this 
differential cancels with some other holomorphic disc. Suppose we have two 
such homotopy classes of discs, $\phi$ and $\phi'$. Splicing the inverse 
of one to the other, $\phi^{-1} \ast \phi'$ must produce a periodic 
domain. 
This periodic domain must evaluate to an element of $\Lambda_{II}$ under 
the application of $n_{\overline{z}}$. However, $\sum 
n_{z_{i}}(\mathcal{P}) = 0$ for every periodic domain. Since classes with 
holomorphic representatives must have non-negative multiplicities, it must 
be the case that when $n_{\overline{z}}(\phi) = 0$ so too 
$n_{\overline{z}}(\phi') = 0$. As the differential missing all marked 
points arises from a subset of the moduli spaces in the Heegaard-Floer 
differential and adding the extra marked points does not eliminate one 
disk in a cancelling pair without eliminating the other, so 
$\Theta^{+}_{std}$ is still closed
\\
\\
\noindent
{\bf Boundary III:} $\al_{0}\al'_{1}, \be_{0}\be_{1}$ \\

This diagram represents the result of composition. Topologically, each of 
the new $\al'$ curves may be slid down a component of $S_{0}$ until it 
reaches a meridian. After sliding across the meridian, and back up the 
diagram, we have the connect sum of the diagrams for $Y_{0}$ and $Y_{1}$. 
Again, each $Spin^{c}$ structure on $Y$ is the sum of structures, 
$\mathfrak{s}_{i}$, from $Y_{0}$ and $Y_{1}$. Additionally, $H_{2}(Y; \Z) 
\cong H_{2}(Y_{0}; \Z) \oplus H_{2}(Y_{1}; \Z)$. However, the lattices now 
combine as $\Lambda_{III} \equiv \big(\Lambda_{Y_{0}} + 
\Lambda_{Y_{1}}\big) \oplus 
\overline{0}$, the span of the two original lattices, and $z_{i}'$ is in 
the same 
domain as $w$. We required that our 
original diagrams be weakly admissible for our $Spin^{c}$ structures. We 
will see below how to extend periodic domains so that they continue to 
have positive and negative multiplicities in the diagram for $Y$. Thus, 
the new diagram  will be weakly admissible.

In the diagram for $S_{0}$ there are $g_{0}$ $\al$ curves and $g_{0}$ 
$\be$ curves. As we have not changed these, all the $\be$ curves that 
intersect an alterred $\al'$ in the ``lower'' diagram must pair with an 
$\al$ curve from the lower diagram when describing a generator. Hence, the 
new $\al$'s must pair with 
$\be$'s from the upper diagram. These intersections have precisely the 
same form as intersections with the meridians in the original diagram for 
$Y_{1}$. This allows us to establish a one-to-one correspondence between 
the product of generators for $Y_{0}$ and $Y_{1}$ and those of $Y$. The 
chain complex, as an abelian group, is the product of the original chain 
complexes. These generators we denote ${\bf x\, y}$.

If we have basepoints for a complete set of paths on $Y_{0}$ and $Y_{1}$ 
for our $Spin^{c}$ structures, we may choose as our basepoint on $Y$, the 
generator corresponding to the product of these basepoints. To specify the 
complete set of paths, consider a class $\phi$ with $n_{w}(\phi) = 0$ in 
either of the original diagrams. If $\phi$ is in the lower diagram, we may 
use $\phi$ in the diagram for $Y$ as the condition on $n_{w}$ implies that 
the domain of the disc does not extend into the upper diagram: a small 
region at the top of each strand lies in the domain containing $w$ in the 
diagram for $Y_{0}$. In $\phi$ occurs in the upper diagram, its domain may 
cross $n_{z'_{i}}$ and include copies of the meridians in its boundary. In 
the diagram for $Y$, we may extend this disc by following the strand down 
to the meridian from $S_{0}$. At crossings, the disc gains a boundary 
component and an intersection point, or a copy of a framed component. 
However, it will not cross $w$, and it crosses  $z_{i}$ the same number of 
times as $\phi$ crossed $z'_{i}$. Periodic domains will continue to have 
the same multiplicities in the regions coming from their respective 
diagrams. 

The class $\phi_{\bf x}^{0}$ may be used to join ${\bf x}_{0} \, {\bf 
y}$ to ${\bf x} \, {\bf y}$ for any ${\bf y}$ coming from the diagram 
for $Y_{1}$. Likewise, $\phi_{\bf y}^{1}$, when extended, may be used to 
join ${\bf x} \, {\bf y}_{0}$ to ${\bf x} \, {\bf y}$. We use 
these for the complete set of paths, and extend to get $\phi_{{\bf 
x}\, {\bf y}}$ by composing $\phi_{\bf y}^{1} \ast \phi_{\bf x}^{0}$. 
In particular, the filtration value for ${\bf x}_{0} \, {\bf y}$ is 
$-n_{\overline{z'}}(\phi_{\bf y}^{1})$, and the difference in filtration 
values
between ${\bf x} \, {\bf y}$ and ${\bf x}_{0} \, {\bf y}$ is, 
$\mathrm{mod\,}\Lambda$, $-n_{\overline{z}}(\phi_{\bf x}^{0})$. Thus, 
given 
filtration indices on the two diagrams, we can construct a filtration 
index on the composite which agrees with the vector sum: 
$\overline{\mathcal{F}}_{0} + \overline{\mathcal{F}}_{1}$. As we are only 
concerned with the ``hatted'' theory we need only identify the coset. 
\\
\\
We now return to the triple $\al_{0}\al_{1}$, $\be_{0}\be_{1}$, and 
$\al_{0}\al_{1}'$. We call the induced four manifold $X$. We choose on $X$ 
the $Spin^{c}$ structure $\mathfrak{u}$ that is $\mathfrak{s} \times I$ 
and restricts to the torsion $Spin^{c}$ structure on the
$\al_{0}\al_{1}, \al_{0}\al'_{1}$- boundary. We then have $\Lambda_{X} 
\equiv \Lambda_{0} \oplus \Lambda_{1} + \Lambda_{II}$. We use a homotopy 
class of triangles to join ${\bf x}_{0} \otimes {\bf y}_{0}$, 
$\Theta^{+}_{std}$, and ${\bf x}_{0}\,{\bf y}_{0}$; a choice made more 
specific below as the unique local holomorphic class. As we assign 
$\Theta^{+}_{std}$ filtration index $\overline{0}$ and this local class 
has $n_{w}(\psi) = n_{\overline{z}}(\psi) = 0$, we have the following 
relation 
for the filtration indices on generators and for some $\lambda \in 
\Lambda_{X}$:
 
$$
\overline{\mathcal{F}}({\bf x\, y}) = \overline{\mathcal{F}}_{0}({\bf 
x})\oplus \overline{0} + \overline{0} \oplus 
\overline{\mathcal{F}}_{1}({\bf y}) + \lambda_{X} 
$$
Topologically, the cobordism, once we fill in the second boundary, is 
$Y_{0}\# Y_{1} \times I$.
If we take the quotient $\mathrm{mod\,}\Lambda_{X}$, we  recover the 
filtration index on $Y$ as $\Z^{2n}/\Lambda_{X} \cong 
\Z^{n}/\Lambda_{III}$ and the filtrations will add correctly. Since 
$z_{i}'$ is in the same domain as $w$ in the diagram for $Y$, there is a 
chain 
isomorphism preserving filtrations which drops their entries in the 
filtration index. Thus, we recover $\widehat{CF}(Y, S_{0} 
\cdot S_{1})$ as a relatively indexed complex (and not, as initially could 
happen,
a quotient of its index group).

The Heegaard triple will be weakly admissible for the doubly periodic 
domains, so we may choose an area form on $\Sigma$ assigning the periodic 
domains signed are equal to zero. As it stands, this may assign large 
portions of the diagram small areas because the periodic domains abutting 
the old meridians from $S_{1}$ in the $\al_{0}\al_{1}, 
\al_{0}\al'_{1}$-boundary are quite substantial. We may address this 
difficulty by handlesliding the portion of the new
$\al'$'s down the diagram for $S_{0}$ until they are close to the 
meridians for $S_{0}$. By doing this, we will have introduced new 
intersections between individual $\al$ and $\be$ curves; however, none of 
these may occur in a generator. Were we to use one of them, there would be 
too few $\al$'s remaining in the upper region of the diagram to pair with 
the $\be$'s found there, and no means to ameliorate this deficiency with 
$\be$'s from the bottom region. Furthermore, nothing in the previous 
analysis will be changed by this alteration. 

In this new diagram, there are obvious holomorphic triangles abutting each 
intersection point ${\bf x} \otimes {\bf y}$ and $\Theta^{+}_{std}$. These 
consist of $g_{0}+g_{1}$ disjoint topological triangles embedded in $\Si$ 
whose domains are contained in the support of the periodic regions from 
the $\al_{0}\al_{1}, \al_{0}\al'_{1}$-boundary. The triangles near the 
meridians for $S_{1}$ are shown in Figure \ref{fig:hhhh}. None of these 
triangles intersect a marked point We may make those periodic domains 
arbitrarily small in unsigned area, forcing our local triangles to have 
area smaller than $\epsilon$. Without the adjustment in the previous 
paragraph,
we would not be able to ensure that only the triangles identified above 
give rise to 
$\epsilon$-``small'' homotopy classes.  Using the induced area 
filtrations, the chain map decomposes into a ``small'' portion, which is 
an isomorphism, and a ``large'' portion:
$$
F( ({\bf x} \otimes {\bf y}) \otimes \Theta^{+}_{std}) = \pm {\bf x\,y} + 
\mathrm{lower\ order}
$$
We see then that the chain map found by counting triangles not crossing 
any marked points induces an injection of $\widehat{CF}(Y_{0} \# Y_{1}, 
S_{0} + S_{1}; \mathfrak{s}, [\overline{j}_{0}] \oplus 
[\overline{j}_{1}])$ into $ \widehat{CF}(Y, S_{0} \cdot 
S_{1};\mathfrak{s}, 
[\overline{j}_{0} + \overline{j}_{1}])$ and that the map is a chain 
isomorphism on $\oplus \widehat{CF}(Y_{0} \# Y_{1}, S_{0} + S_{1}; 
\mathfrak{s}, [\overline{j}'_{0}] \oplus [\overline{j}'_{1}])$ where 
$[\overline{j}_{0} + \overline{j}_{1}] = [\overline{j}'_{0} + 
\overline{j}_{1}]$ $\mathrm{mod\ }\Lambda_{III}$. Together with our 
analysis of boundary I, this proves the result.

Finally, as the small triangles used in the argument each have $n_{w} = 0$ 
and $\mu = 0$, and the cobordism induces the torsion $Spin^{c}$, 
the absolute grading for the image will be the sum of 
the absolute gradings for the original intersection points, when 
$\mathfrak{s}_{i}$ are torsion. Since there are handleslides in the 
$\al_{0}\al'_{1}, \be_{0}\be_{1}$ diagram taking the curves replacing 
meridians in $\al'_{1}$ back to the meridians, and the ``small'' triangles 
in each handleslide map link the corresponding generators, the absolute 
grading for the generators in
this diagram are the same as for $Y_{0}$ shifted by that of $\Theta^{+}$. 
\\
\\
{\bf Note:} $\Theta^{+}$ has grading $\frac{g_{i}}{2}$. But the cobordism 
has $H_{2}$ free, with dimension $g_{1} + g_{2}$ and signature equal to 
zero, so the grading shift formula provides the result.   

\subsubsection{$S_{1}(i,S_{2})$}

\begin{figure}
\begin{center}
\begin{picture}(5778,4140)
\put(942,2719){\makebox(0,0)[lb]{\smash{{\SetFigFont{10}{12.0}{rm}$\beta$'s}}}}
\put(1024,1304){\makebox(0,0)[lb]{\smash{{\SetFigFont{10}{12.0}{rm}$z_{i}$}}}}
\put(2760,891){\makebox(0,0)[lb]{\smash{{\SetFigFont{10}{12.0}{rm}$z_{1}'$}}}}
\put(3555,912){\makebox(0,0)[lb]{\smash{{\SetFigFont{10}{12.0}{rm}$z_{2}'$}}}}
\includegraphics{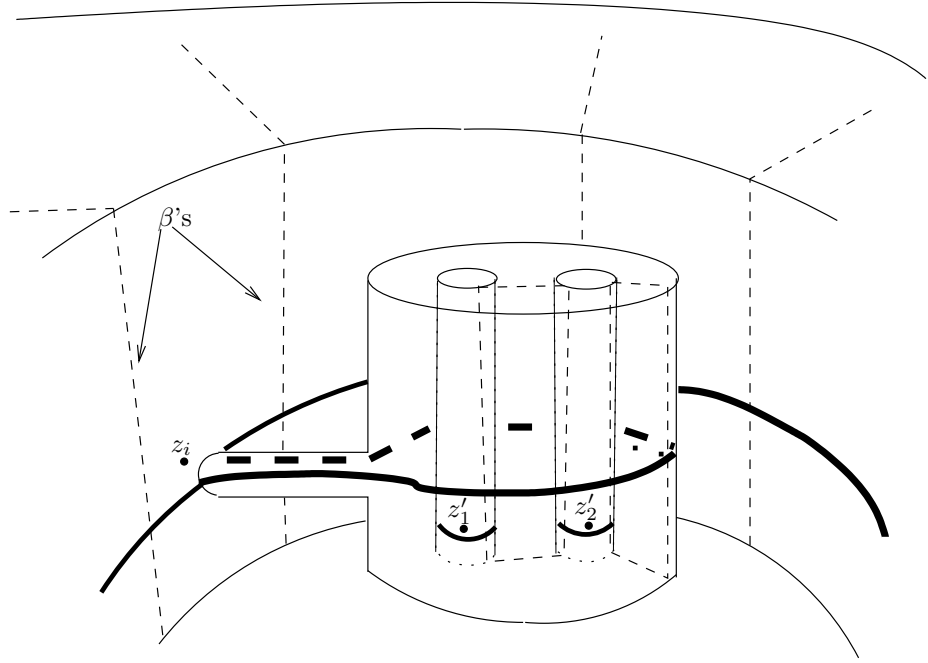}
\end{picture}%
\end{center}
\caption{The Heegaard Diagram for a String Satellite. This is a cutaway 
view of
a diagram for $S_{0}(i, S_{1})$ with the cut through the $i^{th}$ strand 
of $S_{0}$.}
\label{fig:sate}
\end{figure}

The third operation is a form of string satellite to a string link. This 
can be formulated using the Heegaard diagram shown in Figure 
\ref{fig:sate}

Once again, for string links in $S^{3}$, there is only one way to pair 
meridians with $\be$ curves to achieve an intersection point. Indeed, if 
we draw $\al$'s as vertices of a graph and $\be$'s as edges, the meridians 
and the $\be$'s that intersect them form a tree with one edge not 
possessing a vertex on one end. There is only one way to pair edges to end 
points in such a graph. The remaining $\be$'s in the diagram can only 
intersect $\al$'s according to the intersections in the original diagrams. 
However, the construction still applies to string links in more general 
manifolds, presented as surgery on framed links in $D^{2} \times I$. A 
count of $\al$'s and $\be$'s, shows that generators for this new 
diagram occur as products of generators from the old diagrams, even when 
we have wound to achieve some admissibility and possibly increased the 
number of intersections at each meridian. Alternately, we may use the 
standard form to obtain
diagrams to which the argument from $S^{3}$ still applies. Once again, the 
construction is 
a connect sum of two three manifolds, and once again the $Spin^{c}$ 
structures, etc. transfer as expected. 
 
Thus the chain complex is the product of the chain complexes for the 
constituent string links. The filtration indices, however, differ from 
before. To ease the argument, we note that we may think of such a string 
link as the multiplication of one strand in $(Y_{0}, S_{0})$ $n_{1}$ 
times, followed by a composition with $(Y_{1}, S_{1})$ amalgamated with a 
trivial string link on $n_{0}$ strands.  We already know the result of 
composition, hence we need only understand the string satellite where the 
inner constituent is an $n_{1}$ stranded trivial string link in $S^{3}$. 
This has only one intersection point, hence the chain complex, as a group, 
is the same as that for $(Y_{0}, S_{0})$ for each $Spin^{c}$ structure. 

Any class $\phi$ joining two generators, with $n_{w}(\phi) = 0$ can be 
extended to the new string link. It includes the new $\al$ and thus goes 
up the inner string link to the top and back down to the new meridians. 
For each time $\phi$ crosses $z_{i}$, each of the new meridians, $m'_{1}, 
\ldots, m'_{n_{2}}$ will be in the boundary of the new disc. In the 
trivial string link picture there is one generator: $u_{1} \times \cdot 
\times u_{n_{1}}$ with one intersection on each meridian.  In particular, 
the new $\Lambda$ in $\Z^{n_{0} + n_{1} - 1}$ is spanned by vectors $$v' = 
(\lambda_{1}, \ldots, \lambda_{i}, \lambda_{i}, \ldots, \lambda_{i}, 
\ldots, \lambda_{n_{1}})$$ where $\lambda_{i}$ is repeated $n_{1}$ times 
and $(\lambda_{1}, \ldots, \lambda_{n_{0}})$ is a vector in 
$\Lambda_{S_{0}}$. 
We choose the extension of $\phi_{\bf x}$ to $\phi_{{\bf x}\times 
u_{1}\times \cdots \times u_{n_{1}}}$ to give our complete set of paths. 
The filtration index is now measured by  

$$
(\mathcal{F}_{1}, \ldots, \mathcal{F}_{i-1}, \mathcal{F}_{i}, \ldots, 
\mathcal{F}_{i}, \mathcal{F}_{i+1}, \ldots, \mathcal{F}_{n_{1}})   
$$
with $n_{2}$ copies of $\mathcal{F}_{i}$. 

As a consequence, discs with $n_{z_{i}}(\phi) = 0$ in the original diagram 
extend as themselves to the new diagram. In addition, any disc with 
$n_{z'_{i}}(\phi') = 0$, $i = 1, \ldots, n_{2}$ corresponds to a disc in 
the original diagram with $n_{z_{i}} = 0$. The differentials 
$\widehat{\partial}$ and $\widehat{\partial'}$ must be the same, and count 
only classes of discs which do not need to be extended. 

Putting all this together, if we denote the string satellite found by 
substituting $S_{1}$ in the $i^{th}$ strand of $S_{0}$ by $S_{0}(i, 
S_{1})$ then

$$
\begin{array}{c}
\widehat{HF}(Y, S_{0}(i, S_{1}); \mathfrak{s}, [(l_{1}, \ldots, 
l_{n_{1}+n_{2} - 1})] \cong \\ \ \\
\bigoplus_{[\overline{j}'] + [\overline{k}'] = 
[\overline{l}] \mathrm{mod\ }\Lambda'} H_{\ast}(\widehat{CF}(Y_{0}, S_{0}; 
\mathfrak{s}_{0}, [\overline{j}]) \otimes \widehat{CF}(Y_{1}, S_{1}; 
\mathfrak{s}_{1}, [\overline{k}]))
\end{array}
$$
where $\overline{j}' = (j_{1}, \ldots, j_{i-1}, j_{i}, \ldots, j_{i}, 
j_{i+1}, \ldots, j_{n_{1}})$ and $\overline{k}' = (0, \ldots, 0, k_{1}, 
\ldots, k_{n_{2}}, 0, \ldots, 0)$ and $\Lambda' = \Lambda + \overline{0} 
\oplus \Lambda_{1} \oplus \overline{0}$.
\\
\\

We need also to calculate the absolute grading, when appropriate. When we 
have inserted the trivial string link into $S_{0}$ and a torsion 
$Spin^{c}$ structure on $Y_{0}$ we may handleslide the new $\al$ across 
the new meridians to arrive at a picture for a standard connect sum. At 
each handleslide, there is a small $\mu = 0$ homotopy class of triangles 
with $n_{w} = 0$ and admitting holomorphic representative joining each 
intersection point to the corresponding point in the new diagram (the 
product decomposition of generators is unchanged). In the connect sum 
picture, the gradings add -- the grading of the product generator is the 
same as the grading of the generator from $S_{0}$. In the cobordism 
induced by the handleslides, the grading does not change: $gr({\bf x} 
\times u_{1} \times \cdot \times u_{n_{1}}) = $ $gr_{Y_{0}}({\bf x})$
\\
\\
\section{Alexander Invariants for String Links in $D^{2} \times I$}
\label{sec:Euler}
Let $S$ be a string link in $D^{2} \times I$. In this 
section we relate the Euler characteristic of $\widehat{CF}(
S; \overline{v}) \otimes \Q$ to classical Alexander invariants built 
from coverings of $D^{2} \times I - S$. 

\subsection{Alexander Invariants for String Links}

Let $S \subset D^{2} \times I$ be a string link; denote a complement
of its tubular neighborhood by $X = D^{2} \times I - \mathrm{int} N(S)$ 
and 
$X \cap D^{2} \times \{i\}$ by $E_{i}$ for $i=0,\ 1$. By the 
Meyer-Vietoris
sequence we have that $H_{1}(X; \Z) \cong \Z^{k}$ and is generated by the 
meridians of the strands on $D^{2} \times \{0\}$. As string links
are specialized $\thn$-graphs, to construct an Alexander
invariant for the string link we replicate the approach of R. Litherland, 
\cite{Thet}, for the Alexander invariants of $\thn$-graphs (cf. 
\cite{Livi}). Let $\mathfrak{H}$
be the ring $\Z[t_{1}^{\pm 1}, \ldots, t_{k}^{\pm 1}]$; we will consider 
the
torsion properties of the $\mathfrak{H}$-module $H_{1}(\widetilde{X}, 
\widetilde{E}_{1}; \Z)$
where $\widetilde{X}$ is the universal Abelian cover of $X$ (determined by 
the Hurewicz map
$\pi_{1}(X, x_{0}) \ra H_{1}(X, \Z)$ and $\widetilde{E}_{1}$ is the 
pre-image of $E_{1}$ under the covering map.
\\
\\
{\bf Note:}
Taking the lifts to the universal Abelian cover, we consider the long 
exact
sequence for the pair $(\widetilde{X}, \widetilde{E}_{1})$:

$$
\lra H_{1}(\widetilde{E}_{1}) \lra H_{1}(\widetilde{X}) \lra 
H_{1}(\widetilde{X}, \widetilde{E}_{1}) 
\lra H_{0}(\widetilde{E}_{1}) \lra H_{0}(\widetilde{X})
$$
Now, $H_{0}(\widetilde{E}_{1}) \cong \Z$, and it maps isomorphically onto 
the next term.  However, $H_{1}(\widetilde{E}_{1})$ is not generally $0$; 
thus, 
the invariant derived from $H_{1}(\widetilde{X}, \widetilde{E}_{1})$ is 
not an invariant of the complement of $S$ alone. (However, for marked 
knots, i.e. one stranded
string links, $\widetilde{N}$ is an infinite strip, and 
$H_{1}(\widetilde{X}; \Z) \cong H_{1}(\widetilde{X}, \widetilde{N}; \Z)$
as $\mathfrak{H}$-modules.)

We may
construct a relative cell decomposition for $(X, E_{1})$. We think of 
$E_{1}$, a punctured disc, as the portion of the boundary $\partial X - X 
\cap \mathrm{int}
D^{2} \times \{0\}$. Start by constructing a relative cell complex
for $(E_{0}, \partial)$ which consists of $k$ one-cells joining
the internal punctures in a chain to the outer boundary, along the bottom 
of the
projection of $S$. Then add a two-cell to construct the disc. This may
be extended to the entirety of $X$ by attaching one-cells at each of the
crossings in the projection along the axis of the projection and two-cells
for each face in the projection with the exception of the leftmost one, 
called
$U$. The two-cells, $R_{1}, \ldots, R_{k}$, arising from faces that 
intersect the $0^{th}$ level in $I$ must
intersect $D^{2} \times \{0\}$ in one of the one-cells in its 
decomposition.
Otherwise, the two-cells glue to the one-cells at the crossings with the 
remainder of their boundaries
glued into $\partial X$ according to the projection. Finally, the 
complement of this complex is the interior of a three-cell, 
which we glue in to complete $X$. 

We may collapse the two-cell in $E_{0}$ into the union of $E_{1}$ and the 
other
two-cells by contracting the three-cell. Likewise we may collapse the
one-cells in $E_{0}$ into the union of $E_{1}$ and the other
one-cells by contracting $R_{1}, \ldots, R_{k}$ respectively. This 
leaves a relative cell complex with an equal number of $1$- and $2$- 
cells. 
We call this cell complex $Y$. The homotopy and homology
properties of the pair $(X, E_{1})$ are encompassed in this complex.

However, the chain complex for $\widetilde{Y}$ as a relative complex 
becomes

$$
0 \lra C_{2}(\widetilde{Y}) \stackrel{\widetilde{\partial}}{\lra} 
C_{1}(\widetilde{Y}) \lra 0
$$
Thus, $H_{1}(\widetilde{X}, \widetilde{E}_{1}) \cong {\mathrm coker\ 
}\widetilde{\partial}$, 
and the matrix, $P$, for $\widetilde{\partial}$ as
a presentation of the $\mathfrak{H}$-module $H_{1}(\widetilde{X}, 
\widetilde{E}_{1})$ is 
square. By taking the homomorphism $\epsilon: \mathfrak{H} \ra \Z$, 
defined by substituting $1$ for each variable, we see that $\epsilon(P)$ 
is the boundary 
map for the relative chain complex,$Y$. Since $H_{1}(X, E_{1}; \Z) \cong 
0$, as 
$E_{1}$ contains meridians of all types, $P$ has non-zero determinant, 
which we take to be an Alexander 
polynomial of the string link. This is how R. Litherland approaches the abelian
invariants of $\theta_{k}$-graphs.

\subsection{Fox Calculus}

We will calculate the homology of these covering spaces 
is by applying the Fox calculus to the fundamental group of the
complement of the string link, see \cite{Onkn} or \cite{Fox} or
Fox's original articles. 
\\
\\
Let $X$ be a finite, connected $CW$-complex with a single $0$-cell, $p$, 
with no cells 
of  dimension $3$ or greater. Below, we will choose this to be the 
complement of some embedded 
$1$-complex in the three sphere. Let $G = \pi_{1}(X, p)$, be the 
fundamental group based at 
$p$, and presented as $< s_{1}, \ldots, s_{n}\,|\, R_{1}, \ldots, R_{m} 
>$. Let 
$\phi : G \ra H_{1}(X; \Z)$ be the Hurewicz  homomorphism. We consider the 
cover $\widetilde{X}$
determined by this map, and let $\widetilde{X}_{0}$ be the pre-image of 
$p$ 
under the covering map. 

Using the free differentials, we may form 
the matrix 
$\left( \frac{\partial R_{j}}{\partial s_{i}} \right)$ found by first finding
$\frac{\partial}{\partial s_{i}}$ of $R_{j}$, thought of as an element in
the free group, and then applying the quotient maps from the free group to $G$, and
the Hurewicz map to the first homology group. This matrix is a 
presentation matrix for $H_{1}(\widetilde{X}, \widetilde{X}_{0})$ as an $\mathfrak{H}$-module. 

One way to present the fundamental group of the complement $X = D^{2} 
\times I - S$
is to choose a point in $D^{2} \times I - S$ and loops through faces of 
the projection
as generators. Using the crossings to provide the relations, we obtain a 
presentation,
similar to the Dehn presentation for a knot group, with $k$ more 
generators than relations.
For this presentation the Fox calculus produces a presentation matrix with
$k$ more columns than rows. These correspond to the faces 
$R_{1}, \ldots, R_{k}$ and collapsing these faces to obtain 
the
cell complex $Y$ corresponds to eliminating the columns in 
the presentation matrices for Alexander polynomials of knots and links.
 
Furthermore, were we to consider the Heegaard diagrams associated with the 
string link, we would have $k$ meridians along $D^{2} \times \{0\}$. There 
would only be one choice of intersection between these meridians and the 
attaching curves 
derived from the faces, $R_{1}, \ldots, R_{k}$, that could be extended to 
an intersection of $\T_{\al} \cap \T_{\be}$. The collapsing of these faces 
corresponds to this unique choice and to the use of 
$H_{1}(\widetilde{X}, \widetilde{E}_{1})$ as the appropriate classical
analog for the Floer homology. 

This suggests using other presentations of the fundamental group to 
calculate the
Alexander invariant. In \cite{Livi}, P. Kirk, C. Livingston, and Z. Wang 
calculate the invariant
using the analog of the Wirtinger presentation.  
There is a projection of a string link (with kinks at the top 
of each strand, for example) where the Wirtinger presentation is generated 
by 

$$
m_{1}, \ldots, m_{k}, u_{1}, \ldots, u_{s}, m'_{1}, \ldots, m'_{k}
$$
where $m_{i}$ is the meridian of the $i^{th}$ component in 
$D^{2}\times\{1\}$ and $m'_{i}$ is the meridian in $D^{2} \times \{0\}$. 
The kinks
ensure that $m_{i} \neq m'_{i}$. Applying the Fox calculus to 
the relations arising from the crossings we obtain a 
matrix  $(A\ B\ C)$, with entries in $\mathfrak{H}$, where the blocks 
reflect the division in the generators. It is shown in \cite{Livi} that 
$(A\ B)$ is 
invertible over $F = \Q(h_{1}, \ldots, h_{k})$, the quotient field of 
$\mathfrak{H}$. They name the determinant, $\mathrm{det} (A\ B)$, 
the {\em torsion}, $\tau(L)$, of the string link and relate it to the 
Reidemeister
torsion of the based, acyclic co-chain complex $C^{\ast}(X, E_{0}; F)$ 
with coefficients twisted by the map $\pi_{1} \ra \Z^{k}$.  
$S$.

\subsection{Heegaard Splittings and Fox Calculus}

Here we relate the computation of Alexander polynomials to our Heegaard 
diagrams. The upshot will be to identify the Euler characteristic of
$\widehat{CF}(S) \otimes \Q$ with the previously defined Alexander 
invariant of
$H_{1}(\widetilde{X}, \widetilde{E}_{1})$. In particular,
$$
\sum_{\overline{v} \in \Z^{k}} \chi(\widehat{HF}(S, \overline{v}; \Q) 
t_{1}^{v_{1}}\cdots t_{k}^{v_{k}}
$$
is a generator for the order ideal of this module. Various other authors 
have used much the same argument in different
settings; J. Rasmussen provides a very similar argument for Heegaard 
diagrams for 
three manifolds  in \cite{Rasm}.
\\
\\
Let $S$ be a string link in $D^{2} \times I$. We consider the standard
Heegaard decomposition induced from a projection described in section 1.
Let $H_{\al}$ be the handlebody determined by $\seta$, and $H_{\be}$ be 
the handlebody determined by $\setb$. We assume that 
our meridians lie in $D^{2} \times \{0\}$. Take as our basepoint, $p_{0}$ 
for $\pi_1(X)$ the $0$-cell in $H_{\be}$. For each of the faces,
we choose a path $f_i$, the gradient flow line oriented 
from the basepoint to the critical point corresponding to $\be_i$ which 
links the core positively in $S^3$. The other gradient line oriented {\it 
from the index 1 critical point to} $0${\it -cell} will be called 
$\overline{f}_{i}$. The loops $b_{i} =  \overline{f}_{i}\circ f_{i}$  
generate $\pi_{1}(X, p_{0})$.
\\
\\
The $\al$'s, not including the meridians,  induce the relations for a 
presentation of $\pi_{1}$ corresponding to the Dehn presentation of the 
fundamental group, \cite{Onkn}. We
choose an intersection point ${\bf u}\in \T_\al \cap \T_\be$ which 
corresponds to points $u_{i} \in \al_{\sigma(\,i\,)} \cap \be_{i}$ in 
$\Sigma$ for some permutation $\sigma \in S_{g}$. Note that the choice 
along
the meridians is prescribed for each such intersection point: 
there are $k$ meridians and $k+1$ faces intersecting them, but we cast one 
aside. This arrangement implies that our only choices occur on 
non-meridional $\al$'s.  For a 
non-meridional $\al$, let $[\al_{\sigma(\,i\,)}]$ be the path from the 
basepoint, along
$f_{i}$, through the attaching disc for $\be_{i}$ to $u_i$, 
and around $\al_{\sigma(\,i\,)}$ with the its orientation, and then back 
the same way to the basepoint. Each time $[\al_j]$ crosses $\be_i$ 
positively, we append a $b_i$ to the relation;
each time it intersects negatively we append a $b_{i}^{-1}$. The word so 
obtained is called $a_{i}$. We derive this principle by looking at the 
segments $\al_{j}^{s}$ into which the $\be$'s cut $\al_i$ and flowing them 
forwards along the gradient flow. The interior of each segment flows to 
the basepoint, while the endpoints flow to critical points in the 
attaching discs for the $\be$'s.  Thus, the path 
from  one endpoint of the segment, to the critical point corresponding to 
that $\be_s$, then along
some $f_s^{-1}$ or $\overline{f}_s$, and back along one of  $f_t$ or 
$\overline{f}_{t}^{-1}$, then to the
other end of the segment, and back along the segment, is null homotopic. 
This allows us to break
the $\al$ up into the various $\be$'s it crosses. \
\\
\\

Then $\left( \frac{\partial\, a_j}{\partial b_{i}} \right)$, ignoring the 
columns corresponding
to the faces abutting $D^{2} \times \{0\}$, is a 
presentation matrix 
for the $\mathfrak{H}$-module $H_{1}(\tilde{X}, \tilde{E}_{1})$, and hence 
its determinant will provide the Alexander invariant.  
\\
\\
If we consider the free derivative of  $a_{j}$ with respect to a $b_i$ we 
find terms 
which correspond to each intersection point of $\al_j$ with  $\be_i$. The 
term possesses a minus sign when the two intersect negatively, otherwise 
it possesses a positive sign. The terms correspond to paths from the 
basepoint through $f_{\sigma^{-1}(j)}$ to $\be_{\sigma^{-1}(j)}$, through 
the attaching disc to $u_{\sigma^{-1}(j)}$, along $\al_{j}$ to the 
intersection point with $\be_{i}$ and then back along $f_{i}^{-1}$.  This 
can be rewritten as a word in the $b_{i}$'s. Summing over all intersection 
points with $\be_{i}$ equals $\partial_{b_{i}}(a_{j})$.
\\
\\
Let $\mu$ be is the Hurewicz map from the fundamental group to the first 
homology 
group. According to the Fox calculus, the matrix 
$[\mu(\partial_{b_i}{a_j})]$, 
is a presentation matrix for the homology of the cover as an 
$\mathfrak{H}$-module. Again, we ignore the $b_{j}$'s corresponding to the 
faces abutting $D^{2} \times \{0\}$. We calculate the  Alexander invariant 
by computing the 
determinant of this matrix. Each term
in this determinant has the form ${\mathrm sgn}(\sigma)(-1)^{\#} 
h_{1}^{\rho_1} \cdots h_{k}^{\rho_k}$, where
$\rho_i$ is the sum of the powers of $h_i$ over the terms in the 
determinant multiplying to this monomial;   we do not allowing any 
cancellation of terms. This monomial corresponds
to a specific intersection point in $\T_{\al} \cap \T_{\be}$ found from 
the pairing
of rows and columns in the matrix. Likewise 
$\#$ is the number of negative
intersections $\al_{\sigma(\,i\,)} \cap \be_{i}$ in the $g$-tuple 
corresponding to this term.
\\ 
\\
Let ${\bf x}$ and ${\bf y}$ be two intersection points. We will consider 
the differences
$$
\#_{\bf y} - \#_{\bf x} \hspace{1in} \rho_i({\bf y}) - \rho_i({\bf x})
$$
Since we are considering points in $\T_\al \cap \T_\be$ for a diagram of 
$S^3$, there is
a homotopy class of discs $\phi \in \pi_2({\bf x}, {\bf y})$. 
\\
\\
We place marked points into the diagram corresponding to the strands and 
according to the method in section 1. We may measure  
how many times a $2$-chain in $X$, representing a homotopy class 
of  discs, $\phi$, intersects the link components by evaluating $(n_{w} - 
n_{z_{i}})(\phi)$.
We wish to show that 

$$ \rho_i({\bf y}) - \rho_i({\bf x}) = (n_{w} - n_{z_i})(\phi)$$ 
\noindent
The right hand side counts the number of times that that the boundary of 
$\mathcal{D}(\phi)$ winds around
the $i^{th}$ meridian. We need only show that the same is true of the 
left, or, equivalently, that
$\mu_i$, the $i^{th}$-coordinate of the boundary, equals the left hand  
side. 
\\
\\
In the boundary of the disc we have the $\al$'s oriented from the points
in ${\bf x}$ to those in ${\bf y}$. We can take segments starting at $u_i$ 
and travelling along $\al_{\sigma(i)}$ to $x_i$ and $y_i$ so that their 
difference is the oriented segment of the boundary of $\phi$ in  $\al_i$. 
We join this to the basepoint by using paths in the attaching discs for 
the $\be$'s and the preferred paths $f_{j}$ or $f_{j}^{-1}$ at each 
endpoint. Breaking 
this up as before, we can convert this path into a word of $b_i$'s and 
their inverses. If
we look at one of the $\be_i$ boundary segments in $\phi$, we see that the 
concatentation of the words
for the $\al$ segments corresponding to the intersection points with 
$\be_{i}$ homotopes into the 
$\al$ boundary and the $\be$ boundary of $\phi$.

Thus $\mu_i$ of the concatenation equals $\mu_i$ of the boundary of  
$\mathcal{D}(\phi)$. Furthermore, $\mu_i$ applied to each word of the 
concatenation tells us how many more times the segment in one $\al$ 
corresponding to ${\bf y}$, converted into a word of generators,  wraps 
around the $i^{th}$ meridian than does the segment corresponding to ${\bf 
x}$. Taking $\mu_i$ of the concatenation gives the sum of these 
differences, or $\rho_i({\bf y}) - \rho_i({\bf x})$.
\\
\\ 
We now consider the difference in $\#$ between the two intersection 
points. The intersection point determines a permutation $\sigma_{\bf x}$ 
where $x_{i} \in \al_{\sigma(\,i\,)} \cap \be_{i}$. We orient $\T_{\al}$ 
by the projection $\al_{1} \times \cdots \times \al_{g} \ra \T_{\al}$, and 
likewise for $\T_{\be}$. The orientation of $Sym^{g}(\Sigma)$ is given by 
the orientation of $T_{x_{1}}\Sigma \oplus \cdot \oplus T_{x_{g}}\Sigma$. 
Then $\T_{\al} \cap_{\bf x} \T_{\be}$ has local sign 

$${\mathrm sgn} (\sigma_{\bf x}) (-1)^{\frac{g(g-1)}{2}} \cdot 
(\al_{\sigma(\,1\,)} \cap_{x_{1}} 
\be_{1} ) \times \cdots \times (\al_{\sigma(\,g\,)} \cap_{x_{g}} 
\be_{g} )$$
or
$$
 {\mathrm sgn} (\sigma_{\bf x}) (-1)^{\frac{g(g-1)}{2}} (-1)^{\#({\bf x})}
$$
The difference in sign between ${\bf y}$ and ${\bf x}$ is then  
multiplication by  ${\mathrm sgn}(\sigma_{\bf y}) {\mathrm 
sgn}(\sigma_{\bf x}) \cdot ( -1)^{\#_{\bf y} - \#_{\bf x}}$. This is also 
the difference in sign between 
terms in the determinant, and corresponds to the $\Zmod{2}$ grading in 
section 10.4 of \cite{3Man}.
\\
\\
Thus, if we consider those intersection points with $\rho({\bf x}) = 
\overline{v}$, for 
a given vector $\overline{v}$, we recover the intersection points for a 
given filtration
index since $\rho$ satisfies the index relation. In addition, these each 
correspond to the
term $h_{1}^{v_{1}}\cdots h_{k}^{v_{k}}$ and occur with sign given by the 
$\Zmod{2}$ grading
of the Heegaard-Floer homology, which is also the sign of the 
corresponding term in the 
determinant. For rational coefficients, the sum of these generators with 
sign is the Euler characteristic of the homology
group corresponding to $\overline{v}$ for this filtration index. 

\section{State Summation for Alexander Invariants of  String Links
in $S^{3}$}
\label{sec:gencomb}

We consider a generic projection of a string link $S$ in $D^{2} \times I$ 
onto the plane. 
As in \cite{Alte}, \cite{Kauf}, we can 
use this projection
to draw an ancillary rooted, planar graph. The intersection points of 
$\T_{\al} \cap 
\T_{\be}$  will  correspond to a subset of the maximal spanning forests of 
this 
graph, subject to certain constraints imposed by the meridians. From these 
graphs we will prescribe a recipe for computing  functions, 
$\mathcal{F}_{i}$ and $G$, on the intersection points which satisfy the 
same relations relative to homotopy classes $\phi$ as the exponents and 
signs of the Alexander invariants. The $\mathcal{F}_{i}$ will form a  
filtration index, and $G$ will be the grading of our chain complexes. 

Once we have adapted the spanning tree construction to apply to string 
links, the argument is the precise 
analog of the argument in \cite{Alte}. We do not need the results of 
Heegaard-Floer homology 
for the filtration calculation, but the discussion of grading will presume 
some familiarity with them. 

\subsection{Planar Graph Preliminaries}

We consider planar graphs in the unit square, $I^{2}$. Choose a number, 
$k$, and
place $\frac{k}{2}$ vertices, marked by $\ast$, along the bottom edge when 
$k$ is even ($\frac{k+1}{2}$ when $k$ is odd). Place additional vertices, 
labelled by $\bullet$, along the top edge until there are $k$ vertices 
total. Let $\Gamma$ be a {\em connected}, planar graph in $I^{2}$ which 
includes these vertices, but whose other vertices and all its edges are in 
the interior of the square. We let $F$ be a maximal spanning forest for 
$\Gamma$, with a tree component for each $\ast$ on the boundary, rooted at 
$\ast$, and oriented away from its root.

We may define a dual for $\Gamma$ by taking its planar dual inside 
$I^{2}$, $\Gamma^{\ast}$, and placing the vertices that correspond to 
faces 
of $I^{2} - \Gamma$ touching $\partial I^{2}$ 
on $\partial I^{2}$. Since $\Gamma$ is connected, this choice of arc on 
the boundary is
unambiguous. There is one vertex which corresponds to the left side of the 
square. We 
replace it with an $\ast$ and continue counter-clockwise, changing 
boundary vertices to $\ast$'s until we have alterred $\frac{k}{2} + 1$ (or 
$\frac{k+1}{2}$, $k$ odd). This graph must also be connected. 

We say that $F$ admits a dual forest if the edges in $\Gamma^{*}$ 
corresponding to edges
of $\Gamma - F$ form a maximal spanning forest, $F^{\ast}$, with each 
component rooted at a single $\ast$. In that case, we orient the forest 
away from its roots. Not every $F$ admits a dual forest: a component of 
$F^{\ast}$ may contain two $\ast$'s. We consider the set 
$\overline{\mathcal{F}}$ of forests in $\Gamma$ that are part of a dual 
pair $(F, F^{\ast})$.
We will encode $F^{\ast}$ in the diagram for $\Gamma$ by
inscribing the edges in $\Gamma - F$ with a transverse arrow which concurs 
with the 
orientation of $F^{\ast}$.  

\setlength{\unitlength}{3158sp}%
\begin{figure}
\begin{picture}(7654,7029)
\put(5746,5436){\makebox(0,0)[lb]{\smash{{\SetFigFont{12}{14.4}{rm}$\ddots$}}}}
\put(6730,3761){\makebox(0,0)[lb]{\smash{{\SetFigFont{12}{14.4}{rm}$\ast$}}}}
\put(7609,1548){\makebox(0,0)[lb]{\smash{{\SetFigFont{12}{14.4}{rm}$\ast$}}}}
\put(4446,1511){\makebox(0,0)[lb]{\smash{{\SetFigFont{12}{14.4}{rm}$\ast$}}}}
\put(7000,1600){\makebox(0,0)[lb]{\smash{{\SetFigFont{12}{14.4}{rm}$\Gamma^{\ast}$}}}}
\put(1021,3011){\makebox(0,0)[lb]{\smash{{\SetFigFont{8}{14.4}{rm}$L-1$}}}}
\put(5634,3523){\makebox(0,0)[lb]{\smash{{\SetFigFont{12}{14.4}{rm}$\Gamma$}}}}
\includegraphics{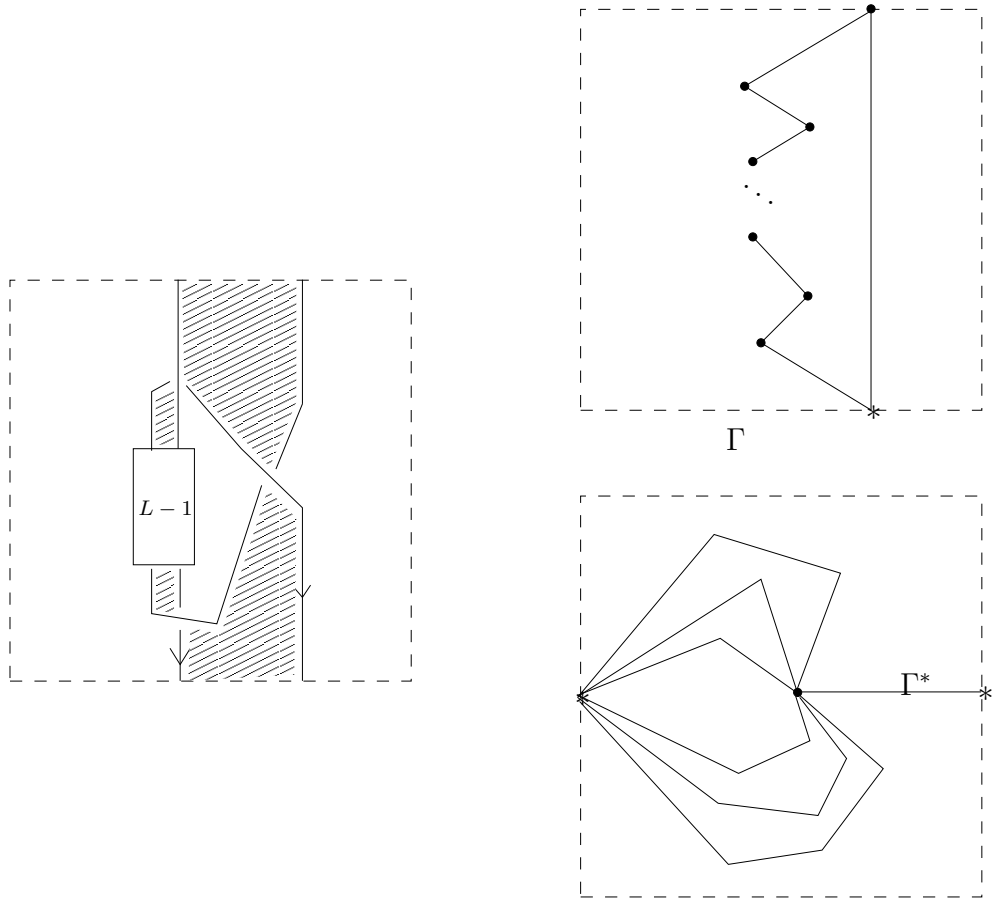}%
\end{picture}%
\caption{An example of the ancillary graph for a simple string link.}
\label{fig:meridians}
\end{figure}
\setlength{\unitlength}{3947sp}%

Now consider a string link, $S$ in $D^{2} \times I$ and a generic 
projection of $S$ into
$I^{2}$. We decompose $S = S_{1} \cup \cdots \cup S_{l}$, where each 
$S_{j}$ consists of 
a maximal string link with connected projection, i.e. one whose projection 
into $I^{2}$ forms 
a connected graph. 

\begin{lemma}
For the Heegaard decomposition of $S^{3}$ defined by the connected 
projection of a string link, $S$, there is a one-to-one correspondence 
between the generators of $\widehat{CF}(S)$ and
the set of dual pairs, $\overline{\mathcal{F}}$, for a planar graph 
$\Gamma \subset I^{2}$
as above.
\end{lemma}

{\bf Proof:}
The regions in $I^{2} - p(S)$ can be colored with $2$ colors as the 
projected graph is 
$4$-valent. We label the leftmost region with the letter ``$U$'' and color 
it white. We
then alternate between black and white across the edges of the projection. 
By using vertices corresponding to the black regions and edges 
corresponding to crossings where two (not 
necessarily distinct) black regions abut, we may form an second planar 
graph. For
those regions touching the border of $I^{2}$, the vertex should be place 
on the boundary. 
 
We replace the vertex of each black region abutting the bottom edge of 
$I^{2}$ with an $\ast$. Thus, we have a graph embedded in $D^{2}$ with  
$k$ vertices on the boundary, $\frac{k+1}{2}$ of which are $\ast$'s when 
$k$ is odd and $\frac{k}{2}$ when $k$ is even. An example is given in 
Figure \ref{fig:meridians}.  In particular, $\Gamma$ is connected: 
starting at a vertex of $\Gamma$ inside a black region of $p(S)$ we can 
take a path to $p(S)$ and follow $p(S)$ to a point in the boundary of any 
other black region. Since any edge of $p(S)$ touches a black region, we 
may then perturb the path into the black regions and find a corresponding 
path in the graph of black regions. If we place $U$'s in all the white 
regions abutting the bottom edge of $I^{2}$, then the graph of white 
regions is connected and, replacing $U$'s with $\ast$'s is the dual graph, 
$\Gamma^{\ast}$, from above. 
\\
\\
{\bf Note:} The dual of $\Gamma$ is always connected, but it is not always 
the graph of white regions. When both $\Gamma$ and the graph of white 
regions are connected, however, the
correspondence is complete. By embedding $\Gamma$ in the closure of the 
black regions, we
can assign to each white region a vertex in $\Gamma^{\ast}$. If two 
distinct white regions
go to the same vertex, there must be a path in the closure of the black 
regions from $\partial I^{2}$ to $\partial I^{2}$ separating the two white 
regions, but not intersecting a vertex of
$p(S)$ (or else the white regions would not map to the same vertex of 
$\Gamma^{\ast}$). Such a path disconnects the graph of white regions. 
\\
\\
Suppose we use the Heegaard decomposition of $S^{3}$ arising from the 
diagram for $S$.
Following \cite{Alte}, we describe an intersection point in $\T_{\al} \cap 
\T_{\be}$ by local data at the vertices of $p(S)$. 
For each non-meridional $\al$ there are four intersection points with 
$\setb$, 
corresponding to the four regions in the projection abutting the crossing 
defined
by $\al$. For each meridian there are one or two intersection points 
depending upon 
whether it intersects the region $U$. However, there can be only one 
choice along all the
meridians which assigns each meridian to a distinct $\be$. We will place a 
$\bullet$ 
in the quadrant corresponding to the intersection point at each vertex.  

Every intersection point corresponds to a pair of dual maximal spanning 
forests in the black and white graphs of the projection. The unique 
choice along the meridians corresponds to the rooting of the forests. We 
then choose the edges in the black graph which join two regions through a 
quadrant marked with a $\bullet$. As each
black region contains a $\bullet$, this produces a subgraph, $F$, 
containing all the vertices of $\Gamma$. We can perform the same operation 
in the white graph to obtain a second sub-graph. 

Furthermore, all the components of these sub-graphs are trees.  A cycle in 
$F$ would bound a disc in $S^{2}$ not containing a region labelled $U$. 
Rounding the crossings of $p(S)$ along the  cycle, we find  a $4$-valent 
planar graph with 
$B_{in}$ crossings and  $B_{in} + 1$ faces not touching the cycle. The 
original intersection point  must form a $1-1$ correspondence between 
these faces and crossings as all the surrounding faces were consumed by 
the cycle. There can be no such identification and thus no cycle in the 
black graph. Similarly, if the intersection point does not produce a 
forest in the
white graph there is a contradiction. Thus we have two maximal spanning 
forests.

Every component must contain precisely one $\ast$. It cannot contain more 
as there
is a $\bullet$ for each edge in the tree and for each root. In order for 
the number of
edges plus roots (the $\al$'s) to equal the number of $\be$'s there must 
be precisely one
root. Thus the two sub-graphs are a dual pair of maximal spanning forests 
for the graphs of black and white regions

Conversely, the arrows on the edges of $\Gamma$ found from a dual pair 
$(F, F^{\ast})$ tell us
how to complete the assignment of $\al$'s to $\be$'s from the unique 
assignment along the meridians: for each non-meridional $\al_{i}$ we 
choose the intersection point in $\al_{i} \cap \be_{\sigma(\,i\,)}$ 
pointed towards by the arrow on the edge corresponding to the crossing 
defined by $\al_{i}$. The existence of $F^{\ast}$ ensures that no arrow 
contradicts 
the assignment along the meridians by pointing into a region labelled with 
$U$.

\subsection{A Variant of the Clock Theorem}

\begin{figure}
\begin{center}
\begin{picture}(5060,1767)(1371,-2448)
\thinlines
\put(2185,-2209){\vector( 0, 1){517}}
\put(2451,-1511){\circle*{70}}
\put(1989,-1561){\circle*{70}}
\put(2151,-1374){\circle*{70}}
\put(2126,-724){\circle*{70}}
\put(3076,-1661){\circle*{70}}
\put(4726,-1739){\circle*{70}}
\put(5763,-1576){\circle*{70}}
\put(5301,-1626){\circle*{70}}
\put(5463,-1439){\circle*{70}}
\put(5438,-789){\circle*{70}}
\put(6388,-1726){\circle*{70}}
\put(3540,-1061){\vector(-1, 0){  0}}
\put(3540,-1061){\vector( 1, 0){893}}
\put(1414,-1674){\circle*{70}}
\put(1421,-1636){\line( 3, 4){705}}
\multiput(1878,-1248)(-6.71429,5.03571){29}{\makebox(1.6667,11.6667){\SetFigFont{5}{6}{rm}.}}
\put(1878,-1248){\vector( 4,-3){0}}
\multiput(2136,-1353)(-5.03571,-6.71429){29}{\makebox(1.6667,11.6667){\SetFigFont{5}{6}{rm}.}}
\put(2139,-738){\line( 0,-1){282}}
\put(2139,-1020){\line( 0,-1){329}}
\put(2139,-1349){\line( 2,-1){282}}
\multiput(5448,-1418)(-5.03571,-6.71429){29}{\makebox(1.6667,11.6667){\SetFigFont{5}{6}{rm}.}}
\put(5440,-762){\line( 1,-1){940}}
\put(5451,-803){\line( 0,-1){282}}
\put(5451,-1085){\line( 0,-1){329}}
\put(5451,-1414){\line( 2,-1){282}}
\put(6102,-1110){\vector(-1,-1){282}}
\put(2139,-724){\line( 1,-1){940}}
\thicklines
\put(3050,-1610){\vector(-1, 1){799}}
\put(4750,-1679){\vector( 3, 4){564}}
\thinlines
\put(4723,-1704){\line( 3, 4){705}}
\put(1374,-2448){\makebox(0,0)[lb]{\smash{{\SetFigFont{12}{14.4}{rm}A 
sub-graph 
wholly contained in the face}}}}
\end{picture}%
\end{center}
\caption{The clock $\ra$ and counter-clock $\leftarrow$ moves on maximal 
forests.} 
\label{fig:clockmove}
\end{figure}
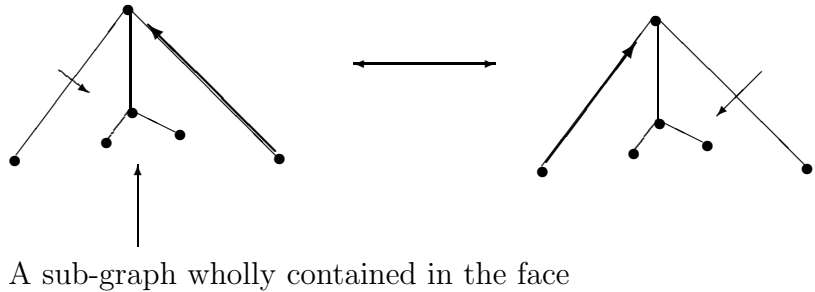

We will now examine the structure of the set of dual maximal spanning 
forests. 
As in \cite{Cloc} we consider two moves performed on the decorations a 
dual pair
inscribes on $\Gamma$: the clock and counter-clock moves. These are
moves interchange the two pictures in Figure \ref{fig:clockmove}. There 
should be 
a face -- not labelled with a $U$ -- of $I^2 - \Gamma$ abutting these two 
edges at 
their common vertex. This allows a portion of $\Gamma$ to be 
wholly contained in the interior of the face.  

A clock move performed on $F$ in $\Gamma$ corresponds to a clock move 
performed 
on $F^{\ast}$ in $\Gamma^{\ast}$. These moves take dual maximal forest 
pairs
into dual maximal forest pairs. If the edge with the transverse arrow 
joins
two distinct components of $F$, then the new oriented sub-graph, $F'$, of 
$\Gamma$
after the clock move is still a forest. The portion of the component 
through the vertex
beyond the vertex is a tree not containing the other vertex of the 
transverse edge. The move
merely prunes this section of one component and glues it to another 
component. If
the edge with transverse arrow joins vertices in the same component of 
$F$, it is conceivable
that a cycle could form. However, this can only happen if the arrow on the 
transverse edge
points out of the disc bounded by this cycle, and thus a root of the dual 
graph must be contained in the cycle. Since those roots lie on $\partial 
I^{2}$, this cannot happen.

For a connected, finite 
planar graph, $\Gamma$ with only one root, the structure
of maximal spanning trees is already understood, \cite{Cloc}. We require 
that the
root be in the boundary of $U$, the unbounded component of $\R^2 
- \Gamma$. Pick one of the trees, $T$ . Each additional edge in $\Gamma$, 
when 
adjoined to the tree, divides
the plane into a bounded and an unbounded component. Draw an arrow 
pointing into the bounded component on each of these edges. This is the 
decoration inscribed
by the dual tree as before.

In \cite{Cloc}, Gilmer and Litherland prove Kauffman's clock theorem:

\begin{theorem}{\bf The Clock Theorem}
The set $\mathcal{T}$ of maximal, spanning trees is
a graded, distributive lattice under the partial order defined by $T \ge 
T'$ if we can move 
from $T$ to $T'$ solely by using clock-moves.
\end{theorem}
\noindent
We  will only need that any $T \in \mathcal{T}$ can be obtained from any 
other $T'$ by making clock and counter-clock moves. It is shown 
in \cite{Cloc} that only a finite
number of clock (or counter-clock) moves can be made successively before 
we reach a tree not admitting 
another. Furthermore, this tree is unique for the type of move. We 
can go from any tree to any 
other by continually making clockwise moves until we reach the unique 
un-clocked tree and 
then make counter-clock moves to get to the other tree. 

\subsubsection{The Clock Theorem for Forests}

The vertices on the boundary of $I^{2}$ divide the boundary 
into arcs. We draw an arrow into the the regions labelled by $U$, across 
the corresponding arcs. Place arrows pointing out along the other edges. A 
dual pair $(F, F^{\ast})$ for the string link $S$ extends these arrows in 
the sense that each face has exactly one arrow pointing
into it.

\setlength{\unitlength}{3158sp}%
\begin{figure}
\begin{center}
\begin{picture}(8148,1848)
\put(5625,959){\makebox(0,0)[lb]{\smash{{\SetFigFont{12}{14.4}{rm}$\ldots$}}}}
\put(1113,934){\makebox(0,0)[lb]{\smash{{\SetFigFont{12}{14.4}{rm}$\ldots$}}}}
\includegraphics{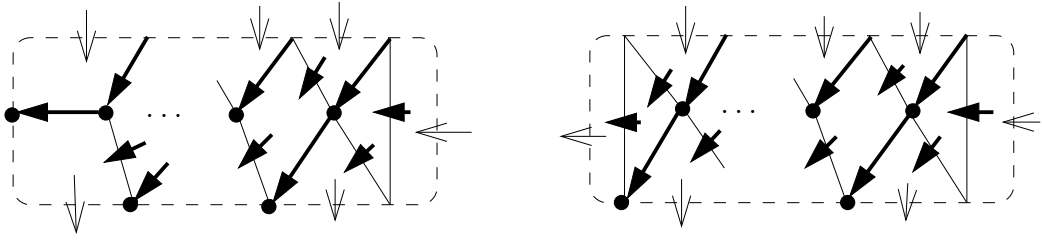}%
\end{picture}%
\end{center}
\caption{Graphs in $I^{2}$ admitting a unique pair of maximal spanning 
forests. In particular, no counter-clock moves can be performed in them. 
Furthermore, embedded appropriately in a planar graph with a maximal tree 
inscribing the same decorations, no counter-clock move in the planar graph 
can alter the decorations in this region. We use the one on the left when 
$k$ is even, the one on the right when $k$ is odd.} \label{fig:extensor}
\end{figure}
\setlength{\unitlength}{3947sp}%

We use the arrows on $\partial I^{2}$ to extend $\Gamma$ to a planar graph 
with a single root, so that each dual pair $(F, F^{\ast})$ corresponds to 
a unique maximal spanning {\em tree} in the resulting graph. To do this we 
consider a new square with the reverse of the decorations on the boundary 
of the old square, with the exception that the one on the right edge and 
the first root from right to left remain the same. We can then extend 
these decorations by 
a graph, $\Gamma'$, and the dual pair as in Figure \ref{fig:extensor}. By 
inspection, these are the unique decorations providing a dual pair of 
forests in this graph and extending the boundary conditions. (No two 
arrows may point into the same white region, if the dual
is to be a forest). 

We can glue this decorated graph to the one from $S$ to obtain a planar 
graph where there is
only one $U$, corresponding to the leftmost edge of $S$, one $\ast$, and 
the pair $(F, F^{\ast})$ becomes the decorations from a maximal spanning 
{\em tree} as in the clock theorem. Furthermore, a  maximal spanning tree 
in the glued graph inscribing the decorations on the $\Gamma'$-portion as 
in Figure \ref{fig:extensor} corresponds to a maximal pair $(F, 
F^{\ast})$ in the graph of black regions for $S$.

Now we perform counter-clock moves until we reach the maximally clocked 
tree. At no time
do these moves disrupt the decorations in the $\Gamma'$-portion. No such 
move can occur on
a face in $\Gamma'$ as there is no vertex with the requisite arrangement 
of tree edge and 
transverse arrow. Furthermore, the $\Gamma'$ region can be disrupted from 
outside only
when a counter-clock move occurs on a face abutting $\Gamma$. Noting that 
the arrows point
out of the vertices on the bottom, and into the vertices on the top, 
inspection
shows that no counter-clock move can occur on such a face at a vertex from 
$\Gamma'$. This is
{\em not} true for clock moves, which can occur on the top left of 
$\Gamma'$. Finally,
since the transverse arrows point into the faces that were formerly 
labelled by $U$, and
these arrows are never altered, no counter-clock move ever involves a face 
formerly labelled
with $U$. However, as any forest pair for $S$ may be extended to a tree 
for the new graph and counter-clocked to the maximal clocked tree, there 
is always a sequence of counter-clock and clock moves, not involving 
$\Gamma'$, which connect any two pairs for $S$.

In short,

\begin{lemma}
The dual pairs $(F, F^{\ast})$ for the graph of black regions, $\Gamma$, 
found for 
a connected projection of a string link may be converted, one into 
another, by
clock and counter-clock moves performed on the decorations coming from the 
rooting
of the string link.
\end{lemma}

\subsection{State Summation}

\begin{figure}
\begin{center}
\begin{picture}(6747,4234)(801,-4800)
\thicklines
\put(1538,-3549){\vector(-1, 1){703}}
\put(2031,-1413){\line(-1,-1){318}}
\put(2120,-1334){\vector( 1, 1){313}}
\put(6021,-1703){\vector( 1, 1){710.500}}
\put(6731,-1719){\line(-1, 1){313}}
\put(6318,-1305){\vector(-1, 1){297}}
\thinlines
\thicklines
\thinlines
\thicklines
\put(4947,-3474){\vector( 1, 1){709.500}}
\thinlines
\put(5351,-3179){\line( 1,-1){307}}
\put(5255,-3084){\vector(-1, 1){309}}
\put(3582,-3557){\vector(-1, 1){705}}
\thicklines
\put(3248,-3170){\vector( 1, 1){338}}
\put(3190,-3222){\line(-1,-1){326.500}}
\put(2429,-1731){\vector(-1, 1){705}}
\thinlines
\put(1233,-3151){\vector( 1, 1){315}}
\put(1146,-3250){\line(-1,-1){303.500}}
\put(6818,-3461){\vector( 1, 1){705}}
\thicklines
\put(7136,-3090){\vector(-1, 1){338}}
\put(7188,-3148){\line( 1,-1){326.500}}
\put(6248,-885){\makebox(0,0)[lb]{\smash{{\SetFigFont{12}{14.4}{rm}$-\frac{1}{2}$}}}}
\put(4659,-3184){\makebox(0,0)[lb]{\smash{{\SetFigFont{12}{14.4}{rm}$-\frac{1}{2}$}}}}
\put(2708,-3244){\makebox(0,0)[lb]{\smash{{\SetFigFont{12}{14.4}{rm}$\frac{1}{2}$}}}}
\put(6974,-2670){\makebox(0,0)[lb]{\smash{{\SetFigFont{12}{14.4}{rm}$-\frac{1}{2}$}}}}
\put(7412,-3169){\makebox(0,0)[lb]{\smash{{\SetFigFont{12}{14.4}{rm}$-\frac{1}{2}$}}}}
\put(2026,-799){\makebox(0,0)[lb]{\smash{{\SetFigFont{12}{14.4}{rm}$\frac{1}{2}$}}}}
\put(1901,-1949){\makebox(0,0)[lb]{\smash{{\SetFigFont{12}{14.4}{rm}$-\frac{1}{2}$}}}}
\put(6376,-1911){\makebox(0,0)[lb]{\smash{{\SetFigFont{12}{14.4}{rm}$\frac{1}{2}$}}}}
\put(1514,-3261){\makebox(0,0)[lb]{\smash{{\SetFigFont{12}{14.4}{rm}$\frac{1}{2}$}}}}
\put(1201,-2786){\makebox(0,0)[lb]{\smash{{\SetFigFont{12}{14.4}{rm}$\frac{1}{2}$}}}}
\put(3214,-2774){\makebox(0,0)[lb]{\smash{{\SetFigFont{12}{14.4}{rm}$\frac{1}{2}$}}}}
\put(5276,-2674){\makebox(0,0)[lb]{\smash{{\SetFigFont{12}{14.4}{rm}$-\frac{1}{2}$}}}}
\end{picture}%
\end{center}
\caption{Filtration Weights Depicted for the Thick Strand and the $i^{th}$ 
Filtration Value. The meridians do not contribute to the weights.}
\label{fig:weights}
\end{figure}
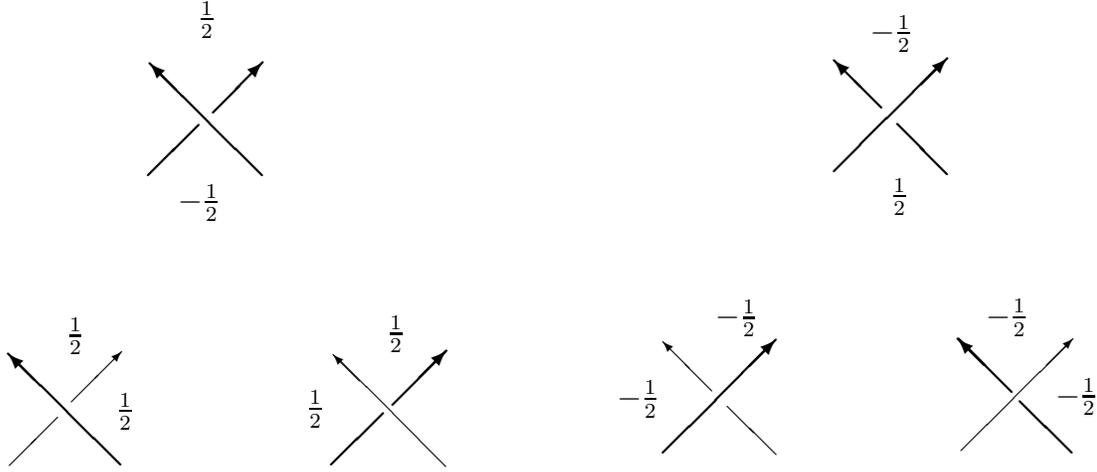

Following \cite{Alte}, we will prescribe weights at the crossings of the 
string link as in Figure \ref{fig:weights}, and extend those weights to 
apply to 
more than one component. For each intersection point in $\T_{\al} \cap 
\T_{\be}$, we consider the associated dual pair $(F, F^{\ast})$ and 
locally place $\bullet$'s in the regions  abutting
each intersection point or merdian according to the direction of the 
trees. We sum the weights from 
the marked region at each crossing and meridian. This prescribes 
$\mathcal{F}_{i}$ for the $i^{th}$ strand
if it is the thick 
strand.  Crossings not involving the thick strand 
contribute nothing to the weight.  When we change from the intersection 
point 
obtained from one tree to that obtained from another tree differing by a 
clock or 
counter-clock move, the change in exponent is given by $(n_w - 
n_{z_i})(\phi)$ for 
any $\phi$ in $\pi_2({\bf x},{\bf y})$ that joins these intersection 
points. We verify 
that this equals the change in the weights. As the Heegaard diagram 
describes $S^3$, an integer homology sphere, such a class $\phi$ must 
exist.  Since the dual pair $(F, F^{\ast})$ coming from an intersection 
point
can be obtained from the dual pair for any other intersection point by 
clock and counter-clock moves, we may compute the difference in exponents 
for any two intersection points by looking at the difference in the 
overall weights.

\begin{figure}
\begin{center}
\begin{picture}(3088,966)(2754,-1828)
\thicklines
\put(3105,-1477){\line(-1,-1){318}}
\put(3194,-1398){\vector( 1, 1){313}}
\put(3041,-1040){\makebox(0,0)[lb]{\smash{{\SetFigFont{12}{14.4}{rm}$1$}}}}
\put(5099,-1779){\vector( 1, 1){710.500}}
\put(5809,-1795){\line(-1, 1){313}}
\put(5396,-1381){\vector(-1, 1){297}}
\put(5341,-1007){\makebox(0,0)[lb]{\smash{{\SetFigFont{12}{14.4}{rm}$-1$}}}}
\put(3503,-1795){\vector(-1, 1){705}}
\end{picture}%
\end{center}
\caption{Weights for the Absolute Grading}
\label{fig:gradweight}
\end{figure}
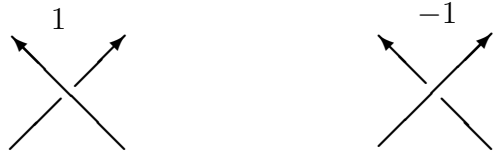

Likewise, there are weights to calculate the grading for each intersection 
point. See Figure
\ref{fig:gradweight}. For each intersection point we add the weights over 
all the 
crossings without considering which strands appear. The meridians do not 
contribute to the grading.
\\
\\

We now prove that the difference in grading and filtration values between 
one maximal forest
and the forest that results after a clock or counter-clock move equals the 
difference in the weights 
defined for each tree above. We denote the local contribution to each 
intersection point by placing a $\bullet$ or $\circ$ at each crossing. We 
will assume that $\bullet$ and $\circ$ are identical for crossings that 
are not depicted. Following the definitions of maximal trees and clock and 
counter-clock moves given above, we can verify that the moves from $\circ$ 
to $\bullet$ fulfill our requirements and exhaust all possible moves. We 
break the argument up into cases:
\\
\\
{\bf Case I:} Figure \ref{fig:casei} shows the cases where a counter-clock 
move joins 
intersection points at two crossings (not meridians). 
These correspond to unique discs, namely squares, ``atop" the Heegaard 
surface.  As the squares do not cross  any of the multi-points, there will 
be no change in exponents corresponding to any of the three strands. This 
equals  the change in the weights. On the other hand. squares always have 
a one dimensional space of holomorphic 
representatives, so the intersection point $\circ$ has grading $1$ greater 
than $\bullet$ in Heegaard-Floer homology. This equals the change in the 
grading weights. 
\\
\\
\setlength{\unitlength}{3000sp}%
\begin{figure}
\begin{center}
\begin{picture}(8163,3446)(695,-3195)
\thinlines
\put(1883,-369){\circle{80}}
\put(1901,-2642){\circle*{80}}
\put(3277,-2648){\circle{80}}
\put(3264,-2336){\circle*{80}}
\put(3402,-1823){\vector( 0,-1){1350}}
\thicklines
\put(1692,-2473){\vector(-1, 0){940}}
\put(1862,-2473){\line( 1, 0){1410}}
\put(3527,-2473){\line( 1, 0){799}}
\thinlines
\put(1777,-3148){\vector( 0, 1){1387}}
\put(1251,-2274){\makebox(0,0)[lb]{\smash{{\SetFigFont{12}{14.4}{rm}$-\frac{1}{2}$}}}}
\put(1839,-2299){\makebox(0,0)[lb]{\smash{{\SetFigFont{12}{14.4}{rm}$-\frac{1}{2}$}}}}
\put(2876,-2774){\makebox(0,0)[lb]{\smash{{\SetFigFont{12}{14.4}{rm}$+\frac{1}{2}$}}}}
\put(3551,-2774){\makebox(0,0)[lb]{\smash{{\SetFigFont{12}{14.4}{rm}$+\frac{1}{2}$}}}}
\put(6384,-2268){\circle{80}}
\put(6389,-2579){\circle*{80}}
\put(7765,-2585){\circle{80}}
\put(7752,-2273){\circle*{80}}
\put(7890,-1760){\vector( 0,-1){1350}}
\thicklines
\put(6180,-2410){\vector(-1, 0){940}}
\put(6350,-2410){\line( 1, 0){1410}}
\put(8015,-2410){\line( 1, 0){799}}
\thinlines
\put(6265,-1698){\vector( 0,-1){1387}}
\put(7177,-2748){\makebox(0,0)[lb]{\smash{{\SetFigFont{12}{14.4}{rm}$+\frac{1}{2}$}}}}
\put(8027,-2786){\makebox(0,0)[lb]{\smash{{\SetFigFont{12}{14.4}{rm}$+\frac{1}{2}$}}}}
\put(6377,-2761){\makebox(0,0)[lb]{\smash{{\SetFigFont{12}{14.4}{rm}$+\frac{1}{2}$}}}}
\put(5752,-2761){\makebox(0,0)[lb]{\smash{{\SetFigFont{12}{14.4}{rm}$+\frac{1}{2}$}}}}
\put(1896,-2331){\circle{80}}
\put(1888,-680){\circle*{80}}
\put(3264,-686){\circle{80}}
\put(3251,-374){\circle*{80}}
\put(6308,-331){\circle{80}}
\put(6313,-642){\circle*{80}}
\put(7689,-648){\circle{80}}
\put(7676,-336){\circle*{80}}
\put(3389,-1211){\vector( 0, 1){1350}}
\thicklines
\put(1679,-511){\vector(-1, 0){940}}
\put(1849,-511){\line( 1, 0){1410}}
\put(3514,-511){\line( 1, 0){799}}
\thinlines
\put(1764,-1186){\vector( 0, 1){1387}}
\put(7814,-1173){\vector( 0, 1){1350}}
\thicklines
\put(6104,-473){\vector(-1, 0){940}}
\put(6274,-473){\line( 1, 0){1410}}
\put(7939,-473){\line( 1, 0){799}}
\thinlines
\put(6189,239){\vector( 0,-1){1387}}
\put(1851,-236){\makebox(0,0)[lb]{\smash{{\SetFigFont{12}{14.4}{rm}$-\frac{1}{2}$}}}}
\put(1314,-236){\makebox(0,0)[lb]{\smash{{\SetFigFont{12}{14.4}{rm}$-\frac{1}{2}$}}}}
\put(2939,-249){\makebox(0,0)[lb]{\smash{{\SetFigFont{12}{14.4}{rm}$-\frac{1}{2}$}}}}
\put(3551,-236){\makebox(0,0)[lb]{\smash{{\SetFigFont{12}{14.4}{rm}$-\frac{1}{2}$}}}}
\put(7926,-274){\makebox(0,0)[lb]{\smash{{\SetFigFont{12}{14.4}{rm}$-\frac{1}{2}$}}}}
\put(7351,-249){\makebox(0,0)[lb]{\smash{{\SetFigFont{12}{14.4}{rm}$-\frac{1}{2}$}}}}
\put(6351,-836){\makebox(0,0)[lb]{\smash{{\SetFigFont{12}{14.4}{rm}$+\frac{1}{2}$}}}}
\put(5739,-824){\makebox(0,0)[lb]{\smash{{\SetFigFont{12}{14.4}{rm}$+\frac{1}{2}$}}}}
\end{picture}%
\end{center}
\caption{Case I: Weights, depicted for the horizontal strand, do not 
change 
under the alteration from $\circ$ to $\bullet$. The thin strands do not 
need to come from the same component.  The thin component on the left
receives weight $0$ from these configurations. The thin strand on the 
right receives the same weight from $\circ$ and $\bullet$. The grading 
change occurs 
along the middle strand, as inspection of the crossings shows. }
\label{fig:casei}
\end{figure}
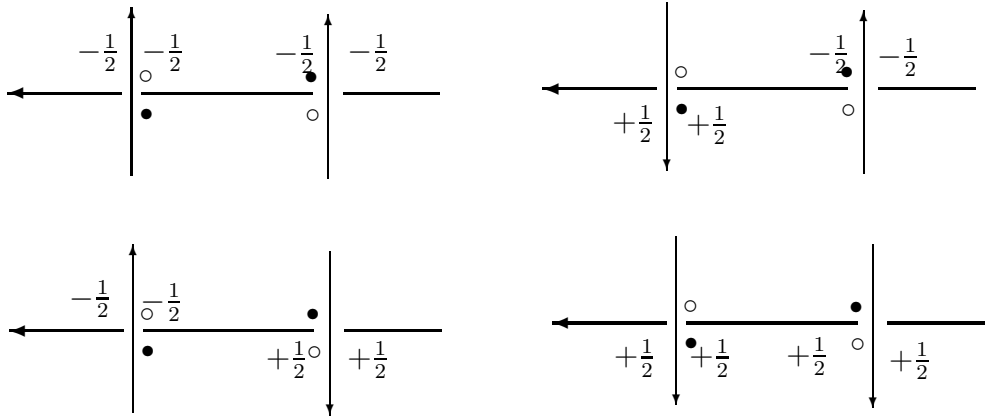
\setlength{\unitlength}{3947sp}%
\noindent
{\bf Case II:} Figure \ref{fig:caseii} shows the same alteration but with 
a different 
configuration
of under and over-crossings. The homotopy class we choose is now a 
``square" with punctures and handles added to it. In particular, the disc 
travels off the end of the figure in the direction of the knot picking up 
punctures
at crossings and joins punctures into handles if it happens to go through 
the same crossing twice. It terminates
on the meridian corresponding to the horizontal strand. Thus, the 
filtrations 
remain unchanged except in the $i^{th}$ component. But the ``disc" passes 
over $z_i$ once, so $\mathcal{F}_i(\bullet) - \mathcal{F}_i(\circ)$ $= 
(n_w - 
n_{z_i})(\phi)$ $= -1$. In \cite{Alte}, \Oz and \Sz show that $\mu(\phi) = 
1$ for such a class, so the grading change equals the change in grading 
weights.
\\
\\
\setlength{\unitlength}{3000sp}%
\begin{figure}
\begin{center}
\begin{picture}(8163,3394)(695,-3195)
\thinlines
\put(3402,-1823){\vector( 0,-1){1350}}
\put(1776,-2386){\vector( 0, 1){564}}
\put(1764,-1173){\line( 0, 1){564}}
\put(1764,-423){\vector( 0, 1){564}}
\put(1883,-369){\circle{80}}
\put(1888,-680){\circle*{80}}
\put(3264,-686){\circle{80}}
\put(3251,-374){\circle*{80}}
\put(6308,-331){\circle{80}}
\put(6313,-642){\circle*{80}}
\put(7689,-648){\circle{80}}
\put(7676,-336){\circle*{80}}
\put(1896,-2331){\circle{80}}
\put(1901,-2642){\circle*{80}}
\put(3277,-2648){\circle{80}}
\put(3264,-2336){\circle*{80}}
\put(6384,-2268){\circle{80}}
\put(6389,-2579){\circle*{80}}
\put(7765,-2585){\circle{80}}
\put(7752,-2273){\circle*{80}}
\put(3389,-1211){\vector( 0, 1){1350}}
\thicklines
\put(1867,-511){\vector(-1, 0){1128}}
\put(1849,-511){\line( 1, 0){1410}}
\put(3514,-511){\line( 1, 0){799}}
\thinlines
\put(7814,-1173){\vector( 0, 1){1350}}
\thicklines
\put(6292,-473){\vector(-1, 0){1128}}
\put(6274,-473){\line( 1, 0){1410}}
\put(7939,-473){\line( 1, 0){799}}
\thinlines
\put(1776,-3136){\line( 0, 1){564}}
\thicklines
\put(1674,-2473){\line( 1, 0){1598}}
\put(3527,-2473){\line( 1, 0){799}}
\thinlines
\put(7890,-1760){\vector( 0,-1){1350}}
\thicklines
\put(6462,-2410){\vector(-1, 0){1222}}
\put(6350,-2410){\line( 1, 0){1410}}
\put(8015,-2410){\line( 1, 0){799}}
\put(1833,-2474){\vector(-1, 0){1081}}
\thinlines
\put(6264,-2522){\vector( 0,-1){564}}
\put(6264,-2336){\line( 0, 1){564}}
\put(6176,-572){\vector( 0,-1){564}}
\put(6176,-386){\line( 0, 1){564}}
\put(1390,-2201){\makebox(0,0)[lb]{\smash{{\SetFigFont{12}{14.4}{rm}$+\frac{1}{2}$}}}}
\put(1954,-2201){\makebox(0,0)[lb]{\smash{{\SetFigFont{12}{14.4}{rm}$+\frac{1}{2}$}}}}
\put(3524,-2865){\makebox(0,0)[lb]{\smash{{\SetFigFont{12}{14.4}{rm}$+\frac{1}{2}$}}}}
\put(2719,-2845){\makebox(0,0)[lb]{\smash{{\SetFigFont{12}{14.4}{rm}$+\frac{1}{2}$}}}}
\put(6383,-2744){\makebox(0,0)[lb]{\smash{{\SetFigFont{12}{14.4}{rm}$-\frac{1}{2}$}}}}
\put(5658,-2765){\makebox(0,0)[lb]{\smash{{\SetFigFont{12}{14.4}{rm}$-\frac{1}{2}$}}}}
\put(7389,-2765){\makebox(0,0)[lb]{\smash{{\SetFigFont{12}{14.4}{rm}$+\frac{1}{2}$}}}}
\put(8134,-2805){\makebox(0,0)[lb]{\smash{{\SetFigFont{12}{14.4}{rm}$+\frac{1}{2}$}}}}
\put(7268,-248){\makebox(0,0)[lb]{\smash{{\SetFigFont{12}{14.4}{rm}$-\frac{1}{2}$}}}}
\put(7973,-248){\makebox(0,0)[lb]{\smash{{\SetFigFont{12}{14.4}{rm}$-\frac{1}{2}$}}}}
\put(6322,-872){\makebox(0,0)[lb]{\smash{{\SetFigFont{12}{14.4}{rm}$-\frac{1}{2}$}}}}
\put(5577,-832){\makebox(0,0)[lb]{\smash{{\SetFigFont{12}{14.4}{rm}$-\frac{1}{2}$}}}}
\put(1853,-208){\makebox(0,0)[lb]{\smash{{\SetFigFont{12}{14.4}{rm}$+\frac{1}{2}$}}}}
\put(2699,-208){\makebox(0,0)[lb]{\smash{{\SetFigFont{12}{14.4}{rm}$-\frac{1}{2}$}}}}
\put(3564,-168){\makebox(0,0)[lb]{\smash{{\SetFigFont{12}{14.4}{rm}$-\frac{1}{2}$}}}}
\put(1249,-208){\makebox(0,0)[lb]{\smash{{\SetFigFont{12}{14.4}{rm}$+\frac{1}{2}$}}}}
\end{picture}%
\end{center}
\caption{Case II: Weights, depicted for the horizontal strand,  reduce by 
$1$ under the alteration from $\circ$ to $\bullet$. The same comments as 
in the caption for case I apply to the thin strands and the grading. }
\label{fig:caseii}
\end{figure}
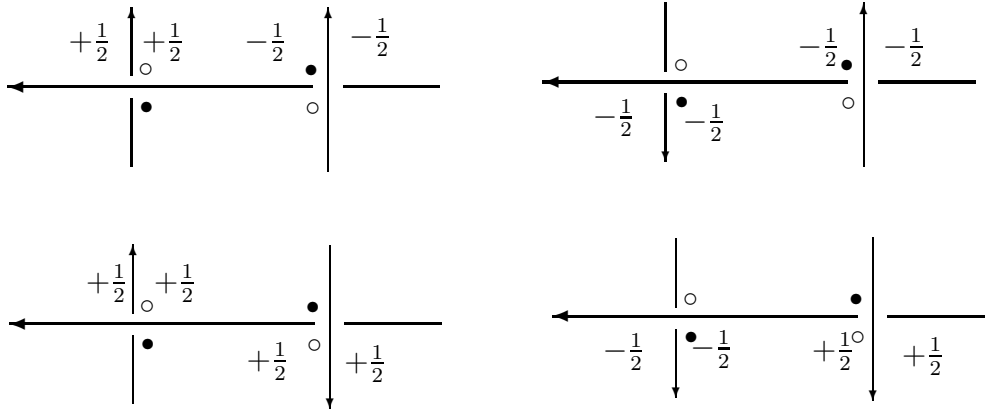
\setlength{\unitlength}{3947sp}%
\noindent
{\bf Case III:} For the other cases with three distinct strands, the 
strand on the right should go 
under the horizontal strand. However, if we rotate the figures in cases 
(1) and 
(2) $180^{\circ}$ using the horizontal strand as an axis, we get precisely 
those cases. The disc also rotates and occurs ``beneath" the Heegaard 
diagram. This disc represents a counter- clock move; However, we will 
still calculate the difference for a clock move from $\circ$ to $\bullet$. 
For this we use $\phi^{-1}$ with $\mathcal{D}(\phi^{-1})=$ 
$-\mathcal{D}(\phi)$.  The weights on the horizontal strands reflect 
across the 
horizontal strand as do the weights for the grading. The weights for the 
thin 
strands remain the same. Thus $\mathcal{F}_{i}^{III}(\bullet) = 
\mathcal{F}_{i}^{II}(\circ)$ and $\mathcal{F}_{i}^{III}(\bullet)  - 
\mathcal{F}_{i}^{III}(\circ) =$ $-\big( \mathcal{F}_{i}^{II}(\bullet) 
- \mathcal{F}_{i}^{II}(\circ) \big )$ $= +1$ for those configurations in 
case II. But
  $(n_w - n_{z_i})(\phi_{III}^{-1}) =$ $- (n_w - n_{z_i})(\phi_{III}) =$ 
$- (n_w - 
n_{z_i})(\phi_{II}) = + 1$ as $\phi_{II}$ and $\phi_{III}$ include the 
$i^{th}$ meridian the same number of times. 
\\
\\
{\bf Case IV:} In the cases where two or more of the above strands are the 
in the same component, we  
employ the following observations: 1) if the two thin strands belong to 
the 
same link component, then 
nothing changes, and 2) if the horizontal strand corresponds to the same 
component 
as one of the thin strands 
(or both), the sum of the weights in each quadrant differs by the same 
amount from that quadrant's weight 
as a self-crossing. Thus the difference between intersection points of the 
sum is
the same as the difference of the self-intersection weights. The grading 
computations don't
change. Inspecting the values in Figures \ref{fig:casei} and 
\ref{fig:caseii}, we see that the
filtration difference still equals the difference in the weights for the 
horizontal strand. 
\\
\\
Finally, we have implicitly assumed that the horizontal strand between the
two intersections is locally unknotted. Local knotting alters the topology 
of the domain $\mathcal{D}(\phi)$ above. Take the square in case I. If we 
knot the horizontal strand,
there is still a class $\phi$ joining the two intersection points, but it 
is a punctured disc with the same four points on its outer boundary, and 
one point that is both 
$\bullet$ and $\circ$ on each of its other boundaries. 
These new boundaries come from the faces in the
projection of the local knot, and consist entirely of $\be$ curves. The 
$\al$
curve at each intersection point on such a boundary joins that boundary to
another boundary $\be$, possibly from the square. Together these form a
forest with vertices being $\be$ curves, and edges being the $\al$'s
that join them. In the Appendix we show that these still have $\mu(\phi) = 
1$,which is all we need in the proof.

\subsection{Heegaard-Floer Gradings and Grading Weights}

In the remainder of this section we connect the gradings from 
Heegaard-Floer homology and those from the weights. It should come as no 
surprise that they are equal and that the argument in \cite{Alte} adapts 
readily.
For the grading, we know that ${\bf gr}({\bf y}) - {\bf gr}({\bf x})$ is 
the same as the difference in weights.
However, we will show that there is an intersection point for which the 
sum of the weights above is the absolute grading obtained from the 
Heegaard-Floer theory. Thus the sum of the grading weights will equal the 
grading for every intersection point since the difference between these is 
equal for distinct intersection points. We will
need a lemma before we proceed.
\\
\begin{lemma} \label{lemma:uniqueness} Let $G \subset \R^2$ be a finite 
graph and let $H_{\al}$ be the handlebody that is its regular neighborhood 
in $S^3$. Then $S^3 - H_{\al}$ is a handlebody
with the co-cores of its one-handles corresponding to the bounded faces of 
$\R^2 - G$. We
choose these co-cores to be the $\setb$ of a Heegaard diagram for $S^3$. 
Suppose that $\seta$
contains $\al$'s which intersect at most $2$ $\be$'s, each geometrically 
once. Furthermore, 
assume that each $\al$ links exactly one edge of $G$. Then there is only 
one point in 
$\T_{\al} \cap \T_{\be}$.
\end{lemma}

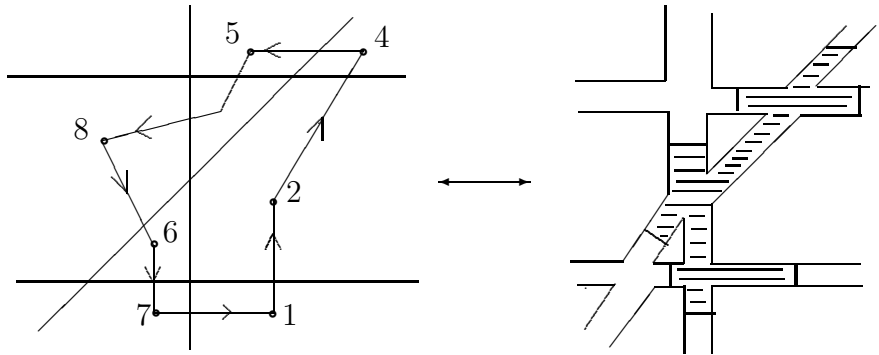
\begin{figure}
\begin{center}
\begin{picture}(5477,2212)(585,-2236)
\thinlines
\put(3301,-1149){\vector(-1, 0){  0}}
\put(3301,-1149){\vector( 1, 0){575}}
\thicklines
\put(2830,-328){\circle{38}}
\put(2124,-328){\circle{38}}
\thinlines
\put(4846,-1801){\line( 0,-1){423}}
\thicklines
\put(1203,-890){\circle{38}}
\put(2265,-1271){\circle{38}}
\put(1518,-1538){\circle{38}}
\put(2261,-1974){\circle{38}}
\put(1528,-1971){\circle{38}}
\thinlines
\put(1739,-36){\line( 0,-1){2162}}
\put(1076,-486){\line( 1, 0){2021}}
\put(1151,-1774){\line( 1, 0){2021}}
\put(2764,-111){\line(-1,-1){1974}}
\put(651,-1774){\line( 1, 0){517}}
\put(1114,-486){\line(-1, 0){517}}
\multiput(1495,-738)(-6.26667,-6.26667){16}{\makebox(1.6667,11.6667){\SetFigFont{5}{6}{rm}.}}
\multiput(1401,-832)(7.83333,-3.91667){13}{\makebox(1.6667,11.6667){\SetFigFont{5}{6}{rm}.}}
\multiput(1250,-1124)(6.26667,-6.26667){16}{\makebox(1.6667,11.6667){\SetFigFont{5}{6}{rm}.}}
\put(1344,-1218){\line( 0, 1){141}}
\multiput(1460,-1689)(3.91667,-7.83333){13}{\makebox(1.6667,11.6667){\SetFigFont{5}{6}{rm}.}}
\multiput(1507,-1783)(3.91667,7.83333){13}{\makebox(1.6667,11.6667){\SetFigFont{5}{6}{rm}.}}
\multiput(1951,-1926)(6.71429,-6.71429){8}{\makebox(1.6667,11.6667){\SetFigFont{5}{6}{rm}.}}
\multiput(1994,-1969)(-6.71429,-6.71429){8}{\makebox(1.6667,11.6667){\SetFigFont{5}{6}{rm}.}}
\multiput(2262,-1479)(-3.91667,-7.83333){13}{\makebox(1.6667,11.6667){\SetFigFont{5}{6}{rm}.}}
\multiput(2215,-1573)(3.91667,7.83333){13}{\makebox(1.6667,11.6667){\SetFigFont{5}{6}{rm}.}}
\multiput(2262,-1479)(3.91667,-7.83333){13}{\makebox(1.6667,11.6667){\SetFigFont{5}{6}{rm}.}}
\multiput(2305,-278)(-7.83333,-3.91667){13}{\makebox(1.6667,11.6667){\SetFigFont{5}{6}{rm}.}}
\multiput(2211,-325)(7.83333,-3.91667){13}{\makebox(1.6667,11.6667){\SetFigFont{5}{6}{rm}.}}
\put(2551,-782){\makebox(1.6667,11.6667){\SetFigFont{5}{6}{rm}.}}
\put(2573,-887){\line( 0, 1){141}}
\multiput(2573,-746)(-6.26667,-6.26667){16}{\makebox(1.6667,11.6667){\SetFigFont{5}{6}{rm}.}}
\put(2830,-325){\line(-1, 0){705}}
\multiput(2125,-325)(-3.76000,-7.52000){51}{\makebox(1.6667,11.6667){\SetFigFont{5}{6}{rm}.}}
\put(1937,-701){\line(-4,-1){752}}
\put(1185,-889){\line( 1,-2){329}}
\put(1514,-1547){\line( 0,-1){423}}
\put(1514,-1970){\line( 1, 0){752}}
\put(2266,-1970){\line( 0, 1){705}}
\put(2266,-1265){\line( 3, 5){564}}
\put(5020,-88){\line( 0,-1){470}}
\put(5020,-558){\line( 1, 0){470}}
\put(4225,-1820){\makebox(1.6667,11.6667){\SetFigFont{5}{6}{rm}.}}
\put(5013,-1294){\line( 1, 1){564}}
\put(5577,-730){\line( 1, 0){376}}
\put(5007,-1666){\line( 1, 0){940}}
\put(5020,-1801){\line( 0,-1){423}}
\put(5026,-1801){\line( 1, 0){940}}
\put(5674,-543){\line( 1, 1){376}}
\put(5681,-550){\line( 1, 0){329}}
\put(4180,-707){\line( 1, 0){564}}
\put(4744,-707){\line( 0,-1){517}}
\put(4744,-1224){\line(-2,-3){282}}
\put(4462,-1647){\line(-1, 0){329}}
\put(4174,-1794){\line( 1, 0){235}}
\multiput(4409,-1794)(-4.70000,-7.05000){41}{\makebox(1.6667,11.6667){\SetFigFont{5}{6}{rm}.}}
\put(5035,-719){\line( 1, 0){329}}
\put(5364,-719){\line(-1,-1){376}}
\put(4988,-1095){\line( 0, 1){376}}
\put(4988,-719){\line( 1, 0){ 94}}
\multiput(4840,-1377)(-4.70000,-7.05000){41}{\makebox(1.6667,11.6667){\SetFigFont{5}{6}{rm}.}}
\put(4652,-1659){\line( 1, 0){141}}
\put(4793,-1659){\line( 1, 0){ 47}}
\put(4840,-1659){\line( 0, 1){282}}
\put(5485,-552){\line( 1, 1){235}}
\multiput(5720,-317)(6.13043,6.13043){24}{\makebox(1.6667,11.6667){\SetFigFont{5}{6}{rm}.}}
\thicklines
\put(5180,-556){\line( 0,-1){141}}
\thinlines
\put(5405,-729){\line( 1, 0){ 94}}
\put(5283,-845){\line( 1, 0){ 94}}
\put(5174,-954){\line( 1, 0){ 94}}
\put(5655,-409){\line( 1, 0){ 94}}
\put(4378,-2186){\line( 3, 4){282}}
\put(4660,-1810){\line( 1, 0){188}}
\put(5533,-550){\line( 1, 0){ 94}}
\put(4879,-1364){\line( 1, 0){ 94}}
\put(4866,-1826){\line( 1, 0){ 94}}
\put(4885,-1903){\line( 1, 0){ 94}}
\put(4693,-1486){\line( 1, 0){ 47}}
\put(5020,-1284){\line( 0,-1){376}}
\put(4802,-1698){\line( 1, 0){658}}
\put(4885,-1621){\line( 1, 0){ 94}}
\put(4885,-1531){\line( 1, 0){ 94}}
\put(4892,-1441){\line( 1, 0){ 94}}
\put(4814,-1756){\line( 1, 0){658}}
\put(4679,-1429){\line( 1, 0){ 94}}
\put(4699,-1358){\line( 1, 0){ 94}}
\put(4725,-1294){\line( 1, 0){282}}
\put(4892,-1204){\line( 1, 0){188}}
\put(4769,-1204){\line( 1, 0){141}}
\put(4796,-1146){\line( 1, 0){329}}
\put(4783,-1076){\line( 1, 0){188}}
\put(4783,-992){\line( 1, 0){188}}
\put(5148,-999){\line( 1, 0){ 94}}
\put(5103,-1044){\line( 1, 0){ 94}}
\put(5058,-1089){\line( 1, 0){ 94}}
\put(5232,-896){\line( 1, 0){ 94}}
\put(5346,-781){\line( 1, 0){ 94}}
\put(5251,-672){\line( 1, 0){611}}
\put(5238,-614){\line( 1, 0){658}}
\put(5591,-479){\line( 1, 0){ 94}}
\put(4727,-41){\line( 0,-1){470}}
\put(4727,-511){\line(-1, 0){564}}
\put(5713,-351){\line( 1, 0){ 94}}
\thicklines
\multiput(4597,-1454)(7.05000,-4.70000){21}{\makebox(6.6667,10.0000){\SetFigFont{7}{8.4}{rm}.}}
\put(4757,-1660){\line( 0,-1){141}}
\put(4853,-1974){\line( 1, 0){188}}
\put(5738,-300){\line( 1, 0){188}}
\put(5944,-543){\line( 0,-1){188}}
\put(5546,-1659){\line( 0,-1){141}}
\put(4757,-915){\line( 1, 0){235}}
\put(2345,-1273){\makebox(0,0)[lb]{\smash{{\SetFigFont{12}{14.4}{rm}2}}}}
\put(2897,-300){\makebox(0,0)[lb]{\smash{{\SetFigFont{12}{14.4}{rm}4}}}}
\put(1959,-274){\makebox(0,0)[lb]{\smash{{\SetFigFont{12}{14.4}{rm}5}}}}
\put(1009,-878){\makebox(0,0)[lb]{\smash{{\SetFigFont{12}{14.4}{rm}8}}}}
\put(1569,-1522){\makebox(0,0)[lb]{\smash{{\SetFigFont{12}{14.4}{rm}6}}}}
\put(1399,-2044){\makebox(0,0)[lb]{\smash{{\SetFigFont{12}{14.4}{rm}7}}}}
\put(2319,-2044){\makebox(0,0)[lb]{\smash{{\SetFigFont{12}{14.4}{rm}1}}}}
\end{picture}%
\end{center}
\caption{The Contradiction  in Lemma \ref{lemma:uniqueness} Arising from a 
Cycle.}\label{fig:contract}
\end{figure}

{\bf Proof:} Supppose ${\bf x}$ and ${\bf y}$ are points in $\T_{\al} \cap 
\T_{\be}$ and 
${\bf x} \neq {\bf y}$. Then ${\bf x} = \{ x_1, \ldots, x_g \}$ and ${\bf 
y} = \{y_1, \ldots,
y_g\}$ where $x_i \in a_{\sigma(i)} \cap \be_{i}$ and $y_i \in a_{\psi(i)} 
\cap \be_{i}$ with
$\sigma, \psi \in S_n$. If ${\bf x} \neq {\bf y}$ then $\psi^{-1} \circ 
\sigma$
is not the identity. It must therefore have a decomposition into cycles 
with at least
one of length greater than or equal to $2$. In the
planar graph formed by placing a vertex in each bounded face of $\Gamma$, 
and
an edge between each pair of faces abutted by an $\al$ curve, each 
non-trivial cycle in the cycle decomposition corresponds to a 
cycle in the graph. Each cycle in the graph implies the existence of a 
collection of $\al$'s that are null-homologous in
$\Sigma$, contradicting the Heegaard assumption (see Figure 
\ref{fig:contract}). Thus
at most one intersection point exits. Since $\widehat{HF}(S^3) \simeq 
\Z_{(0)}$, there
is at least one intersection point, ${\bf x}_{0}$.
\\
\\

Starting with a string link in $S^3$, we use handleslides to construct a 
Heegaard diagram 
as in the lemma. The intersection point in this Heegaard diagram will 
correspond to one in the 
diagram for the string link. We will ensure that the handleslides do 
not alter the absolute grading. The  intersection 
point for the string link must then have absolute grading equal to 
$0$. This will also be the value of $G({\bf x})$ for
that intersection point. Since we know that $G({\bf y}) - G({\bf x})$ 
$= gr({\bf y}) - gr({\bf y})$, the weights above will give the absolute 
grading of each intersection point. 
\\
\\
To obtain an acceptable diagram, it suffices  that the new $\al$  
link an edge entering the crossing defined by the old $\al$ before 
the handlesliding as the weights for the grading only occur in the 
quadrant between 
the two edges exiting a crossing. The unique intersection point will then 
correspond to 
an old intersection point with $G({\bf x}) = 0$. We ignore the points
$\{ z_1, \ldots, z_k \}$ to find a pointed diagram for $S^3$. 

We order the crossings in the string link projection by the following 
conditions.
Each crossing of $L_{k}$ is larger than any crossing of strands $L_{i}$ 
and $L_{j}$
with $i, j < k$. For any $k$, the crossings with $L_{i}$ for $i \leq k$
are enumerated from largest to smallest
by the first time they are encountered while travelling backwards along 
$L_{k}$ from
the meridian. We adjust each crossing in increasing order by
starting at that crossing (either for $L_j$ with itself or between $L_i$ 
and 
$L_j$ with $i \leq j$) and isotoping and
handelsliding the $\al$-curve in the direction of $L_{j}$, going 
over all the $\al$'s along that route. At self-crossings we choose to 
follow the edge which exits the crossing and arrives at a meridian without 
returning to the crossing.  When we arrive at the meridian we handleslide 
across it, and then repeat 
the procedure in reverse. This produces an $\al$ linking the penultimate 
edge through the original crossing. Furthermore, for each crossing the 
ordering
implies that there is a path to the meridian along the orientation of one 
or 
other strand, along which the crossing does not recur, and such that none 
of 
the $\al$'s encountered have been previously alterred. 
\\
\\
The handleslides do not take the $\al$'s across $w$. By inspecting the 
standard handleslide diagram, \cite{Hom3}, we can see that the new 
intersection
point is the in the image of one of the original intersection points, the 
one 
with marking in the same quadrant at each crossing,  under the composition 
of the handleslide maps. The cobordism induced by the handleslides is 
$S^3 \times I$, so the formula calculating the change in absolute grading 
imples that the grading does not change. Moreover, since the triangle does 
not 
cross $w$ we know that the image is $[{\bf x_0}, 0]$ which has absolute 
grading zero. This equals the absolute grading assigned by the weights.
\\
\\
{\bf Note:} 
For a single knot in $S^{3}$, a string link consists solely in the 
choice of a point on the knot. This gives a preferred position for the 
meridian and the two points $w, z$ used in \cite{Knot} to calculate the 
knot Floer homology. The weights above agree with those of \cite{Alte}.

\subsection{Example}

\setlength{\unitlength}{3158sp}%
\begin{figure}
\begin{center}
\begin{picture}(7356,3229)
\put(5431,1625){\makebox(0,0)[lb]{\smash{{\SetFigFont{12}{14.4}{rm}$\ddots$}}}}
\put(1223,1625){\makebox(0,0)[lb]{\smash{{\SetFigFont{12}{14.4}{rm}$\ddots$}}}}
\includegraphics{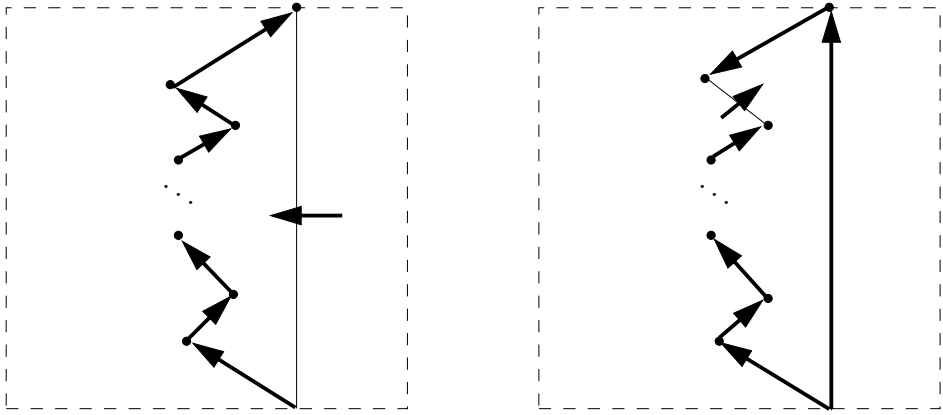}%
\end{picture}%
\end{center}
\caption{Examples of dual maximal forests for the graph $\Gamma$ from 
Figure \ref{fig:meridians}}
\label{fig:mtn}
\end{figure}
\setlength{\unitlength}{3947sp}%

In the graph of black regions for Figure \ref{fig:meridians} we may form 
dual forests in
two ways, see Figure \ref{fig:mtn}. First, we may have arrows pointing
up the segments on the left side of the graph of black regions. Second, we 
may have
some arrows pointing up that side, then a transverse arrow, then arrows 
pointing
down the remainder of the segments on the left side. To join the vertices 
at top and bottom
we have an arrow up the single segment on the right side.

When all the arrows go up the left side, they place the 
local intersection
points on the opposite side of $L_{2}$ from that pointed to by the 
orientation of $L_{1}$.
Looking at the weights indicates that, regardless of the crossing type, 
such a tree contibutes
$0$ to the filtration index for $L_{2}$ and $0$ to the grading weight. On 
the other
hand, when the crossings are positive all the segments on the left are 
assigned
$-\frac{1}{2}$ for the weight on $L_{1}$. The transverse arrow on the 
right will
contribute nothing to either sum. A similar analysis for negative linking 
implies
that this forest contributes $h_{1}^{-\,lk(L_{1}, L_{2})}$ to the Euler 
characteristic.

For the second type of forest, let $m$
be the number of arrows pointing up on the left side. Then $m$ can vary 
between
$m=0$ and $m=2\,lk(L_{1}, L_{2}) - 1$ $= 2\,L - 1$ when the linking number 
is positive. The monomial for this state is $h_{1}^{-s}h_{2}^{- L + s}$ 
when $m = 2\,s$ and $-h_{1}^{-s - 1} h_{2}^{-L + s}$ when $m = 2\,s + 1$. 
The minus sign in the latter comes from the generator having grading $-1$. 
Together these imply a polynomial of the form

$$
\big(h_{1}^{- L} + h_{1}^{-L + 1}h_{2} + \cdots + h_{2}^{-L}\big) - \big 
(h_{1}^{-L}h_{2}^{-1} + \cdots h_{1}^{-1}h_{2}^{-L} \big)
$$

Since the string link is alternating, and the minus signs occur from 
grading $-1$, we can
determine the homology for the string link from the above polynomial (see 
the next section). It
is ($v_{1}, v_{2} \leq 0$)

$$
\begin{array}{cc}
\widehat{HF}(S, (v_{1}, v_{2})) \cong \Z_{(0)} & v_{1} + v_{2} = 
-lk(L_{1}, L_{2}) 
\end{array}
$$

$$
\begin{array}{cc}
\widehat{HF}(S, (v_{1}, v_{2})) \cong \Z_{(-1)} & v_{1} + v_{2} = 
-lk(L_{1}, L_{2}) - 1
\end{array}
$$

\subsection{Braids and Alternating String Links}

\subsubsection{Triviality for Braids}
In section \ref{sec:Euler}, we saw that the Euler characteristic of the 
chain groups for a string link in $S^{3}$ produces the torsion of the 
string link $\tau(S_{L})$ defined in \cite{Livi}.  
The torsion is known to be trivial if  the string link is a pure braid. 
Analogously, we may prove:

\begin{lemma} $\widehat{HF}(S_{L}) \cong \Z_{0}$ when $S_{L}$ is a pure 
braid. Using the weights from the previous section, the homology is 
non-trivial
only for the index $(0, \ldots, 0,0)$
\end{lemma}

{\bf  Proof:} There are two ways to prove this statement. The first 
analyzes the combinatorics of Kauffman states in the diagram. There is 
only one state:
The meridians consume the first row of $\be$-curves. Each 
crossing then has three regions already claimed, so it must use the 
fourth. This consumes the second set, and we proceed up the diagram.  
However, a more conceptual 
explanation may also be given. We know that the invariant we have defined 
does not depend upon how the 
strands move about on $D^{2} \times \{1\}$. We may simply undo the braid 
to obtain the trivial string link. 
 
The statement about filtration indices follows from the weights
defined in the previous section and the observation that the unique state 
assigns
its local contributions to the quadrant between two edges pointing into 
the
crossing as the strands are oriented down the page.

\subsubsection{Vanishing Differential for Alternating String Links}
\label{sec:Altern}

We call a string link, $S$, {\em alternating} if there is a projection of 
$S$
where proceeding along any strand in $S$ from $D^{2} \times \{0\}$ to 
$D^{2} \times \{1\}$ encounters alternating over and under-crossings. 
In a projection, we may place a small kink (formed using the first
Reidemeister move) in any strand which does not initially
participate in any crossing (self or otherwise) without
changing whether the projection is alternating. 

We call the $i^{th}$ 
strand, $s_{i}$
when counting from left to right on $D^{2} \times \{0\}$ and $s^{i}$ when
counting on $D^{2} \times \{1\}$. If we follow the $s_{i}$ from $D^{2} 
\times \{0\}$
to $D^{2} \times \{1\}$ find that it is $s^{\sigma(i)}$ for some 
permutation 
$\sigma \in S_{k}$.  

At each end of $D^{2} \times I$ we may label the strands as $u$ or $o$
for whether the strand is the over or under strand in the crossing
immediately preceding (following) the end of the strand on $D^{2} \times 
\{1\}$ ($D^{2} \times \{0\}$). This assigns a $k$-tuple called the {\em 
trace} of the string link at that end. We denote the trace at $D^{2} 
\times \{i\}$ by $T_{i}(S)$.  
There is an inversion on such $k$-tuples found by interchanging
$u$'s and $o$'s, which we denote $v \ra \overline{v}$. To compose 
alternating string 
links to have an alternating result requires $T_{1}(S_{1}) = 
\overline{T_{0}(S_{2})}$.
\\
\\
Our goal is to show

\begin{thm}
Let $S$ be a string link with alternating projection. Then the chain 
complex, $\widehat{CF}(S, \overline{j})$, arising from the projection 
Heegaard splitting has trivial differential for every index 
$\overline{j}$. In fact, all the generators have the same grading.
\end{thm}
\noindent
This generalizes the result in \cite{Alte} for alternating knots. We will
need the knot case for the general result. The result in \cite{Alte} is 
somewhat
stronger (by identifying the grading).
\\
\ \\
{\bf Note:} We call Morse critical points of index $1$ for the projection 
of $s_{i}$ to the $I$ factor a ``cap''. Critical points of index $0$ are 
called ``cups''.

\begin{lemma}
We may draw an alternating projection for $S$ so that 1) Every crossing 
occurs with
both strands oriented up, 2) No two crossings, caps, or cups occur in the 
same level in the $I$
factor.
\end{lemma}

\noindent
{\bf Proof:} Rotate every crossing without the correct orientations 
so that both strands go up. This can be done in a small neighborhood of 
the 
crossing at the expense of introducing cups and caps.
The second condition is achieved by a small perturbation of the diagram. 
\\
\\
Given a diagram in this form, we proceed with a few combinatorial 
lemmas. These are devoted to showing that an alternating string link
may be closed up to give an alternating knot with certain additional 
properties. At each time $t$ in the parametrization of the $i^{th}$ 
strand, let $f_{i}(t)$ be the number of strands strictly to the left of 
that 
point on $s_{i}$.

\begin{lemma} 
Let $t$ be a time when $s_{i}(t)$ is in a level (in $I$) not containing 
any caps or cups, and
not at a crossing. The total number of crossings encountered along the 
strand $s_{i}$ by the time $t$ is $\equiv$ $|f_{i}(t) - i + 1| {\mathrm 
mod\ } 2$ when $s_{i}$ is oriented up and
$\equiv |f_{i}(t) - i| {\textrm mod\ } 2$ when $s_{i}$ is oriented down. 
\end{lemma}

\noindent
{\bf Proof:} The number of caps plus the number of cups encountered in the 
$i^{th}$ strand, by time $t$, is even when the strand is oriented up, and 
odd when the strand is oriented down. $f_{i(t)}$ changes value by $1$
as $s_{i}(t)$ goes through a cap, cup or crossing. It changes value by $2$
at levels where a cap or cup occurs to the left of $s_{i}(t)$. By reducing
modulo $2$ we eliminate the latter variation. Since there are
$f_{i}(t)$ strands to the left, having started with $i - 1$ strands to the 
left, 
the number of cups, caps and crossings must be congruent to $|f_{i}(t) - i 
+ 1|$ 
modulo $2$. Removing the parity of the number of caps and cups gives the 
result. 

\begin{cor}
The total number of crossings encountered by the $i^{th}$ strand is 
congruent modulo $2$ to
$|\sigma(i) - i|$.
\end{cor}

\noindent
This lemma has the following consequence:
\begin{lemma}
Suppose $s_{i}$ and $s_{j}$ cross somewhere in $S$. If $T_{0}(s_{i}) = 
T_{0}(s_{j})$ then
$i \equiv j\ {\textrm mod\ }2$. If $T_{0}(s_{i}) \neq T_{0}(s_{j})$ then 
$i \equiv j\,+\,1\ {\textrm mod\ }2$
\end{lemma}

\noindent
{\bf Proof:} Consider the first time that they cross in the ordering on 
$s_{i}$. Suppose that $i < j$, and that there are $k$ strands to the left 
of the point in $s_{i}$ just before the crossing. Then there must be $k 
\pm 1$ strands to the left of $s_{j}$. We label each point on the strands, 
except at crossings, by a $u$ or an $o$ depending
upon whether an over, or under, crossing must occur next. The labels of 
the points on the two strands just before the crossing of $s_{i}$ and 
$s_{j}$ must be different. If $s_{i}$ has encountered an even number of 
crossings prior to this point, it will have label $T_{0}(s_{i})$, 
otherwise it has label $\overline{T_{0}(s_{i})}$. The same will be true of 
$s_{j}$. We have assumed that both strands are oriented up just before the 
crossing. Thus, the parity of the number of crossings involving $s_{i}$ is 
that of $|f_{i}(t_{i}) - i + 1|$ $ = |k - i + 1|$ and the same parity for 
$s_{j}$ is $|f_{j}(t_{j}) - j + 1|$ $\equiv | k - j |$. If $T_{0}(s_{i}) = 
T_{0}(s_{j})$, then one strand must have experienced an even number of 
crossings, while the other experienced an odd number. This happens when $i 
\equiv j$. Otherwise, both strands must enounter the same parity of 
crossings and $i \equiv j + 1$.
\\
\\
We decompose $S = S_{1} \cup \cdots \cup S_{l}$, where $S_{j}$ consists of 
a maximal string link with connected projection. We apply the following 
lemma
to each $S_{j}$.

\begin{lemma}
For a connected, alternating string link $T_{0}(S)$ must be either $(u, o, 
u, \ldots)$ with
alternating entries, or $(o, u, o, \ldots)$
\end{lemma}

\noindent
{\bf Proof:} For $s_{i}$ and $s_{j}$ there is a sequence $s_{i_{0}}, 
\ldots, s_{i_{r}}$ with
$s_{i_{0}} = s_{i}$ and $s_{i_{r}} = s_{j}$ and where consecutive entries 
cross one another. The result follows from induction using the conclusion 
of the previous lemma, which also proves the base case. 
\\
\\
Since any strand in $S_{j}$ divides the projection into two parts, we see 
that each $S_{j}$ must consist of consecutive strands in the diagram 
(along both ends). Furthermore, we have shown    

\begin{lemma}
For an alternating projection of $S$, $T_{0}(S) = \overline{T_{1}(S)}$.
\end{lemma}

\noindent
This lemma guarantees that the usual closure of the string link (join 
$s^{i}$ to $s_{i}$ for all $i$) is alternating.
\\
\\
{\bf Proof:} We divide $S = S_{1} \cup \cdots \cup S_{l}$ and apply the 
corollary above to
each maximal sub-string link. By the preceding lemma, we have an 
alternating trace for the
end $D^{2} \times \{0\}$. By the corollary, the other end of $s_{i}$ has 
the same label
when $\sigma(i) - i$ is even, and different labels when $\sigma(i) - i$ is 
odd. The result follows directly.
\\
\\
In fact, a repetition of $u$ or $o$ in $T_{0}(S)$ implies that $S = S_{1} 
\cup S_{2}$ for an alternating string link. $S_{1}$ consists of those 
strands including 
and to the left of the first $u$, and $S_{2}$ consists of those strands 
including and to the right of the second. 
\\
\\ 
{\bf Note:} Recall that we add kinks by the first Reidemeister move to 
strands who don't cross
any other strand (including themselves). This is how they get labelled.
\\
\\
The chain complex of a string link $S$ which decomposes as $S_{1} \cup 
\cdots \cup S_{l}$
for the index $(\overline{j}_{1}, \cdots, \overline{j}_{l})$ is  
$\widehat{CF}(S_{1}, \overline{j}_{1}) \otimes \cdots \otimes 
\widehat{CF}(S_{2}, \overline{j}_{l})$ with the standard tensor product 
differential. Thus, proving the theorem for connected, alternating 
string links will prove the general result. Our strategy is to compare 
the chain group from out projection to the chain group of an alternating 
knot.
\\
\\
We make a few observations about braids.
First, their projection Heegaard diagrams possess only one 
generator. 
Second, given the projection of a braid, forgetting over and under 
crossings, there is precisely one set of crossing data with $T_{0}(B) = 
(u, o, u, o, \ldots)$. Write the braid as a product of generators or their 
inverses. The traces picks out which (generator or inverse) must occur. 
Pushing up the diagram, the traces $(u, o, u, o, \ldots)$ repeats as 
$T_{0}(B')$ for the remainder of the braid. 

\begin{figure}
\begin{center}
\begin{picture}(7053,2779)(389,-2999)
\thinlines
\put(5646,-2075){\line( 0,-1){705}}
\put(5646,-2780){\line( 1, 0){940}}
\put(6586,-2780){\line( 0, 1){141}}
\put(6738,-1423){\line( 0, 1){423}}
\put(6841,-1423){\line( 0, 1){423}}
\put(7114,-1423){\line( 0, 1){423}}
\put(6730,-2388){\line( 0, 1){423}}
\put(6832,-2380){\line( 0, 1){423}}
\put(7114,-2380){\line( 0, 1){423}}
\put(4459,-663){\line( 0, 1){423}}
\put(4459,-240){\line( 1, 0){2773}}
\put(7232,-240){\line( 0,-1){517}}
\put(7232,-757){\line( 0, 1){  0}}
\put(7233,-1423){\line( 0, 1){705}}
\put(2401,-1961){\framebox(1075,537){}}
\put(2622,-1423){\line( 0, 1){705}}
\put(2784,-1415){\line( 0, 1){423}}
\put(2887,-1415){\line( 0, 1){423}}
\put(3160,-1415){\line( 0, 1){423}}
\put(2776,-2380){\line( 0, 1){423}}
\put(2878,-2372){\line( 0, 1){423}}
\put(3160,-2372){\line( 0, 1){423}}
\put(1726,-638){\line( 1, 0){893}}
\put(2619,-638){\line( 0,-1){141}}
\put(1666,-2081){\line( 0,-1){564}}
\put(1666,-2645){\line( 1, 0){987}}
\put(2648,-2639){\line( 0, 1){658}}
\put(505,-655){\line( 0, 1){423}}
\put(505,-232){\line( 1, 0){2773}}
\put(3278,-232){\line( 0,-1){517}}
\put(3278,-749){\line( 0, 1){  0}}
\put(539,-2047){\line( 0,-1){799}}
\put(539,-2846){\line( 1, 0){2726}}
\put(3265,-2846){\line( 0, 1){188}}
\put(3279,-1415){\line( 0, 1){705}}
\put(2690,-2175){\line(-1, 0){ 94}}
\put(1491,-1582){\line( 1, 1){235}}
\put(1726,-1347){\line( 0, 1){705}}
\put(539,-2038){\line( 0, 1){423}}
\multiput(539,-1615)(3.81081,7.62162){38}{\makebox(1.6667,11.6667){\SetFigFont{5}{6}{rm}.}}
\put(680,-1333){\line( 1, 0){ 94}}
\put(505,-664){\line( 0,-1){564}}
\multiput(505,-1228)(3.76000,-7.52000){26}{\makebox(1.6667,11.6667){\SetFigFont{5}{6}{rm}.}}
\put(1444,-1389){\line( 0,-1){705}}
\put(1017,-1364){\line( 0,-1){705}}
\put(1256,-1628){\line( 0,-1){423}}
\put(804,-1611){\line( 0,-1){470}}
\put(846,-1330){\line( 1, 0){376}}
\multiput(641,-1483)(3.91667,-7.83333){13}{\makebox(1.6667,11.6667){\SetFigFont{5}{6}{rm}.}}
\put(688,-1577){\line( 1, 0){282}}
\put(970,-1577){\line( 0, 1){  0}}
\put(1077,-1577){\line( 1, 0){329}}
\put(1256,-681){\line( 0,-1){846}}
\put(1444,-680){\line( 0,-1){611}}
\put(1017,-680){\line( 0,-1){611}}
\put(804,-681){\line( 0,-1){846}}
\put(1306,-1327){\line( 1, 0){188}}
\multiput(1494,-1327)(6.26667,-6.26667){16}{\makebox(1.6667,11.6667){\SetFigFont{5}{6}{rm}.}}
\put(1631,-1475){\line( 0, 1){  0}}
\multiput(1631,-1475)(6.71429,-6.71429){8}{\makebox(1.6667,11.6667){\SetFigFont{5}{6}{rm}.}}
\put(1678,-1522){\line( 0,-1){423}}
\put(5596,-798){\line( 0, 1){282}}
\put(5596,-516){\line( 1, 0){987}}
\put(6583,-516){\line( 0,-1){282}}
\put(6591,-1431){\line( 0, 1){705}}
\put(6355,-1969){\framebox(1075,537){}}
\put(6584,-2654){\line( 0, 1){658}}
\put(4411,-2047){\line( 0,-1){940}}
\put(4411,-2987){\line( 1, 0){2820}}
\put(7231,-2987){\line( 0, 1){329}}
\put(7233,-2678){\line( 0, 1){705}}
\put(3268,-2670){\line( 0, 1){705}}
\put(4455,-670){\line( 0,-1){423}}
\multiput(4455,-1093)(3.81081,-7.62162){38}{\makebox(1.6667,11.6667){\SetFigFont{5}{6}{rm}.}}
\put(4596,-1375){\line( 1, 0){ 94}}
\put(4421,-2044){\line( 0, 1){564}}
\multiput(4421,-1480)(3.76000,7.52000){26}{\makebox(1.6667,11.6667){\SetFigFont{5}{6}{rm}.}}
\put(5360,-1319){\line( 0, 1){705}}
\put(4933,-1344){\line( 0, 1){705}}
\put(5172,-1080){\line( 0, 1){423}}
\put(4720,-1097){\line( 0, 1){470}}
\put(4762,-1378){\line( 1, 0){376}}
\multiput(4557,-1225)(3.91667,7.83333){13}{\makebox(1.6667,11.6667){\SetFigFont{5}{6}{rm}.}}
\put(4604,-1131){\line( 1, 0){282}}
\put(4886,-1131){\line( 0, 1){  0}}
\put(4993,-1131){\line( 1, 0){329}}
\put(5172,-2027){\line( 0, 1){846}}
\put(5360,-2028){\line( 0, 1){611}}
\put(4933,-2028){\line( 0, 1){611}}
\put(4720,-2027){\line( 0, 1){846}}
\put(5222,-1381){\line( 1, 0){188}}
\multiput(5410,-1381)(6.26667,6.26667){16}{\makebox(1.6667,11.6667){\SetFigFont{5}{6}{rm}.}}
\put(5547,-1233){\line( 0, 1){  0}}
\multiput(5547,-1233)(6.71429,6.71429){8}{\makebox(1.6667,11.6667){\SetFigFont{5}{6}{rm}.}}
\put(5594,-1186){\line( 0, 1){423}}
\put(5410,-1149){\line( 1,-1){235}}
\put(5645,-1384){\line( 0,-1){705}}
\put(6611,-2201){\line(-1, 0){ 94}}
\put(6593,-1807){\makebox(0,0)[lb]{\smash{{\SetFigFont{12}{14.4}{rm}$S 
\cdot B$}}}}
\put(6935,-1244){\makebox(0,0)[lb]{\smash{{\SetFigFont{12}{14.4}{rm}...}}}}
\put(6892,-2225){\makebox(0,0)[lb]{\smash{{\SetFigFont{12}{14.4}{rm}...}}}}
\put(5953,-1551){\makebox(0,0)[lb]{\smash{{\SetFigFont{12}{14.4}{rm}$U$}}}}
\put(2639,-1799){\makebox(0,0)[lb]{\smash{{\SetFigFont{12}{14.4}{rm}$S 
\cdot B$}}}}
\put(2981,-1236){\makebox(0,0)[lb]{\smash{{\SetFigFont{12}{14.4}{rm}...}}}}
\put(2938,-2217){\makebox(0,0)[lb]{\smash{{\SetFigFont{12}{14.4}{rm}...}}}}
\put(1999,-1543){\makebox(0,0)[lb]{\smash{{\SetFigFont{12}{14.4}{rm}$U$}}}}
\put(728,-2209){\makebox(0,0)[lb]{\smash{{\SetFigFont{10}{12.0}{rm}$u$}}}}
\put(982,-2209){\makebox(0,0)[lb]{\smash{{\SetFigFont{10}{12.0}{rm}$o$}}}}
\put(1208,-2223){\makebox(0,0)[lb]{\smash{{\SetFigFont{10}{12.0}{rm}$u$}}}}
\put(1398,-2223){\makebox(0,0)[lb]{\smash{{\SetFigFont{10}{12.0}{rm}$o$}}}}
\put(389,-2103){\makebox(0,0)[lb]{\smash{{\SetFigFont{10}{12.0}{rm}$o$}}}}
\put(1638,-2082){\makebox(0,0)[lb]{\smash{{\SetFigFont{10}{12.0}{rm}$u$}}}}
\put(4262,-1900){\makebox(0,0)[lb]{\smash{{\SetFigFont{10}{12.0}{rm}$u$}}}}
\put(4694,-2161){\makebox(0,0)[lb]{\smash{{\SetFigFont{10}{12.0}{rm}$o$}}}}
\put(5149,-2172){\makebox(0,0)[lb]{\smash{{\SetFigFont{10}{12.0}{rm}$o$}}}}
\put(5637,-1991){\makebox(0,0)[lb]{\smash{{\SetFigFont{10}{12.0}{rm}$o$}}}}
\put(4910,-2172){\makebox(0,0)[lb]{\smash{{\SetFigFont{10}{12.0}{rm}$u$}}}}
\put(5365,-2161){\makebox(0,0)[lb]{\smash{{\SetFigFont{10}{12.0}{rm}$u$}}}}
\end{picture}%
\end{center}
\caption{The closures to construct an alternating knot from an alternating 
string link. The
case for an odd number of strands is similar. The basing (meridian) for 
the new knot is
shown as a small line across the knot and intersecting $U$. Of course, the
strands should meet at top and bottom as for the closure of a braid.}
\label{fig:alt}
\end{figure}
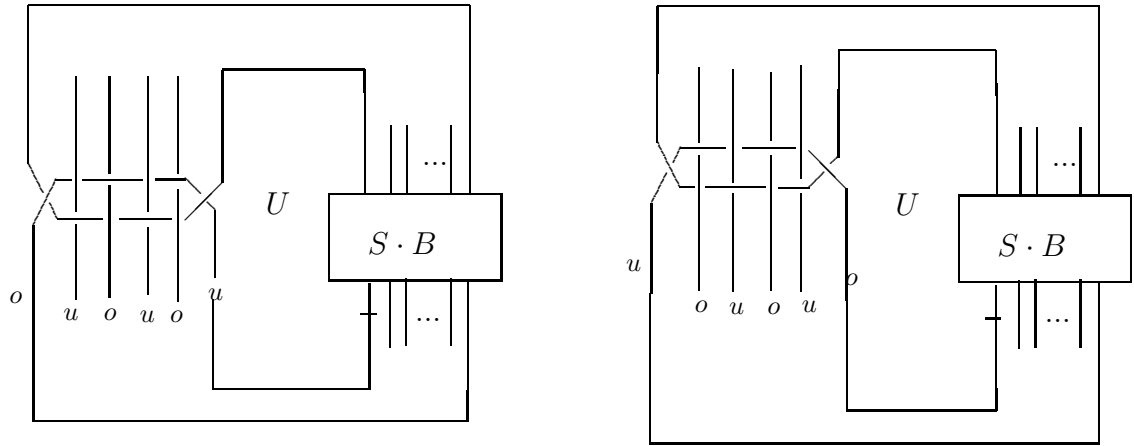

\begin{figure}
\begin{center}
\begin{picture}(4098,2929)
\put(3774,1392){\makebox(0,0)[lb]{\smash{{\SetFigFont{12}{14.4}{rm}$U$}}}}
\put(3037,479){\makebox(0,0)[lb]{\smash{{\SetFigFont{12}{14.4}{rm}$m_1$}}}}
\includegraphics{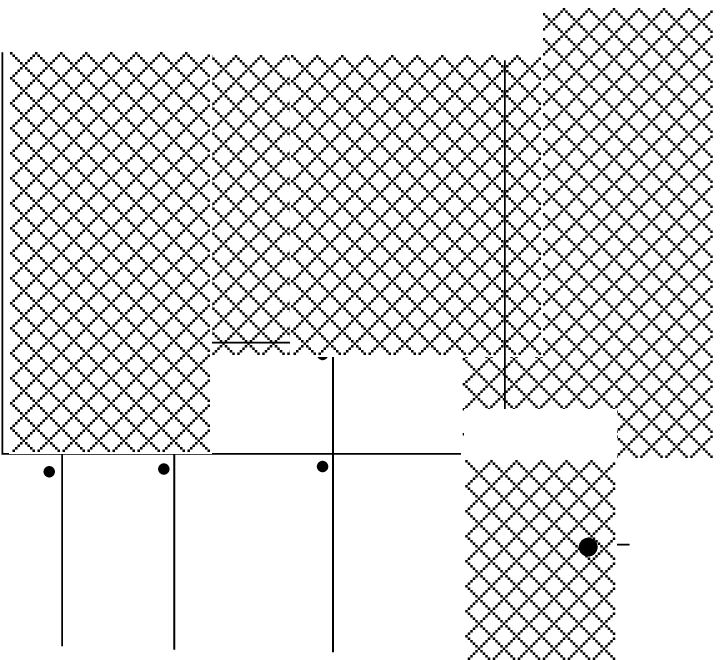}
\end{picture}%
\end{center}
\caption{The Unique Extension of Generators. The shaded regions have 
already been 
assigned to a crossing by the generator in $S \# B$. This corresponds to 
the extension of dual pairs used to extend the clock theorem to string 
links.}
\label{fig:State}
\end{figure}

\noindent
We choose a braid representing a permutation $\tau$ such that $\tau \circ 
\sigma$ is a cyclic
permutation taking $1$ to $k$. As we have just seen, we may use $T_{0}(B) 
= \overline{T_{1}(S)}$ to choose a braid so that $S \# B$ is still an 
alternating string link. We complete the construction by closing the 
new string link as in one of the diagrams in Figure \ref{fig:alt}. The 
trace at the bottom of $S \# B$ determines which to use. Due to our 
assumption about
$\tau$, the result is a knot, K.
\\
\\
By construction, the knot is alternating. We draw a Heegaard diagram for 
this knot by using $m_{1}$ from the string link diagram to give a meridian 
for the knot. This meridian intersects only one $\be$. It is a basing for 
the knot in the Kauffman state picture of \cite{Alte}. We analyze the 
generators of this knot. In particular, we have seen that a generator 
${\bf x}$ of $\widehat{CF}(S, \overline{j})$ extends uniquely to a 
generator ${\bf x'}$ for $S \# B$. We would like to extend ${\bf x'}$ to a 
generator of $\widehat{CF}(K)$. However, when we forget the states on the 
other meridians, we have that the shaded regions of Figure \ref{fig:State} 
receive an assignment, but the others do not. There is a unique way to 
complete the figure to a generator ${\bf x''}$ in the knot complex, and 
this corresponds to a Kauffman state.
\\
\\

\noindent
Suppose that there is a $\phi \in \pi_{2}({\bf x}, {\bf y})$ with 
$\widehat{\mathcal{M}}(\phi) \neq \emptyset$ and with $n_{w}(\phi) = 
n_{z_{i}}(\phi) = 0$ in the Heegard diagram for $S$. 
The condition at the marked points, and that for a string link 
there is only one choice of assignment along the meridians, 
implies that the homotopy class $\phi$ gives a homotopy class 
$\phi'' \in \pi_{2}({\bf x''}, {\bf y''})$ in the diagram for $K$. That 
there is no variation on the meridians allows us to alter the generators 
as in the extension. The marked point conditions ensure that 
$\mathcal{D}(\phi)$ is wholly
insulated from the additional braid and the closure construction. 
Furthermore, 
a nighborhood of each meridian is eliminated by 
the condition that $n_{z_{i}}(\phi) = 0$ as $\phi$ must contain whole 
multiples of the meridians in its boundary.
There must then exist two intersection points ${\bf x''}$ and ${\bf y''}$ 
in 
the {\em same filtration index} for the knot (since $\phi$ doesn't cross 
the marked points) but which {\em differ in grading} by $1$. 
This contradicts the statement in 
\cite{Alte} that, for an alternating knot, the grading of every Kauffman 
state 
representing the same filtration index is the same (determined by the 
signature of the knot!).
Hence, all moduli spaces contributin to the differential must be trivial, 
$\widehat{\partial}_{S} \equiv 0$.  The same argument holds for $\phi$ as 
above 
with $\mu(\phi) \neq 0$, thus all the 
generators of the chain complex in a given filtration index in fact lie in 
the same grading.

\subsection{Euler Characteristic Calculations}
\label{sec:AlexCon}

Let $G$ and $\mathcal{F}_{i}$ be the sums if weights giving the grading
and the filtration indices. To each intersection point in $\T_{\al} \cap 
\T_{\be}$ 
assign the monomial $$(-1)^{G({\bf 
x})} h_{1}^{\mathcal{F}_{1}({\bf x})} \cdots h_{k}^{\mathcal{F}_{k}({\bf 
x})}$$ in  $\Z[h_{1}^{\pm 1}, \ldots, h_{k}^{\pm 1}]$. 
From the description of Alexander invariants for string links, the Laurent 
polynomial, $$\nabla_{S}(h_{1}, \ldots, h_{k}) = \Sigma_{\bf x} 
(-1)^{G({\bf 
x})} h_{1}^{\mathcal{F}_{1}({\bf x})} \cdots h_{k}^{\mathcal{F}_{k}({\bf 
x})}$$ is the torsion of the string link. 
Furthermore, it satisfies the Alexander-Conway 
skein relation in $h_{i}$ at self-crossings of  $L_{i}$. The proof is a 
comparison of the weights assigned to forests in the projections for the 
positive,
negative, and resolved crossings. In
fact, the three Reidemeister moves also preserve this summation, when we 
fix
the ends of the string link.  

We may ask how this sum, when restricted to a certain strand, $L_{i}$ 
relates to the Alexander-Conway polynomial of the knot, $\widehat{L}_{i}$ 
formed by closing the strand. It is shown in \cite{Alte} 
that for $1$-stranded string links (marked knots) they are identical. In 
this section we prove the following lemma, relating the above 
multi-variable
Laurent polynomial to its single variable specifications:

\begin{lemma}
The polynomial $\nabla_{S}(h_{1}, \ldots, 
h_{k})$ evaluates to $\Delta_{\widehat{L}_{i}}$, the Alexander-Conway 
polynomial of 
$\widehat{L}_{i}$, upon setting $h_{j} = 1$ for $j \neq i$.
\end{lemma}
\noindent
In \cite{Alte} the
proof of the statement for knots follows from two observations: I) from 
Kauffman, that the 
polynomial formed by the weights at crossings is the Alexander-Conway 
polynomial, and II) that the polynomial formed by using the first Chern 
class as the 
filtration index is 
symmetric due to the symmetries of $Spin^{c}$ structures on the three 
manifold found from
$0$-surgery on the knot. Since both schemes assign values to intersection 
points that
satisfy the filtration relation, and both produce symmetric polynomials 
under $h \ra h^{-1}$,
they must be equal. We do not have at hand an analog of the first 
observation, and thus will
resort to model calculations.

To prove the lemma, we first observe that

\begin{figure}
\begin{center}
\begin{picture}(2952,3671)
\put(1252,2938){\makebox(0,0)[lb]{\smash{{\SetFigFont{10}{12.0}{rm}$l_{24}$}}}}
\put(1227,1825){\makebox(0,0)[lb]{\smash{{\SetFigFont{10}{12.0}{rm}$l_{23}$}}}}
\put(1239,688){\makebox(0,0)[lb]{\smash{{\SetFigFont{10}{12.0}{rm}$l_{12}$}}}}
\includegraphics{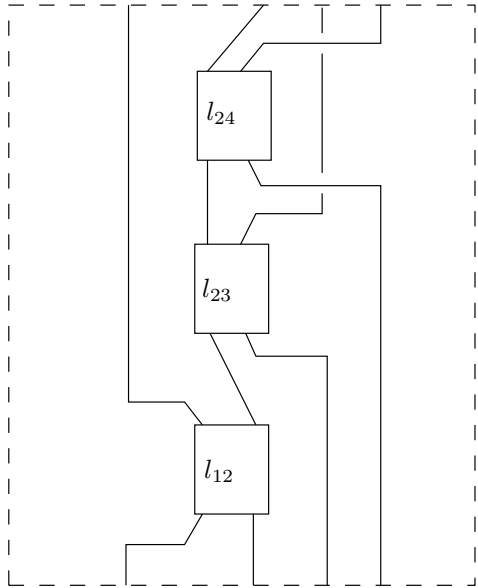}%
\end{picture}%
\end{center}
\caption{Reduced Form for a String Link when $i=2$. The numbered boxes 
indicate the 
number of full twists between the two strands. Notice that this is 
actually a braid.}\label{fig:EStd}
\end{figure}

\begin{lemma} \label{lem:Changes} Suppose we may interchange 
\ben
\item Crossings of $L_{i}$ with itself
\item Crossings of $L_{j}$ with $L_{k}$ when $j, k \neq i$
\een
Then the string link $S$ may be put in the form of a braid as found in 
Figure 
\ref{fig:EStd}. 
\end{lemma}

\begin{figure}
\begin{center}
\begin{picture}(5183,2400)
\includegraphics{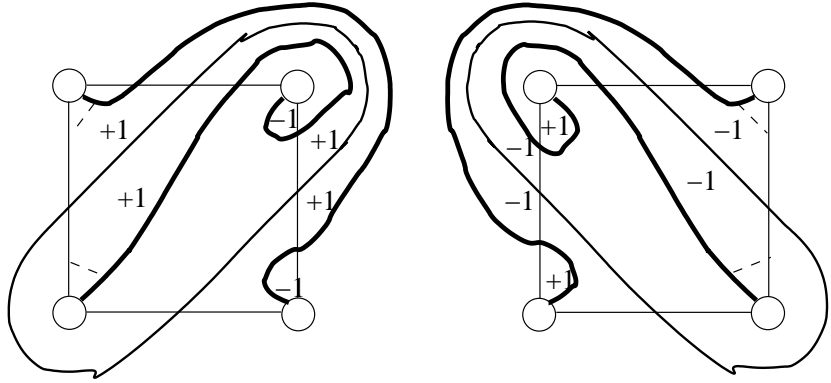}
\end{picture}%
\end{center}
\caption{Local Heegaard diagrams for self-crossings of $L_{i}$. The thick 
curves are the 
$\al$-curves. The thickest is the longitude, used as a surgery curve; 
while the middle thickness
depicts the crossing. There is a positive crossing on the left, and a 
negative crossing on 
the right. The multiplicities are for the band added at the crossing 
according to Seifert's
algorithm.}
\label{fig:cross} 
\end{figure}

{\bf Proof:}  By interchanging self-crossings of $L_{i}$ we may arrange 
for this strand to be unknotted. We may then isotope so that it is a 
vertical strand. Consider
$D^{2} \times I$ to be $I^{2} \times I$ with the $i^{th}$ strand given as 
$(\frac{1}{2}, \frac{3}{4}) \times I$. By isotoping the other strands 
vertically, switching crossings where necessary, we can ensure that each 
is contained in a narrow band $I^{2} \times (a_{j} - \epsilon, a_{j} + 
\epsilon)$ except when coming from or going to the boundary, when we 
assume that they are vertical. If we look in the band  $I^{2} \times 
(a_{j} - \epsilon, a_{j} + \epsilon)$ for the $j^{th}$ strand, we see that 
all the other strands are vertical segments and
the $j^{th}$ strand is some strand in a string link. We may isotope 
$L_{j}$ past all the vertical segments except that for $L_{i}$. By 
effecting self-crossing changes in $L_{j}$ we can make it unknotted. It 
can then be isotoped so that it reaches the level $I^{2} \times a_{j}$ 
from the bottom, winds around $L_{i}$ in that plane exactly $lk(L_{i}, 
L_{j})$ times, and then
proceeds vertically to the top. We un-spool this winding veritcally so 
that the result is a braid, where all the clasps are once again in $I^{2} 
\times (a_{j} - \epsilon, a_{j} + \epsilon)$. Doing this to each of the 
strands produces a string link with projection as
in Figure \ref{fig:EStd} when we choose $a_{j} < a_{k}$ if $j < k$ for $j, 
k \neq i$.
\\
\\
Consider the closure, $\widehat{L}_{i}$ formed by taking the $i^{th}$ 
strand and joining
the two ends in $S^{3}$ by a simple, unknotted arc. The intersection
points in the Heegaard diagram give rise to $Spin^{c}$ structures on the 
three manifold formed by taking $0$-surgery on this knot, \cite{Knot}.  
Let $\mathcal{P}_{i}$ be the periodic
domain corresponding to the Seifert surface in the Heegaard diagram for 
the $0$-surgery manifold. We already know that the sum of the filtration 
weights assigned to ${\bf x}$ by $L_{i}$ differs from 
$\frac{1}{2} <c_{1}(\mathfrak{s}({\bf x}), \mathcal{P}_{i}>$ by a 
constant, 
independent of ${\bf x}$, but which may depend upon $S$. 
Each of these quantities satisfies the 
same difference relation for a class $\phi$ joining two distinct 
intersection points. The crossings of 
$L_{j}$ with $L_{k}$
do not contribute to either calculation as the periodic domain 
$\mathcal{P}$ does not change
topology or multiplicities when we interchange such crossings. It 
consists of puctured cylinders
arising from the linking of $L_{i}$ with $L_{j}$ or $L_{k}$ and 
terminating on the meridian for
that strand. However, the punctures and multiplicities remain the same 
regardless of the type of crossing between $L_{j}$ and $L_{k}$. Since only 
the topology and multiplicities contribute to the calculation of the first 
Chern class for the
intersection points, and the intersection points correspond under crossing 
changes, the value of the first Chern class does not change. Nor do the 
weights change as such
crossings are assigned a weight of $0$. As a result, interchanging 
crossings of $L_{i}$ and $L_{j}$ does not affect the value of the 
constant.
\\
\\
It requires more effort to see that interchanging self-crossings of 
$L_{i}$ will not alter the
constant. But presuming that, we see that the polynomial 
$\nabla_{S}(h_{1}, 1, \ldots, 1)$ 
must be $h_{i}^{C_{i}(S)}\nabla_{L_{i}}$, since by \cite{Alte} the 
polynomial determined by the first 
Chern class is the Alexander-Conway polynomial of the knot. Furthermore, 
we may calculate 
$C_{i}(S)$ by finding the polynomial assigned to the reduced form above, 
since it does not change under the moves of Lemma \ref{lem:Changes}. 
The $i^{th}$ strand is then the unknot, with Alexander-Conway polynomial 
equal to $1$,
and the value of $C_{i}(S)$ is $0$ from our calculation for braids
in the previous section. 
\\
\\
That $C_{i}(S)$ does not change when interchanging self-crossings of 
$L_{i}$ follows by seeing
that the first Chern class calculation changes in the same way as the 
filtration weights. Notice that we may do the calculation for the 
multiplicities shown 
in Figure \ref{fig:cross}

By adding multiples of $[\Sigma]$ and $\mathcal{P}$ we may realize a 
periodic 
region with the 
multiplicities for the crossing as shown in the figure. We fix the 
intersection point to be 
considered. On the left
we have a local contribution to $n_{\bf x}$ of $+2$ on top, $+1$ on left 
and right, and
$0$ on bottom. We may always pick and intersection point where the 
longitude pairs with
a $\be$ intersecting the meridian, and thus does not incur an additional 
contribution from any 
crossing. On the right, the contribution is $-2$, $-1$, and $0$, 
respectively. All other
local contributions are equal as is the value of $n_{w}$ on the periodic 
region. The only
other variation that can occur will be in the Euler measure, 
$\hat{\chi}(\mathcal{P})$. Outside
of the depicted region, $\mathcal{P}$ does not change. Indeed, if we 
divide $\mathcal{P}$ 
by cutting the corners of the large $+1$ or $-1$ region along the dashed 
lines, we have a 
disc of multiplicity $\pm 1$ and the remainder,$R$, of $\mathcal{P}$, 
which is the same in 
both diagrams. Then on the left $ \hat{\chi}(\mathcal{P}) = \hat{\chi}(R) 
+ 1 - 2$ while on 
the right $ \hat{\chi}(\mathcal{P}') = \hat{\chi}(R) - 1 + 2$. The sign 
changes occur because 
of the Euler measure; in the first we add a $+1$ disc joined along two 
segments 
to a $+1$ region and
at two points to a $-1$ region, which contributes nothing. In the second 
case the multiplicities are reversed, and the Euler measure must be 
calculated differently.
Taking the difference of these, and adding the difference of the 
contribution from each quadrant
gives $+2 - (-2) + (-1 - 1) = 2$ on top, $+1 - (-1) + (-1 - 1) = 0$ on 
left and right, and 
$0 - 0 + (-1 -1) = -2$ on bottom. This is $-2$ times the difference in the 
weights for
these quadrants, but that is precisely the factor we divide into the first 
Chern class to get
a filtration index. Hence, the weights for any intersection point
before and after a crossing change 
differ from the 
first Chern class calculation for the corresponding intersection point by 
the same amount. In
particular, $C_{i}(S)$ does not change.
\\
\\
Implicit in this discussion is the following corollary:

\begin{cor} For the filtration indices, $\mathcal{F}_{i}$, calculated
from the weights
$$
<c_{1}(\mathfrak{s}({\bf x}), \mathcal{P}_{i}> = 2\,\mathcal{F}_{i}({\bf 
x}) 
$$
\end{cor}

\subsection{Tangles}

Recall the definition of a  tangle:

\begin{defn}
A {\em tangle}, $\tau$,  in $D^{2} \times I$ is an oriented smooth 
$1$-sub-manifold whose boundary lies in $D^{2} \times \{0, 1\}$. Two 
tangles are isotopic if there exists a self-diffeomorphism of $D^{2} 
\times I$, isotopic to the identity, fixing the ends $D^{2} \times \{0, 
1\}$, and carrying one tangle into the other while preserving the 
orientations of the components. 
\end{defn}

Suppose we have a tangle with $n$-components. We will call it $m$-colored 
if there is a function, $c$, from the components of the tangle, 
$\tau_{j}$, to $\{1, \ldots, m\}$. We require our isotopies to preserve 
the value of $c$. We restrict to those tangles for which exactly one 
component 
in $c^{-1}(i)$ is open for each $i \in \{1, \ldots, m\}$, and this 
component
is oreiented from $D^{2} \times \{1\}$ to $D^{2} \times \{0\}$. The 
collection
of open components will then form a string link.

To construct a Heegaard diagram we convert the tangle, $\tau$, to an 
associated string link, $S(\tau)$ in another three manifold. First, 
connect the 
components with the same color by  paths in the complement of the tangle 
so that, with the paths as edges and the components as vertices, we have a 
tree 
rooted at the open component for that color. These paths may be used to 
band sum the components together, using any number of half twists in the 
band which will match the orientations on components correctly. 
We then perform $0$-surgery on an unknot linking the band once. 
The resulting manifold, $Y$, is a connect 
sum of $s$ copies of  $S^{1} \times S^{2}$'s where $s$ is the number of 
closed components. 
The image of  $\tau$ after perfoming the band sums is a multi-component 
``string link'', $S(\tau)$, in the 
complement of a ball. We have a color map $c$ on this string link which we 
use to index the components. The homology group for index $\overline{j}$ 
for
the colored tangle, $\tau$,  is defined to be $\widehat{HF}( Y, S(\tau); 
\mathfrak{s}_{0}, \overline{j})$ where $\mathfrak{s}_{0}$ is the torsion 
$Spin^{c}$
structure. Underlying this construction is the following lemma:

\begin{lemma}
The isotopy class of $S(\tau)$ is determined by that of $\tau$
\end{lemma}
\noindent
{\bf Proof:} This follows as for Proposition 2.1 of \cite{Knot}.
 It suffices to show that the choices made in performing the 
band sums do not affect the isotopy class of $S(\tau)$. 
Once we add the $0$-framed handles, the choice of the bands no longer 
matters. In each band, we may remove full twists by using the belt trick 
to replace them with a self-crossing of the band. We may then isotope the 
$0$-framed circle to the self-crossing and slide {\em one} of the strands 
in the band across the handle.  Done appropriately this will undo a full 
twist in the band.  Furthermore, by sliding across the $0$-framed handles, 
we may move the bands past each other or any component in the string link. 
By using the trick from \cite{Knot} illlustrated in Figure 
\ref{fig:tangle}, we may arrange that all the closed components of the 
same color are linked 
in a chain to a single open component. In addition, we may interchange any 
two components along the chain. The combination of these moves provides 
Heegaard equivalences, not involving the meridians, between any two ways 
of joining the closed components in a color to that with boundary, 
regardless of the paths for the band sums, or twists in the bands.\\

\setlength{\unitlength}{3000sp}%
\begin{figure}
\begin{center}
\begin{picture}(7385,6341)(1343,-6028)
\thinlines
\put(5968,-644){\line( 1, 0){470}}
\put(6438,-644){\line( 0, 1){329}}
\put(3843,-1844){\vector( 1, 0){987}}
\put(6918,-2169){\vector( 0,-1){846}}
\put(3984,-4019){\vector( 1, 0){846}}
\put(1480,-5744){\line( 0, 1){705}}
\put(1480,-5039){\line( 1, 0){1175}}
\put(2655,-5039){\line( 0,-1){705}}
\put(3792,-3406){\line(-1, 0){752}}
\put(3040,-3406){\line( 0,-1){1316}}
\put(3040,-4722){\line( 1, 0){752}}
\put(1355,-3431){\line( 0,-1){752}}
\put(1355,-4183){\line( 1, 0){1128}}
\put(2483,-4183){\line( 0, 1){799}}
\put(2005,-4719){\line( 0,-1){329}}
\put(2005,-5048){\line( 0, 1){  0}}
\put(2642,-3669){\line( 0, 1){141}}
\put(2642,-3528){\line( 1, 0){141}}
\put(2783,-3528){\line( 0,-1){470}}
\put(2783,-3998){\line(-1, 0){141}}
\put(2642,-3998){\line( 0, 1){188}}
\put(2480,-3744){\line( 1, 0){235}}
\put(2817,-3744){\line( 1, 0){235}}
\put(1992,-4557){\line( 0, 1){235}}
\put(1992,-4322){\line( 1, 0){1034}}
\put(1905,-4409){\line(-1, 0){282}}
\put(1623,-4409){\line( 0,-1){235}}
\put(1623,-4644){\line( 1, 0){705}}
\put(2328,-4644){\line( 0, 1){235}}
\put(2328,-4409){\line(-1, 0){282}}
\put(6142,-4404){\line( 0,-1){376}}
\put(6142,-4780){\line( 0, 1){  0}}
\put(6142,-4942){\line( 0,-1){329}}
\put(6142,-5271){\line( 0, 1){  0}}
\put(5617,-5967){\line( 0, 1){705}}
\put(5617,-5262){\line( 1, 0){376}}
\put(6142,-5254){\line( 1, 0){470}}
\put(6612,-5254){\line( 0,-1){752}}
\put(5979,-4942){\line( 0,-1){329}}
\put(5979,-5271){\line( 0, 1){  0}}
\put(5867,-4550){\line(-1, 0){282}}
\put(5585,-4550){\line( 0,-1){329}}
\put(5585,-4879){\line( 1, 0){893}}
\put(6478,-4879){\line( 0, 1){329}}
\put(6478,-4550){\line(-1, 0){282}}
\put(6196,-4550){\line( 0, 1){  0}}
\put(5967,-4392){\line(-1, 0){470}}
\put(5497,-4392){\line( 0, 1){752}}
\put(6142,-4392){\line( 1, 0){470}}
\put(6612,-4392){\line( 0, 1){329}}
\put(6617,-3842){\line( 0, 1){235}}
\put(5967,-4392){\line( 0,-1){376}}
\put(5967,-4768){\line( 0, 1){  0}}
\put(6617,-3842){\line( 1, 0){235}}
\put(6617,-4067){\line( 1, 0){235}}
\put(6954,-4067){\line( 1, 0){235}}
\put(6955,-3842){\line( 1, 0){235}}
\put(7179,-4067){\line( 0,-1){846}}
\put(7179,-4913){\line( 1, 0){752}}
\put(7179,-3829){\line( 0, 1){235}}
\put(7179,-3594){\line( 1, 0){752}}
\put(6779,-3867){\line( 0, 1){  0}}
\put(6729,-3844){\line( 1, 0){282}}
\put(6729,-4069){\line( 1, 0){282}}
\put(5480,-3269){\line( 0,-1){235}}
\put(5305,-3393){\line( 0,-1){188}}
\put(5305,-3581){\line( 1, 0){1410}}
\put(6715,-3581){\line( 0, 1){188}}
\put(6715,-3393){\line(-1, 0){ 94}}
\put(5317,-3394){\line( 1, 0){141}}
\put(5593,-3394){\line( 1, 0){940}}
\put(5968,-656){\line( 0,-1){376}}
\put(5968,-1032){\line( 0, 1){  0}}
\put(5968,-1194){\line( 0,-1){329}}
\put(5968,-1523){\line( 0, 1){  0}}
\put(5443,-2219){\line( 0, 1){705}}
\put(5443,-1514){\line( 1, 0){376}}
\put(5968,-1506){\line( 1, 0){470}}
\put(6438,-1506){\line( 0,-1){752}}
\put(5805,-1194){\line( 0,-1){329}}
\put(5805,-1523){\line( 0, 1){  0}}
\put(5693,-802){\line(-1, 0){282}}
\put(5411,-802){\line( 0,-1){329}}
\put(5411,-1131){\line( 1, 0){893}}
\put(6304,-1131){\line( 0, 1){329}}
\put(6304,-802){\line(-1, 0){282}}
\put(6022,-802){\line( 0, 1){  0}}
\put(5793,-644){\line(-1, 0){470}}
\put(5323,-644){\line( 0, 1){752}}
\put(6593,-3189){\line( 0,-1){329}}
\put(6443,-94){\line( 0, 1){235}}
\put(5793,-644){\line( 0,-1){376}}
\put(5793,-1020){\line( 0, 1){  0}}
\put(6443,-94){\line( 1, 0){235}}
\put(6443,-319){\line( 1, 0){235}}
\put(6780,-319){\line( 1, 0){235}}
\put(6781,-94){\line( 1, 0){235}}
\put(7005,-319){\line( 0,-1){846}}
\put(7005,-1165){\line( 1, 0){752}}
\put(7005,-81){\line( 0, 1){235}}
\put(7005,154){\line( 1, 0){752}}
\put(6593,-36){\line( 0, 1){ 94}}
\put(6593, 58){\line( 1, 0){141}}
\put(6734, 58){\line( 0,-1){564}}
\put(6605,-119){\line( 0, 1){  0}}
\put(6605,-356){\line( 0,-1){141}}
\put(6605,-497){\line( 1, 0){141}}
\put(1493,-2219){\line( 0, 1){705}}
\put(1493,-1514){\line( 1, 0){1175}}
\put(2668,-1514){\line( 0,-1){705}}
\put(3805,119){\line(-1, 0){752}}
\put(3053,119){\line( 0,-1){1316}}
\put(3053,-1197){\line( 1, 0){752}}
\put(1368, 94){\line( 0,-1){752}}
\put(1368,-658){\line( 1, 0){1128}}
\put(2496,-658){\line( 0, 1){799}}
\put(2018,-656){\line( 0,-1){376}}
\put(2018,-1032){\line( 0, 1){  0}}
\put(2018,-1194){\line( 0,-1){329}}
\put(2018,-1523){\line( 0, 1){  0}}
\put(1930,-856){\line(-1, 0){282}}
\put(1648,-856){\line( 0,-1){235}}
\put(1648,-1091){\line( 1, 0){705}}
\put(2353,-1091){\line( 0, 1){235}}
\put(2353,-856){\line(-1, 0){282}}
\put(2655,-144){\line( 0, 1){141}}
\put(2655, -3){\line( 1, 0){141}}
\put(2796, -3){\line( 0,-1){470}}
\put(2796,-473){\line(-1, 0){141}}
\put(2655,-473){\line( 0, 1){188}}
\put(2493,-219){\line( 1, 0){235}}
\put(2830,-219){\line( 1, 0){235}}
\put(4043,-2081){\makebox(0,0)[lb]{\smash{{\SetFigFont{12}{14.4}{rm}band 
sum}}}}
\put(7218,-2606){\makebox(0,0)[lb]{\smash{{\SetFigFont{12}{14.4}{rm}handleslide}}}}
\put(4080,-4231){\makebox(0,0)[lb]{\smash{{\SetFigFont{12}{14.4}{rm}band 
sum}}}}
\put(2642,-3481){\makebox(0,0)[lb]{\smash{{\SetFigFont{12}{14.4}{rm}$0$}}}}
\put(2392,-4569){\makebox(0,0)[lb]{\smash{{\SetFigFont{12}{14.4}{rm}$0$}}}}
\put(2367,-5956){\makebox(0,0)[lb]{\smash{{\SetFigFont{12}{14.4}{rm}$L_{i}$}}}}
\put(3942,-4794){\makebox(0,0)[lb]{\smash{{\SetFigFont{12}{14.4}{rm}$L_{j}$}}}}
\put(1430,-3794){\makebox(0,0)[lb]{\smash{{\SetFigFont{12}{14.4}{rm}$L_{k}$}}}}
\put(6529,-4792){\makebox(0,0)[lb]{\smash{{\SetFigFont{12}{14.4}{rm}$0$}}}}
\put(8079,-5017){\makebox(0,0)[lb]{\smash{{\SetFigFont{12}{14.4}{rm}$L_{j}$}}}}
\put(5567,-4017){\makebox(0,0)[lb]{\smash{{\SetFigFont{12}{14.4}{rm}$L_{k}$}}}}
\put(5942,-5954){\makebox(0,0)[lb]{\smash{{\SetFigFont{12}{14.4}{rm}$L_{i}$}}}}
\put(6880,-3569){\makebox(0,0)[lb]{\smash{{\SetFigFont{12}{14.4}{rm}$0$}}}}
\put(6355,-1044){\makebox(0,0)[lb]{\smash{{\SetFigFont{12}{14.4}{rm}$0$}}}}
\put(7905,-1269){\makebox(0,0)[lb]{\smash{{\SetFigFont{12}{14.4}{rm}$L_{j}$}}}}
\put(5393,-269){\makebox(0,0)[lb]{\smash{{\SetFigFont{12}{14.4}{rm}$L_{k}$}}}}
\put(6543,119){\makebox(0,0)[lb]{\smash{{\SetFigFont{12}{14.4}{rm}$0$}}}}
\put(5768,-2206){\makebox(0,0)[lb]{\smash{{\SetFigFont{12}{14.4}{rm}$L_{i}$}}}}
\put(2655, 44){\makebox(0,0)[lb]{\smash{{\SetFigFont{12}{14.4}{rm}$0$}}}}
\put(2405,-1044){\makebox(0,0)[lb]{\smash{{\SetFigFont{12}{14.4}{rm}$0$}}}}
\put(2380,-2431){\makebox(0,0)[lb]{\smash{{\SetFigFont{12}{14.4}{rm}$L_{i}$}}}}
\put(3955,-1269){\makebox(0,0)[lb]{\smash{{\SetFigFont{12}{14.4}{rm}$L_{j}$}}}}
\put(1443,-269){\makebox(0,0)[lb]{\smash{{\SetFigFont{12}{14.4}{rm}$L_{k}$}}}}
\end{picture}%
\end{center}
\caption{Kirby Calculus Interchange of  Connecting Paths for a Tangle}
\label{fig:tangle}
\end{figure}
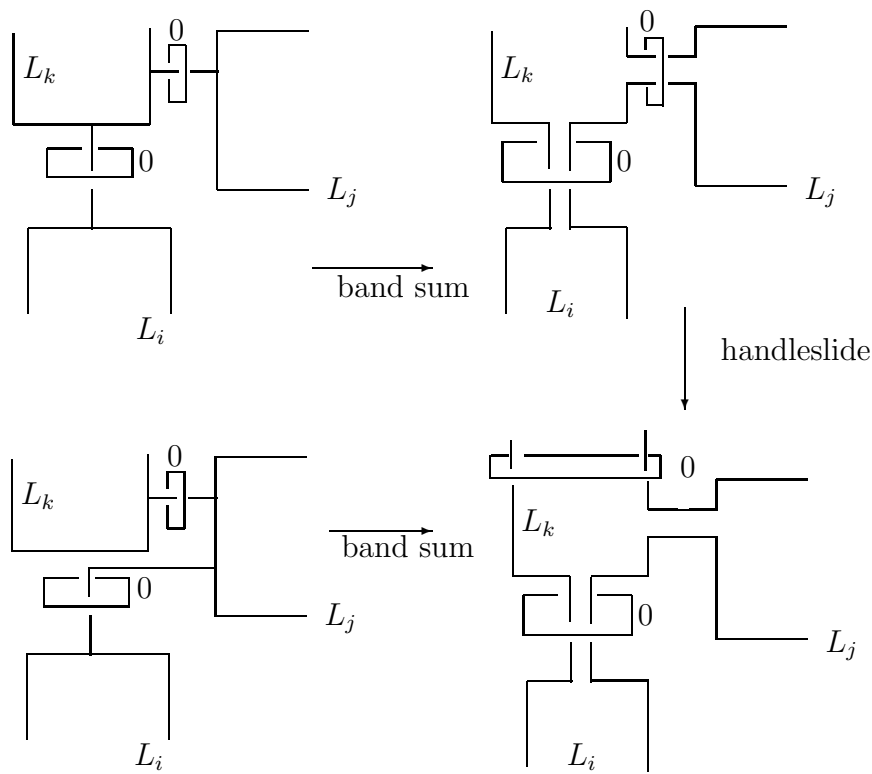

\subsubsection{Long Exact Sequences}
\ \\
The tangle formulation where there is one open component for
each color permits the introduction of the skein long exact 
sequence found in 
\cite{Knot} where we may resolve crossings of components with the same 
color, but
not crossings involving different colors. We should think of the resolved 
crossing, arising from
$0$-surgery on an unknot in the long exact surgery sequence, in the 
context of tangles. We let $L_{-}^{c}$ be a tangle with a negative 
self-crossing in color $c$, $L_{+}^{c}$ be the tangle with a positive 
self-crossing instead, and $L_{0}^{c}$ be the tangle resulting from 
resolving the crossing. As in \cite{Knot} there are two sequences. If the 
crossing is  a self-crossing of a component then

$$
\ra \widehat{HF}(L_{-}^{c}, \overline{j}) \ra \widehat{HF}(L_{0}^{c}, 
\overline{j}) \ra \widehat{HF}(L_{+}^{c}, \overline{j}) \ra
$$
whereas if the crossing occurs between different components of the same 
color we have
$$
\ra \widehat{HF}(L_{-}^{c}) \ra \widehat{HF}(L_{0}^{c}) \otimes V\ra 
\widehat{HF}(L_{+}^{c}) \ra
$$
where $V = V_{-1} \oplus V_{0} \oplus V_{+1}$ and $V_{-1}$ consists of a 
$\Z$ in filtration index $-1$ for the color $c$ and $0$ for all others, 
$V_{0}$ consists of $\Z^{2}$ with filtration index $0$ for all colors, and 
$V_{+1}$ consists of  a $\Z$ in color $c$ filtration index $+1$. The maps 
preserve the filtration indices with the tensor product index defined as 
the sum of those on the two factors. The proof is identical to that in 
\cite{Knot}.

{\bf Note:} Since the theory for tangles arises from thinking of 
them as string links in an another manifold, the result for composition of
string links extends to composition of this sub-class of tangles. The 
sub-class condition disallows the formation of new closed components when
composing.

\end{document}